\documentclass[onside,onecolumn]{sn-jnl}
\usepackage[margin=1in]{geometry}
\usepackage[todo,check]{style}
\usepackage{overpic}
\usepackage{algorithm,algorithmic}

\title{Efficient, Nonlinear Second Moment Methods for Multigroup Thermal Radiative Transfer}
\author*[1]{\fnm{Samuel} \sur{Olivier}}\email{solivier@lanl.gov}
\author[1]{\fnm{James S.} \sur{Warsa}}
\author[2]{\fnm{HyeongKae} \sur{Park}}
\affil[1]{\orgdiv{Computational Physics and Methods}, \orgname{Los Alamos National Laboratory}}
\affil[2]{\orgdiv{Fluid Dynamics and Solid Mechanics}, \orgname{Los Alamos National Laboratory}}

\begin{document}
\abstract{
Thermal radiative transfer (TRT) presents significant computational challenges due to the stiff, nonlinear coupling between radiation and material energy, particularly in multigroup, high-fidelity transport models. 
In this work, we develop an efficient nonlinear acceleration framework for TRT based on the Second Moment (SM) method. 
Our approach couples high-order discrete ordinates transport to a gray, diffusion-based low-order system that implicitly resolves the stiff absorption–emission physics, isolating this stiffness from the high-order system.
The resulting algorithm alternates between transport sweeps and a Newton-type solution of the coupled low-order and material energy balance equations, utilizing nonlinear temperature elimination for improved robustness. 
Crucially, our approach is the first moment-based TRT algorithm with a symmetric and positive definite (SPD) low-order system enabling scalable linear solves via algebraic multigrid-preconditioned conjugate gradient. 
We investigate both consistent and independent low-order discretizations within a discontinuous Galerkin framework and assess their performance on one and two-dimensional gray and multigroup benchmark problems. 
A discrete reference approach is used to assess numerical error in space-time convergence studies of challenging TRT problems. 
Results demonstrate that these algorithms achieve robust nonlinear convergence and significant reductions in transport iterations compared to unaccelerated schemes, resulting in large speedups in overall runtime. 
While the independent formulation offers improved iteration counts on under-resolved meshes, the consistent method provides superior solution quality and robustness. 
Overall, this work establishes the SM method as an effective and scalable approach for nonlinear multigroup TRT and provides insight into the interplay between discretization consistency, accuracy, and solver performance arising in moment-based acceleration algorithms. 
}


\maketitle

\section{Introduction}
Thermal radiative transfer (TRT) models the exchange of energy between thermal radiation and the interacting medium through which it propagates.
TRT is a central component of high energy density physics, inertial confinement fusion, and astrophysics, where extreme temperatures make radiative heat transfer the dominant mechanism of energy exchange. In these regimes, radiation transport is strongly coupled to the evolution of the material temperature through absorption and emission processes that depend nonlinearly on temperature.
This nonlinear coupling introduces significant computational challenges. 
In particular, absorption–emission interactions are stiff on the macroscopic time scales of interest, making explicit time integration impractical. 
As a result, implicit schemes are required.
However, straightforward iterative treatments of the coupled system often exhibit slow convergence and limited robustness requiring time step control to maintain efficiency. 
These challenges are further exacerbated in the frequency-dependent case, where spectral dependence introduces richer radiation-material interaction physics that increase stiffness on top of the increase in dimensionality. 
Finally, high-fidelity models use a Boltzmann description of radiation transport, resulting in a high-dimensional, stiff nonlinear system that is expensive to solve. 

To address these issues, a broad class of moment-based acceleration techniques has been developed including the Variable Eddington Factor (VEF) \cite{mihalas1999foundations}, Quasidiffusion (QD) \cite{goldin}, High-Order Low-Order (HOLO) \cite{CHACON201721}, and Second Moment (SM) \cite{lewis_miller} methods. 
These approaches introduce a low-order system that captures the dominant behavior of the slow-to-converge physics, such as absorption-emission in TRT \cite{anistratov1996nonlinear} and fission in neutronics \cite{Hardy03042025}, in a reduced-dimensional model where it is computationally tractable to resolve such physics implicitly. 
The low-order problem includes discrete closures computed from the high-order system, allowing the low-order system to reproduce the physics fidelity of the high-order system. 
High and low-order solves are alternated with the high-order problem providing approximate closures to the low-order problem and the low-order problem fully resolving the expensive physics. 
The closures are designed to weakly depend on the high-order system; thus, given an efficient low-order solver for the stiff physics, this iterative approach converges rapidly and robustly independent of numerical resolution parameters and material properties. 
Notably, SM uses \emph{additive} closures making the left-hand-side of the low-order system equivalent to radiation diffusion which can thus be discretized in a symmetric and positive definite manner amenable to solution via preconditioned conjugate gradient. 
In contrast,  VEF, QD, and HOLO all have advective-like closure terms that destroy the symmetry of the low-order system such that more expensive non-symmetric-compatible Krylov solvers such as GMRES or BiCGStab are required.
Recent work has expanded this paradigm as applied to TRT through multilevel-in-frequency acceleration \cite{anistratov2020nonlinear,ANISTRATOV2019186,GHASSEMI2020109315}, reduced memory schemes \cite{COALE2023108458}, high-order finite element discretizations \cite{yee_mc21}, high-order implicit-explicit time integration \cite{SOUTHWORTH2024113349}, and tightly coupled radiation-hydrodynamics algorithms \cite{LOU2021110393,Wollaber02012017}.

In this work, we develop a nonlinear acceleration framework that employs the \emph{gray} SM low-order system to accelerate multigroup TRT. 
This framework is centered around the discontinuous Galerkin (DG) spatial and semi-implicit backward Euler temporal discretizations of the multigroup, discrete ordinates (\Sn) transport equations for the high-order problem. 
While SM has previously been studied in the context of steady-state, linear transport \cite{csmm,smm_olivier}, its extension to time-dependent, nonlinear TRT introduces several nontrivial challenges, including the treatment of the nonlinear absorption–emission coupling, multigroup effects, and solver robustness.
We elect to use a gray low-order system to avoid the large cost associated with solving multigroup diffusion-type problems especially since linearization of the Planck emission term induces dense group-to-group coupling that is difficult to resolve efficiently. 
While multilevel-in-frequency approaches reduce this cost \cite{anistratov2020nonlinear}, gray preconditioners, such as gray transport acceleration \cite{LARSEN1988459} and linear multi-frequency gray \cite{MOREL2007244}, are known to be remarkably effective at reducing transport iterations in classical schemes suggesting that a gray moment-based method could retain acceleration efficacy at minimal cost. 
The resulting algorithm combines transport sweeps with a Newton-type solve of the coupled low-order and energy balance equations and utilizes nonlinear elimination of the temperature to improve robustness. 
Since we use a gray low-order system, it is practical to assemble and solve the Jacobian of the Newton-type iteration, resulting in an efficient nonlinear low-order solve.
In multigroup problems, the SPD property of the SM low-order system is preserved by using a novel, \emph{additive} multigroup collapse closure, resulting in the first moment-based TRT algorithm that resolves absorption-emission implicitly with only SPD low-order diffusion solves. 
In contrast, Park et al. \cite{doi:10.1080/00411450.2012.671224} and Anistratov \cite{anistratov2020nonlinear} use a \emph{multiplicative} multigroup closure which introduces an additional non-symmetric, advective-type term to the low-order system. 

Within this framework, we consider two classes of low-order discretizations, consistent and independent, which differ in whether the low-order system is constructed to be algebraically equivalent to the moments of the high-order discretization. 
We compare the performance of the TRT algorithms arising from extending the Local Discontinuous Galerkin (LDG)-based consistent \cite{csmm} and independent \cite{smm_olivier} SM methods from steady-state, linear transport to time-dependent, multigroup, nonlinear TRT.
While many independent and consistent low-order discretizations compatible with our choice of high-order discretization exist \cite{me,dgvef_olivier,rtvef_olivier,dima_dfem,two-level-independent-warsa}, the LDG-based SM methods approximate the energy density and radiative flux (zeroth and first moments) in the same DG finite element space as the high-order system and such methods were shown to be scalably solved by algebraic multigrid-preconditioned conjugate gradient. 
In addition, use of LDG for both the independent and consistent methods minimizes the differences between the consistent and independent approaches, enabling a more direct comparison of their inherent properties.
Note that, while these approaches have been studied for linear transport, their behavior in nonlinear, time-dependent, multigroup TRT has not been previously investigated.
In particular, on steady-state, linear transport problems, the independent method was observed to converge rapidly and robustly but produce spatially diffusive errors and a poor approximation of the first moment (i.e.~radiative flux) while the consistent method produced high-quality solutions but converged slowly on spatially under-resolved problems. 
Here, we investigate the independent method's solution quality and robustness and the consistent method's iterative convergence rate where the nonlinearity, time dependence, and often severe under-resolution in space associated with TRT problems will exacerbate these observed behaviors. 
Note that comparison of the consistent and independent formulations serves not merely as a numerical study, but as a means to elucidate the interplay between discretization consistency, nonlinear convergence, solver conditioning, and solution robustness on challenging TRT problems.
Such insights will inform the design of SM methods for other challenging physics problems, such as extensions to radiation-hydrodynamics, and, due to the close connection between these methods, will apply to VEF, QD, and HOLO algorithms as well. 

Finally, we conduct space-time convergence studies with respect to a spatially and temporally resolved \emph{discrete} reference solution to ascertain the accuracy of these methods. 
This discrete reference approach allows investigation of numerical error on physically realistic problems, avoiding the physics simplifications required to compare to analytic or semi-analytic solutions.
This enables comparison of the efficiency, in accuracy per cost, of the consistent and independent approaches on realistic benchmark problems. 
The discrete reference is used to evaluate the accuracy of certain algorithmic choices relating to the extension of the SM method to time and frequency-dependent problems. 
Computational performance is compared to an unaccelerated scheme that does not use a low-order problem to accelerate absorption-emission to characterize the performance of the methods with respect to a common control method and to ensure that the acceleration provided by the low-order system sufficiently outweighs the additional cost of solving the low-order system. 

To our knowledge, this work represents the first:
\begin{enumerate}
	\item application of the SM method to multigroup, time-dependent TRT, 
	\item moment-based acceleration technique for TRT with an SPD low-order system,
	\item systematic investigation of the consistent and independent approaches in the nonlinear TRT setting. 
\end{enumerate}
The paper proceeds as follows. 
We begin by deriving the multigroup, \Sn TRT system with semi-implicit backward Euler time discretization. 
The corresponding semi-discrete low-order problem is introduced, detailing the SM closures and the collapse to gray. 
The high-level details of the solution procedure are outlined. 
In \S\ref{sec:transport_disc}, the transport and energy balance equations are discretized in space with DG.
The low-order systems are then extended to the TRT system in \S\ref{sec:smm}.
We derive an efficient Newton-type iteration for solving the coupled low-order and energy balance equations in \S\ref{sec:nonlinear_low_order_solve}. 
In addition, we show that the resulting linearized low-order problem is symmetric and positive definite, allowing the use of algebraic multigrid-preconditioned conjugate gradient to solve the linearized diffusion problem at each Newton iteration. 
We then present numerical results comparing accuracy and performance of the consistent and independent methods on a gray, one-dimensional Marshak wave \cite{GANAPOL1983311}, Larsen's problem \cite{LARSEN1988459}, the gray crooked pipe benchmark \cite{2001JCoPh.172..543G}, and a multigroup lattice problem \cite{osti_2280904}. 
Finally, \S\ref{sec:conclusions} gives conclusions and recommendations for future work.

\section{Thermal Radiative Transfer Algorithm} \label{sec:trt}
We consider the thermal radiative transfer system where physical scattering is neglected. 
Frequency-dependence is modeled with the multigroup approximation where the frequency domain is subdivided into $G$ conservation equations with group boundaries defined by the sequence $\{\nu_1,\ldots, \nu_{G+1}\}$ where $\nu_g < \nu_{g+1}$.
The direction variable is modeled with Discrete Ordinates (\Sn) where the solution is collocated at a set of discrete directions taken from a quadrature rule for the unit sphere, $\{\Omegahat_d,w_d\}_{d=1}^{N_\Omega}$.
This multigroup, \Sn TRT system is given by: 
\begin{subequations}
\begin{equation} \label{eq:transport}
	\frac{1}{c} \pderiv{I_{d,g}}{t} + \Omegahat_d\cdot\nabla I_{d,g} + \sigma_g I_{d,g} = \frac{1}{4\pi}\sigma_g B_g(T) \,, \quad \x \in \D\,,\ 1 \leq g \leq G \,,\ 1\leq d \leq N_\Omega \,,
\end{equation}
\begin{equation}
	I_{d,g}(\x,t) = \frac{1}{4\pi} B_g(T_\text{bdr}) \,, \quad \x \in \partial \D \ \mathrm{and} \ \Omegahat_d\cdot\n < 0 \,,\ 1 \leq g \leq G\,,\ 1\leq d \leq N_\Omega \,,  
\end{equation}
\begin{equation} \label{eq:energy_balance}
	C_v \pderiv{T}{t} = \sum_{g=1}^G \paren{c\sigma_g E_g - \sigma_g B_g(T)} \,, \quad \x \in \D \,,
\end{equation}
\end{subequations}
where $\D \subset \R^{\dim}$ denotes the $\dim$-dimensional spatial domain with $\partial\D$ its boundary, $\x \in \R^{\dim}$, $\Omegahat \in \mathbb{S}^2$, and $t \in \R^+$ the spatial, direction-of-flight, and temporal variables, respectively,
$$I_{d,g}(\x,t) = \int_{\nu_g}^{\nu_{g+1}} I(\x,\Omegahat_d,\nu,t) \ud \nu$$
the intensity in group $g$, $T(\x,t)$ the material temperature, $\sigma_g(\x,T)$ the absorption opacity in group $g$, $c$ the speed of light, $C_v$ the heat capacity, and $T_\text{bdr}(\x)$ the prescribed boundary temperature. 
Here, 
	\begin{equation}
		B_g(T) = \int_{\nu_g}^{\nu_{g+1}} \frac{2k^4}{c^2 h^3} \frac{\nu^3}{e^{\nu/T} - 1} \ud \nu 
	\end{equation}
is the Planck emission function integrated over group $g$ with $k$ the Boltzmann constant and $h$ the Planck constant.
Note that we use the convention that $h\nu \rightarrow \nu$ and $kT \rightarrow T$ such that $\nu$ and $T$ both have units of energy in electronvolts. 
Following Clark \cite{CLARK1987311}, the Planck emission is computed as 
	\begin{subequations}
	\begin{equation} \label{eq:planck}
		B_g(T) = B(T)\,b_g(T) \,,
	\end{equation}
	\begin{equation} \label{eq:black_body}
		B(T) = ac T^4 \,, 
	\end{equation}
	\end{subequations}
where $a = \SI{137}{\erg\per\cm\cubed\per\eV^4}$ is the radiation constant and 
	\begin{equation}
		b_g(T) = \int_{\nu_{g}/T}^{\nu_{g+1}/T} \frac{15}{\pi^4}\frac{x^3}{e^x - 1} \ud x \,,
	\end{equation}
the normalized Planck spectrum such that $\sum_{g=1}^G b_g(T) = 1$. 
The normalized spectrum in each group is computed using polylogarithmic and Taylor series expansions of the frequency-dependent Planck emission function \cite{CLARK1987311}. 
In addition, we define the energy density in group $g$ as 
	\begin{equation}
		E_g(\x,t) = \frac{1}{c} \sum_{d=1}^{N_\Omega} w_d\,I_{d,g}(\x,t) \,. 
	\end{equation}
The SM algorithm extends this system by including a low-order, transport-corrected diffusion model. 
The low-order system's energy density solution is used in the absorption term in the energy balance equation (Eq.~\ref{eq:energy_balance}), isolating the stiff absorption-emission physics from the high-dimensional transport equation. 
High-fidelity transport solutions are obtained through the use of the SM closures which linearly couple the high and low-order systems. 
In this paper, we use a gray low-order problem to accelerate the multigroup TRT system. 
This further reduces the dimensionality of the low-order system, making computing and inverting the Jacobian of the low-order system associated with a Newton-type iteration practical. 
The SM TRT system is: 
\begin{subequations}
\begin{equation} \label{eq:smm_transport}
	\frac{1}{c} \pderiv{I_{d,g}}{t} + \Omegahat\cdot\nabla I_{d,g} + \sigma_g I_{d,g} = \frac{1}{4\pi}\sigma_g B_g(T) \,, \quad \x \in \D\,,\ 1 \leq g \leq G \,,\ 1\leq d \leq N_\Omega \,, 
\end{equation}
\begin{equation}
	I_{d,g}(\x,t) = \frac{1}{4\pi} B_g(T_\text{bdr}) \,, \quad \x \in \partial \D \ \mathrm{and} \ \Omegahat\cdot\n < 0 \,,\ 1 \leq g \leq G\,,\ 1\leq d\leq N_\Omega \,, 
\end{equation}
\begin{equation} \label{eq:zeroth}
	\pderiv{E}{t} + \nabla\cdot\vec{F} + c \sigma_E E = \sigma_P B(T) \,, \quad \x \in \D \,,
\end{equation}
\begin{equation} \label{eq:first}
	\frac{1}{c}\pderiv{F}{t} + \frac{c}{3}\nabla E + \sigma_F \vec{F} = \sum_g (\sigma_F - \sigma_g) \vec{F}_g - \nabla\cdot\T \,, \quad \x \in \D\,,
\end{equation}
\begin{equation} \label{eq:smm_bc_analytic}
	\vec{F}\cdot\n = c \alpha E + \beta + 2 F_\text{in} \,, \quad \x \in \partial \D \,, 
\end{equation}
\begin{equation} \label{eq:smm_meb}
	C_v \pderiv{T}{t} = c\sigma_E E - \sigma_P B(T) \,, \quad \x \in \D \,,
\end{equation}
\end{subequations}
where the gray energy density, $E$, and flux, $\vec{F}$, are 
	\begin{equation}
		E(\x,t) = \frac{1}{c}\sum_{g=1}^G\sum_{d=1}^{N_\Omega} w_d\,I_{d,g}(\x,t) \,,
	\end{equation}
	\begin{equation}
		\vec{F}(\x,t) = \sum_{g=1}^G \sum_{d=1}^{N_\Omega}\Omegahat_d\, w_d\, I_{d,g}(\x,t) \,. 
	\end{equation}
The \emph{gray} SM volumetric and boundary corrections are: 
\begin{subequations}
\begin{equation}
	\mat{T}(I) = \sum_{g=1}^G \sum_{d=1}^{N_\Omega} w_d\paren{\Omegahat_d\otimes\Omegahat_d - \frac{1}{3}\I} I_{d,g} \,, 
\end{equation}
\begin{equation}
	\beta(I) = \sum_{g=1}^G \sum_{d=1}^{N_\Omega} w_d\paren{|\Omegahat_d\cdot\n| - \alpha} I_{d,g} \,, 
\end{equation}
\end{subequations}
respectively, where $\n$ is the outward unit normal to an interface in the mesh and 
	\begin{equation} \label{eq:alpha}
		\alpha = \frac{\int |\Omegahat\cdot\n| \ud \Omega}{\int \ud\Omega} \rightarrow \frac{\sum_{d=1}^{N_\Omega} w_d\,|\Omegahat_d\cdot\n|}{\sum_d w_d} \approx \frac{1}{2} \,. 
	\end{equation}
Following \cite[Rem.~3]{smm_olivier}, we approximate $\alpha$ with \Sn quadrature to improve consistency with the discrete high-order transport equation. 
The boundary conditions in Eq.~\ref{eq:smm_bc_analytic} are of the form derived in \cite{QDBC} for QD and the SM analog is discussed in \cite{smm_olivier,csmm}. 
Further details on the derivation of the SM low-order system can be found in \cite{AL,lewis_miller,smm_olivier,csmm}. 

Collapsing to a gray low-order problem requires defining suitable gray opacities, $\sigma_E$ and $\sigma_F$, for the zeroth and first moment, respectively. 
Here, we use 
	\begin{subequations}
	\begin{equation} \label{eq:sigmaE}
		\sigma_E = \frac{\sum_{g=1}^G E_g \sigma_g}{\sum_{g=1}^G E_g} \,, 
	\end{equation}
	\begin{equation} \label{eq:sigmaF}
		\sigma_F = \frac{\sum_{g=1}^G r_g(T) \sigma_g}{\sum_{g=1}^G r_g(T)} \,,
	\end{equation}
	\end{subequations}
where $r_g(T)$ is the normalized Rosseland spectrum such that 
	\begin{equation}
		r_g(T) = \int_{\nu_g/T}^{\nu_{g+1}/T} \frac{15}{4\pi^4} \frac{x^4 e^x}{(e^x-1)^2} \ud x \,. 
	\end{equation}
The normalized spectrum has the property that $\pderiv{B_g}{T} = 4 a c T^3\,r_g(T)$. 
Observe that when $\sum_g E_g$ matches its low-order counterpart, $E$, the denominator of Eq.~\ref{eq:sigmaE} cancels with $E$ in the product $\sigma_E E$ in the zeroth moment equation (Eq.~\ref{eq:zeroth}), leaving the original high-order absorption term, $\sum_g \sigma_g E_g$. 
For the first moment's opacity, $\vec{F}$ cannot be used as a weight function in an analogous manner to $\sigma_E$ since $\vec{F}$ is vector-valued and can be negative.  
Instead, we use the Rosseland-collapsed opacity and include an \emph{additive} correction of the form: 
	\begin{equation} \label{eq:additive_mg_correction}
		\sum_{g=1}^G (\sigma_F - \sigma_g) \vec{F}_g \,. 
	\end{equation}
As with $\sigma_E$, when $\sum_g \vec{F}_g$ and $\vec{F}$ match, the terms involving $\sigma_F$ cancel in Eq.~\ref{eq:first}, leaving the high-order term $\sum_g \sigma_g \vec{F}_g$. 
In this way, the use of $\sigma_E$ and the additive correction for the choice of $\sigma_F$ make the \emph{gray} low-order problem consistent with the frequency and direction-integrated high-order transport equation. 
Note that $\sigma_F$ is the Rosseland-collapsed opacity and not the Rosseland opacity defined as $\sigma_R = \sum_g r_g / \sum_g r_g/\sigma_g$ which arises in diffusion theory. 
Within our specific choices of nonlinear solvers, spatial discretizations, and the additive correction for the gray first moment opacity, our choice of $\sigma_F$ as in Eq.~\ref{eq:sigmaF} was observed to be the most efficient and robust. 
Finally, we define the Planck-collapsed gray opacity 
	\begin{equation}
		\sigma_P = \frac{\sum_{g=1}^G b_g(T) \sigma_g}{\sum_{g=1}^G b_g(T)} \,,
	\end{equation}
such that $\sigma_P B(T) \equiv \sum_{g=1}^G \sigma_g B_g(T)$ given that $B_g(T) \equiv B(T) b_g(T)$. 
Use of the Planck-collapsed gray opacity, $\sigma_P$, allows the SM low-order and energy balance equations to be written entirely independent of frequency and is useful in the design of the efficient nonlinear solver for the low-order problem discussed in \S\ref{sec:nonlinear_low_order_solve}. 

Equations \ref{eq:smm_transport}-\ref{eq:smm_meb} are discretized in time using semi-implicit backward Euler. 
Here, semi-implicit refers to the implicit treatment of absorption-emission but the explicit treatment of the temperature dependence of the opacity. 
In other words, the opacities are evaluated using the temperature from the previous time step and iteratively fixed while absorption-emission is resolved implicitly. 
Letting superscript denote time step index and $\Delta t$ the time step at time step $k$, the semi-implicit backward Euler TRT system for advancing from time $t_k$ to $t_{k+1}$ is: 
\begin{subequations} \label{eq:smm_be}
\begin{equation} \label{eq:smm_be_transport}
	\Omegahat\cdot\nabla I_{d,g}^{k+1} + \tilde{\sigma}_g I_{d,g}^{k+1} = \frac{1}{4\pi}\sigma_g B_g(T^{k+1}) + \frac{1}{c\Delta t} I_{d,g}^k\,, \quad \x \in \D\,,\ 1 \leq g \leq G \,,\ 1\leq d\leq N_\Omega\,,
\end{equation}
\begin{equation}
	I_{d,g}^{k+1}(\x) = \frac{1}{4\pi} B_g(T_\text{bdr}) \,, \quad \x \in \partial \D \ \mathrm{and} \ \Omegahat_d\cdot\n < 0 \,,\ 1 \leq g \leq G\,,\ 1\leq d\leq N_\Omega\,, 
\end{equation}
\begin{equation} \label{eq:smm_be_zeroth}
	\nabla\cdot\vec{F}^{k+1} + c \tilde{\sigma}_E E^{k+1} = \sigma_P B(T^{k+1}) + \frac{1}{\Delta t}E^k\,, \quad \x \in \D \,,
\end{equation}
\begin{equation} \label{eq:smm_be_first}
	\frac{c}{3}\nabla E^{k+1} + \tilde{\sigma}_F \vec{F}^{k+1} = \sum_g (\sigma_F - \sigma_g) \vec{F}_g^{k+1} - \nabla\cdot\T^{k+1} + \frac{1}{c\Delta t}\vec{F}^k \,, \quad \x \in \D\,,
\end{equation}
\begin{equation} \label{eq:smm_bc}
	\vec{F}^{k+1}\cdot\n = c \alpha E^{k+1} + \beta^{k+1} + 2 F_\text{in} \,, \quad \x \in \partial \D \,,
\end{equation}
\begin{equation} \label{eq:smm_be_meb}
	\paren{\frac{C_v}{\Delta t} + \sigma_P B(\cdot)}T^{k+1} = c\sigma_E E^{k+1} + \frac{C_v}{\Delta t} T^k \,, \quad \x \in \D \,,
\end{equation}
\end{subequations}
where 
\begin{equation}
	\tilde{\sigma}_g = \sigma_g + \frac{1}{c\Delta t} \,, \quad \tilde{\sigma}_E = \sigma_E + \frac{1}{c\Delta t} \,, \quad \tilde{\sigma}_F = \sigma_F + \frac{1}{c\Delta t} \,, 
\end{equation}
are the multigroup and gray zeroth and first moment opacities including pseudo-absorption from the backward Euler temporal discretization. 
Here, the multigroup opacity $\sigma_g(T^k)$ is evaluated explicitly at the previous time step's temperature. 
While the multigroup opacity is explicit, the collapsing weight functions for the gray opacities are implicit. 
This results in the gray opacities having a hybrid implicit-explicit functional dependence of the form: 
	\begin{subequations}
	\begin{equation}
		\sigma_E(E^{k+1}, T^k) = \frac{\sum_{g=1}^G E_g^{k+1} \sigma_g(T^k)}{\sum_{g=1}^G E_g^{k+1}} \,, 
	\end{equation}
	\begin{equation}
		\sigma_F(T^{k+1}, T^k) = \frac{\sum_{g=1}^G r_g(T^{k+1}) \sigma_g(T^k)}{\sum_{g=1}^G r_g(T^{k+1})} \,,
	\end{equation}
	\begin{equation}
		\sigma_P(T^{k+1}, T^k) = \frac{\sum_{g=1}^G b_g(T^{k+1}) \sigma_g(T^k)}{\sum_{g=1}^G b_g(T^{k+1})} \,. 
	\end{equation}
	\end{subequations}
This treatment of the gray opacities is required to ensure the frequency-integrated high-order and gray low-order problems have equivalent statements of particle balance. 
For example, if $\sigma_E$ were evaluated fully explicitly as $\sigma_E(E^k,T^k)$, the high and low-order systems would differ on the order of the time discretization error due to the explicit spectrum of the energy density used to weight the gray opacity. 
Note that the closures, $\mat{T}^{k+1} = \mat{T}(I^{k+1})$ and $\beta^{k+1} = \beta(I^{k+1})$, are treated implicitly. 
\begin{rem} \label{rem:time}
We have chosen to time integrate the low-order energy density and flux as independent variables. 
This manifests in use of the \emph{low-order} solution as the time-edge sources $E^k$ and $\vec{F}^k$. 
An alternative choice is to use the high-order solution such that the time-edge sources are replaced with $\sum_g 
E_g^k$ and $\sum_g \vec{F}_g^k$. 
The first choice arises from closing the continuous time derivative and the latter from discretizing in time and then closing. 
Both choices converge at backward Euler's first-order temporal order of accuracy and are thus equivalent in the limit of small time steps. 
\end{rem}

We now discuss an iterative solution algorithm for the coupled, nonlinear system of equations given in Eq.~\ref{eq:smm_be}.  
The algorithm is designed such that, upon iterative convergence, Eq.~\ref{eq:smm_be} will hold up to iterative tolerances. 
The algorithm proceeds as follows. 
The high-order transport equation (Eq.~\ref{eq:smm_be_transport}) is inverted on a lagged emission source. 
The resulting intensity is used to compute the closures, $\mat{T}$ and $\beta$. 
The low-order SM system and energy balance equations are solved simultaneously using a Newton-type iteration where temperature is nonlinearly eliminated. 
Finally, the temperature from the low-order solve is used to compute a new emission source for the high-order problem. 
This process is repeated until the low-order energy density and temperature converge to a specified tolerance. 
This algorithm is effectively a fixed-point iteration on the linear coupling between the high and low-order problems through the closures $\mat{T}$ and $\beta$ and a nested, Newton-type iteration to solve the coupled low-order and energy balance equations. 
Acceleration occurs because the nonlinear radiation-material coupling is fully resolved using the low-order system. 
Since the low-order problem is independent of frequency and direction, it is practical to compute and invert the Jacobian associated with each Newton iteration.
Furthermore, the use of nonlinear elimination, where the temperature is updated at each iteration by nonlinearly solving $C_v/\Delta t + B(\cdot)$ locally at each point in space, improves the robustness of the nested, inner iteration by ensuring a physical temperature is generated even at intermediate stages of the solve.  
Converging the nested iteration ensures that the convergence of the outer iteration depends only on the convergence of the closures, as is the case for the simpler linear transport case. 
The closures are designed to be weak functions of the intensity, resulting in a robust and efficient solution procedure \cite{goldin}. 

\section{Transport and Energy Balance Spatial Discretizations} \label{sec:transport_disc}
The semi-discrete transport and energy balance equations in Eqs.~\ref{eq:smm_be_transport} and \ref{eq:smm_be_meb} are discretized in space with linear discontinuous Galerkin finite elements. 
For brevity, we drop the superscript $(\cdot)^{k+1}$ denoting the solution at $t_{k+1}$ and replace the $(\cdot)^k$ superscript, denoting a solution variable at the previous time step, with $(\cdot)^*$. 
Let the domain $\D$ be split into finite elements, $K$, such that $\D = \cup K$. 
The intensities in each discrete energy group and direction and the temperature are approximated in the finite element space $Y_1$, the space of element-wise, piecewise-discontinuous linear polynomials. 
The multigroup absorption opacity, $\sigma_g$, is represented as piecewise-constant and is computed from the average temperature in each element. 
\begin{rem} \label{rem:opacity}
The weight functions used to collapse the multigroup opacity to gray have spatial dependence on each element, embedding spatial dependence into the gray opacities. 
Thus, even though the multigroup opacity is constant in each element, the gray opacities are not. 
In particular, due to the normalization, $\sigma_E = \sum_g E_g \sigma_g / \sum_g E_g$ is a rational polynomial on each element. 
In addition, the Planck and Rosseland spectra have nonlinear dependence on the piecewise-linear temperature, imbuing nonlinear spatial dependence in the gray first moment and Planck-collapsed opacities, $\sigma_F$ and $\sigma_P$, respectively. 
We have found that neglecting the spatial dependence of the collapsing weight functions by approximating the gray opacities as piecewise-constant significantly degrades the quality of the low-order solution. 
Solution quality is commensurate with the high-order solution when $\sigma_E$, $\sigma_F$, and $\sigma_P$ are approximated as piecewise-linear. 
\end{rem}

The spatial discretization for the transport equation is derived by multiplying the transport equation for each $\Omegahat_d$ and group $g$ by a test function $u \in Y_1$ and integrating over each element in the mesh. 
On interior interfaces, upwinding is used such that 
	\begin{equation}
		\Omegahat_d\cdot\n\, \widehat{I}_{d,g} = \begin{cases}
			\Omegahat_d\cdot\n\, I_{d,g,1} \,, & \Omegahat_d\cdot\n > 0 \\
			\Omegahat_d\cdot\n\, I_{d,g,2} \,, & \Omegahat_d\cdot\n < 0 
		\end{cases} \,,
	\end{equation}
where $I_{d,g,i} = I_{d,g}|_{K_i}$ and the indices ``1'' and ``2'' correspond to two arbitrary neighboring elements in the mesh following the convention that the normal vector, $\n$, points from $K_1$ to $K_2$. 
Note that the upwind flux can be equivalently written using the switch functions
	\begin{equation}
		\Omegahat_d\cdot\n\, \widehat{I}_{d,g} = \frac{1}{2}\paren{\Omegahat_d\cdot\n + |\Omegahat_d\cdot\n|}I_{d,g,1} + \frac{1}{2}\paren{\Omegahat_d\cdot\n - |\Omegahat_d\cdot\n|}I_{d,g,2} \,. 
	\end{equation}
Combining like terms, this is equivalent to 
	\begin{equation}\label{eq:upwind}
		\Omegahat_d\cdot\n\, \widehat{I}_{d,g} = \frac{\Omegahat\cdot\n}{2}(I_{d,g,1} + I_{d,g,2}) + \frac{|\Omegahat\cdot\n|}{2}(I_{d,g,1} - I_{d,g,2}) = \Omegahat\cdot\n\avg{I_{d,g}} + \frac{|\Omegahat\cdot\n|}{2}\jump{I_{d,g}} 
	\end{equation}
where, for $u \in Y_1$, 
	\begin{equation}
		\jump{u} = u_1 - u_2 \,, \quad \avg{u} = \frac{1}{2}(u_1 + u_2) \,, 
	\end{equation}
are the jump and average, respectively, with an analogous definition for vector-valued arguments. 
On the boundary of the domain, the numerical flux is 
	\begin{equation}
		\Omegahat_d\cdot\n\, \widehat{I}_{d,g} = \begin{cases}
			\Omegahat_d\cdot\n\, I_{d,g} \,, & \Omegahat_d\cdot\n > 0 \\ 
			\Omegahat_d\cdot\n\, \frac{1}{4\pi}B_g(T_\text{bdr}) \,, & \Omegahat_d\cdot\n < 0 
		\end{cases} \,. 
	\end{equation}
That is, for an element on the boundary, the inflow is specified by the boundary source and the outflow is specified by the numerical solution. 
Note that $|\Omegahat\cdot\n|$ has a discontinuous derivative at $\Omegahat\cdot\n=0$ and thus \Sn quadrature will be slow to accurately compute angular moments of the upwind flux. 
The consistent method we consider in the following section faithfully computes these moments with \Sn quadrature, resulting in exact agreement with the high-order discretization for any \Sn quadrature rule. 
This results in more expensive and detailed correction terms on interior and boundary faces compared to the independent method. 

With these definitions, the weak form of the transport equation is: for each $\Omegahat_d$ and $1\leq g \leq G$, find $I_{d,g} \in Y_1$ such that 
\begin{multline} \label{eq:dgsn_transport}
	\int_{\Gamma_0} \Omegahat\cdot\n\, \jump{u}\avg{I_{d,g}} \ud s + \frac{1}{2}\int_{\Gamma_0} |\Omegahat\cdot\n|\, \jump{u}\jump{I_{d,g}} \ud s + \int_{\Gamma_{b,d}^+} \Omegahat\cdot\n\, u I_{d,g} \ud s - \int \Omegahat\cdot\nablah u\, I_{d,g} \ud \x \\+ \int \tilde{\sigma}_g\, uI_{d,g}\ud \x = \frac{1}{4\pi}\int u\,\sigma_gB_g(T) \ud \x + \int \frac{1}{c\Delta t}\,u I_{d,g}^* - \frac{1}{4\pi} \int_{\Gamma_{b,d}^-} \Omegahat\cdot\n\, u B_g(T_\text{bdr}) \ud s \,, \quad \forall u \in Y_1 \,, 
\end{multline}
where $\Gamma_0$ is the set of unique interior faces in the mesh, $\Gamma_{b,d}^\pm$ are the outflow/inflow portions of the boundary of the domain, and $\nablah$ the element-local gradient. 
The energy balance discretization is also derived by multiplying by a test function $u \in Y_1$. 
The weak form is then: find $T \in Y_1$ such that 
\begin{equation}
	\int u \paren{\frac{C_v}{\Delta t} + \sigma_P B(\cdot)} T \ud \x = \int c\sigma_E\, u E \ud \x + \int \frac{C_v}{\Delta t}\, u T^* \ud \x \,, \quad \forall u \in Y_1 \,. 
\end{equation}
Note that special attention must be given to the bilinear and linear forms depending on $B_g(T)$ and $B(T)$ since these functions are nonlinear in space and thus cannot be integrated exactly with numerical quadrature. 
However, we have found that inexact spatial quadrature, in particular tensor products of two-point quadrature rules, are accurate enough to preserve the spatial convergence rate of the discretization and produce high-quality solutions. 
Note that, in the SM algorithm, the emission source is evaluated at a temperature determined by the low-order system and thus each high-order solve involves inverting the left hand side of Eq.~\ref{eq:dgsn_transport} only. 

\section{Second Moment Discretizations} \label{sec:smm}
In this section, we briefly derive local Discontinuous Galerkin (LDG)-based consistent and independent spatial discretizations for the low-order diffusion system in Eqs.~\ref{eq:smm_be_zeroth} and \ref{eq:smm_be_first}. 
We follow the procedures outlined in Olivier et al.~\cite{csmm} and Olivier and Haut \cite{smm_olivier} for the consistent and independent methods, respectively. 
The consistent discretization is derived via a discrete residual approach -- where the closures are formed by taking the discrete moments of the high-order equation and adding and subtracting the desired diffusion system -- that ensures the low and high-order solutions differ below iterative tolerances and are independent of the discretization resolution upon iterative convergence. 
The key benefit of the discrete residual is that it allows the flexibility to choose any of the DG-type diffusion discretizations with known efficient linear solvers developed in \cite{csmm} without sacrificing consistency with the high-order equation. 
However, in Olivier et al.~\cite{csmm}, it was shown on linear transport problems that the choice of diffusion discretization does impact the iterative convergence rate, especially on under-resolved problems where the terms that account for the discretization difference between high and low-order are large. 

On the other hand, the independent method discretizes the low-order system without regard for the high-order discretization. 
This class of methods also allows the flexibility to discretize the diffusion system such that it can be scalably solved. 
Relaxing the consistency requirement eliminates the discretization-dependent degradation in the outer convergence rate observed for the consistent methods. 
In particular, the independent method was faster to converge than any consistent method on linear transport problems \cite{csmm}. 
However, the independent LDG method was shown to be less accurate than the consistent method, resulting in issues with solution quality on under-resolved meshes. 

Spatial gradients are the primary challenge in deriving suitable discretizations for the low-order system. 
In the extension to TRT, no new derivative terms are introduced, allowing the direct extension of the discretizations derived in \cite{smm_olivier} and \cite{csmm}. 
In particular, the emission and time-edge source terms in the zeroth and first moments can be viewed as the fixed-sources of the steady-state, linear transport schemes. 
We must then only replace the cross sections with the gray opacities, include time-absorption from the backward Euler discretization, formulate more detailed fixed-sources, and include the multigroup opacity correction term now present in the first moment. 

We proceed by outlining details of the LDG discretization that motivates its use here and then extend the discretizations in \cite{csmm} and \cite{smm_olivier} to the TRT equations. 
LDG discretizes the energy density and flux in the spaces $Y_1$ and $[Y_1]^{\dim}$, respectively, where $[Y_1]^{\dim}$ denotes the space of $\dim$-dimensional, vector-valued functions where each component belongs to $Y_1$. 
The LDG diffusion discretization is derived by: discretizing the zeroth and first moment equations with DG, defining numerical fluxes that decouple the flux across interior elements, and forming a Schur complement for the energy density through element-by-element block elimination of the flux. 
Generally, decoupling the flux across interior interfaces in the mesh makes the resulting Schur complement unstable as the mesh size is reduced. 
In \cite{CASTILLO20061307}, it is shown that stability in the limit of mesh refinement is guaranteed with an upwinding procedure where an arbitrary upwinding of the numerical flux for the flux is exactly counterbalanced by the opposite upwinding in the numerical flux for the energy density. 
It is commonplace in DG discretizations of diffusion to include a penalty bilinear form, an additional term in the discretization that penalizes the piecewise-discontinuous discretization toward continuous \cite{Arnold2002}, whose magnitude is controlled by the penalty parameter, $\kappa$. 
With the LDG numerical flux, stability is guaranteed with any choice, including zero, for the interior penalty parameter \cite{10.1007/s10915-007-9130-3}. 
Following \cite[Eq.~4.30]{csmm}, the numerical fluxes are: 
	\begin{subequations}
	\begin{equation}
		\widehat{\vec{F}}\cdot\n = \avg{\vec{F}\cdot\n} + \frac{s}{2}\jump{\vec{F}\cdot\n} + c \kappa \jump{E} \,,
	\end{equation}
	\begin{equation}
		\widehat{\P}\n = \frac{\n}{3}\paren{\avg{E} - \frac{s}{2}\jump{E}} \,, 
	\end{equation}
	\end{subequations}
where $$\P(\x) = \frac{1}{c}\sum_{g=1}^{G} \sum_{d=1}^{N_\Omega} w_d\,\Omegahat_d\otimes \Omegahat_d\, I_{d,g}(\x)$$ is the pressure. 
Here, $s \in \{+1, -1\}$ controls the upwinding. 
When $s=+1$, $\avg{\vec{F}\cdot\n} + s/2\jump{\vec{F}\cdot\n} = \vec{F}_1\cdot\n$ with $\vec{F}_2\cdot\n$ chosen when $s=-1$. 
The minus sign in $\avg{E} - s/2\jump{E}$ enforces that the energy density is upwinded in the direction opposite of the upwinding of $\vec{F}$, ensuring stability. 
The upwinding parameter $s$ is chosen according to 
	\begin{equation}
		s = \begin{cases}
			+1 \,, & \vec{w}\cdot\n > 0 \\ 
			-1 \,, & \vec{w}\cdot\n < 0 
		\end{cases}
	\end{equation}
for a non-zero vector $\vec{w} \in \R^{\dim}$. 
While $\kappa = 0$, the so-called ``minimally dissipative'' choice, is stable, $\kappa \equiv \alpha/2$ with $\alpha$ defined in Eq.~\ref{eq:alpha} was found to be more consistent with the high-order transport system \cite{csmm}. 
The LDG numerical fluxes do not couple the flux across interior mesh interfaces, allowing the flux to be locally eliminated on each element. 
The resulting Schur complement for the energy density represents a diffusion problem cast in second-order differential form. 
The symmetric and positive definite LDG Schur complement can be scalably preconditioned with algebraic multigrid. 
Once the energy density is known through linearly solving the globally coupled Schur complement, element-local back substitution is applied to solve for the flux.
LDG is a natural choice in this application since the energy density and flux are approximated in the same spaces as the high-order energy density and flux. 
In addition, LDG effectively discretizes the hyperbolic conservation law form of the low-order system, where both $E$ and $\vec{F}$ are advanced in time, while allowing an efficient solution procedure within an implicit time integration scheme via algebraic multigrid. 
This discretization then balances being close to the moments of the high-order discretization while avoiding the ill-conditioning associated with equivalence to the moments of the high-order discretization. 
Finally, LDG is a parameter-free discretization and thus avoids the need to hand-tune parameters, such as the penalty parameter used in interior penalty methods \cite{Arnold2002}, in dynamically evolving and algorithmically complex implementations such as the Newton-type iteration discussed in \S\ref{sec:nonlinear_low_order_solve}. 

In \cite{csmm}, the discrete residual approach is used to derive volumetric correction sources and extensions of the LDG numerical fluxes that force the low-order solution to exactly match the moments of the high-order solution. 
Since the low-order problem is diffusion, the resulting correction sources can be interpreted as the sum of model corrections in the form of the SM closures $\mat{T}$ and $\beta$ and corrections for the difference in discretization techniques for the high and low-order systems.  
As discussed in \cite[\S4.4]{csmm}, the discretization correction terms have the property of tending toward zero as the mesh is refined, reflecting the fact that the moments of the high-order discretization and the low-order discretization become equivalent in the limit of mesh refinement. 
Extending \cite[Eqs.~4.43 and 4.44]{csmm} to the TRT system, the resulting low-order problem is: find $(E,\vec{F}) \in Y_1 \times [Y_1]^{\dim}$ such that 
	\begin{subequations}
	\begin{multline}
		\int_{\Gamma_0} \jump{u}\paren{\avg{\vec{F}\cdot\n} + \frac{s}{2}\jump{\vec{F}\cdot\n}} \ud s + \frac{c\alpha}{2} \int_{\Gamma_0} \jump{u}\jump{E} \ud s + \int_{\Gamma_b} c\alpha\, uE  \ud s - \int \nablah u\cdot\vec{F} \ud \x + \int c\tilde{\sigma}_E\, u E \ud \x \\= \int \sigma_P\, u B(T) \ud \x + \int \frac{1}{\Delta t}\, u E^* \ud \x - 2\int_{\Gamma_b} u\,\vec{F}_\text{in} \ud s + \mathcal{R}_\text{cons}^0(u) \,, \quad \forall u \in Y_1 \,, 
	\end{multline}
	\begin{multline}
		\frac{c}{3} \int_{\Gamma_0} \jump{\vec{v}\cdot\n}\paren{\avg{E} - \frac{s}{2}\jump{E}} \ud s + \frac{c}{3}\int_{\Gamma_b} \vec{v}\cdot\n\, E \ud s - \frac{c}{3}\int \nablah \cdot\vec{v}\, E \ud \x + \int \tilde{\sigma}_F\,\vec{v}\cdot\vec{F} \ud \x \\= \int \frac{1}{c\Delta t}\,\vec{v}\cdot\vec{F}^* \ud \x - \int \vec{v}\cdot\vec{P}_\text{in} \ud s + \mathcal{R}_\text{cons}^1(\vec{v})\,, \quad \forall \vec{v} \in [Y_1]^{\dim} \,, 
	\end{multline}
	\end{subequations}
where the correction source terms are
	\begin{subequations}
	\begin{equation}
		\mathcal{R}_\text{cons}^0(u) = -\frac{1}{2}\int_{\Gamma_0} \jump{u}\jump{\beta} \ud s + \frac{s}{2} \int_{\Gamma_0} \jump{u}\jump{\Jho\cdot\n} \ud s - \int_{\Gamma_b} u\paren{\Jhoplus - c\alpha \phiho - \Jin} \ud s \,,
	\end{equation}
	\begin{multline}
		\mathcal{R}_\text{cons}^1(\vec{v}) = \int \vec{v} \cdot \sum_g\paren{\sigma_F - \sigma_g}\vec{F}_{\smallHO,g} \ud \x + \int \nablah\vec{v} : \T \ud \x - \int_{\Gamma_0} \jump{\vec{v}}\cdot\avg{\T\n} \ud s \\- \frac{c}{2}\int_{\Gamma_0} \jump{\vec{v}}\cdot\jump{\Phoplus - \Phominus + \frac{\n}{3}s\phiho} \ud s - \int_{\Gamma_b} c\,\vec{v} \cdot \paren{\Phoplus - \frac{\n}{3}\phiho} \ud s \,. 
	\end{multline}
	\end{subequations}
Here, the ``HO'' subscript denotes a quantity is computed from the high-order intensity and $(\cdot)^{\pm}$ denotes the half-range-integrated quantities 
	\begin{subequations}
	\begin{equation}
		F_{\smallHO,n}^\pm(\x) = \sum_{g=1}^{G}\sum_{\Omegahat_d\cdot\n\gtrless 0} w_d\,\Omegahat_d\,I_{d,g}(\x) \,,
	\end{equation}
	\begin{equation}
		\Phopm(\x) = \frac{1}{c}\sum_{g=1}^{G}\sum_{\Omegahat_d\cdot\n\gtrless 0} w_d\, \Omegahat_d\otimes\Omegahat_d\, I_{d,g}(\x) \,. 
	\end{equation}
	\end{subequations}
Note that, in addition to the terms given in \cite[Eq.~4.44]{csmm}, $\mathcal{R}_\text{cons}^1(\vec{v})$ includes the discretized multigroup opacity correction.
This discretization uses the LDG numerical fluxes resulting in a well-conditioned left hand side operator but includes right-hand-side sources that correct the LDG numerical fluxes to match the moments of the upwind flux used in the transport discretization upon iterative convergence. 

In \cite{smm_olivier}, independent methods are derived by discretizing the continuous low-order problem (e.g.~Eqs.~\ref{eq:smm_be_zeroth} and \ref{eq:smm_be_first}) directly and using knowledge of the high-order discretization only in the evaluation of the closures. 
Compared to the consistent low-order system, this results in a low-order discretization that uses different boundary conditions and neglects the discretization-dependent correction terms. 
An independent, LDG-based low-order problem for TRT is: find $(E,\vec{F}) \in Y_1 \times [Y_1]^{\dim}$ such that 
	\begin{subequations}
	\begin{multline}
		\int_{\Gamma_0} \jump{u}\paren{\avg{\vec{F}\cdot\n} + \frac{s}{2}\jump{\vec{F}\cdot\n}} \ud s + \frac{c\alpha}{2} \int_{\Gamma_0} \jump{u}\jump{E} \ud s + \int_{\Gamma_b} c\alpha\, uE  \ud s - \int \nablah u\cdot\vec{F} \ud \x + \int c\tilde{\sigma}_E\, u E \ud \x \\= \int \sigma_P\, u B(T) \ud \x + \int \frac{1}{\Delta t}\, u E^* \ud \x - 2\int_{\Gamma_b} u\,\vec{F}_\text{in} \ud s \,, \quad \forall u \in Y_1 \,, 
	\end{multline}
	\begin{multline}
		\frac{c}{3} \int_{\Gamma_0} \jump{\vec{v}\cdot\n}\paren{\avg{E} - \frac{s}{2}\jump{E}} \ud s + \frac{c}{3}\int_{\Gamma_b} \vec{v}\cdot\n\, E \ud s - \frac{c}{3}\int \nablah \cdot\vec{v}\, E \ud \x + \int \tilde{\sigma}_F\,\vec{v}\cdot\vec{F} \ud \x \\= \int \frac{1}{c\Delta t}\,\vec{v}\cdot\vec{F}^* \ud \x + \mathcal{R}^1_\text{ind}(\vec{v})\,, \quad \forall \vec{v} \in [Y_1]^{\dim} \,, 
	\end{multline}
	\end{subequations}
where 
	\begin{equation}
		\mathcal{R}^1_\text{ind}(\vec{v}) = \int \vec{v} \cdot \sum_g\paren{\sigma_F - \sigma_g}\vec{F}_{\smallHO,g} \ud \x + \int \nablah\vec{v} : \T \ud \x - \int_{\Gamma_0} \jump{\vec{v}}\cdot\avg{\T\n} \ud s 
	\end{equation}
is the independent method's SM correction source. 
Compared with consistent, the independent method does not have any half-range closure terms, only applies closures to the first moment, and has simpler boundary conditions. 
The additional terms present in the consistent closures were shown to degrade the iterative efficiency of the consistent scheme compared to independent but significantly improved the accuracy and quality of the low-order solution \cite{csmm}. 
However, since these additional terms have magnitude on the order of the low-order discretization's spatial order of accuracy, consistent and independent iteratively converge at equivalent rates when the mesh is resolved enough. 

\section{Nonlinear Low-Order Solve} \label{sec:nonlinear_low_order_solve}
The efficient and robust resolution of the nonlinear absorption-emission physics within the low-order system is a key component of the SM TRT algorithm. 
At each outer iteration, a nonlinear system of the form 
	\begin{equation} \label{eq:absem}
	\begin{bmatrix}
		\mat{M}_F & -\frac{1}{3}\mat{D}^T & \\
		\mat{D} & \mat{M}_E + \mat{P} & -\mat{B}(\cdot) \\ 
		& -\mat{M}_a & \tilde{\mat{B}}(\cdot) 
	\end{bmatrix} 
	\begin{bmatrix} 
		\vec{F} \\ E \\ T 
	\end{bmatrix}
	= \begin{bmatrix} 
		q_F \\ q_E \\ q_T 
	\end{bmatrix}
	\end{equation}
must be solved where $\mat{M}_F$, $\mat{M}_E$, and $\mat{M}_a$ are mass matrices with coefficients $\tilde{\sigma}_F$, $\tilde{\sigma}_E$, and $\sigma_E$, respectively, $\mat{D}$ the divergence operator, $\mat{P}$ the penalty bilinear form with $\kappa \equiv \alpha/2$, and 
	\begin{subequations}
	\begin{equation}
		u^T \mat{B}(T) = \int \sigma_P\, u B(T) \ud \x \,, \quad u \in Y_1 \,,
	\end{equation}
	\begin{equation}
		u^T \tilde{\mat{B}}(T) = \int u \paren{\frac{C_v}{\Delta t} + \sigma_P B(\cdot)}T \ud \x \,, \quad u \in Y_1 \,,
	\end{equation}
	\end{subequations}
are the discrete emission operators. 
Here, $(\cdot)$ indicates an operator is nonlinearly dependent on the argument in the solution vector associated with its location in the $3\times3$ matrix. 
The source terms $q_F$, $q_E$, and $q_T$ include the discretized boundary conditions, time-edge sources, and SM correction sources. 
This system is solved with a Newton-type iteration where temperature is nonlinearly eliminated locally at each point in space. 
At a high level, the algorithm proceeds as follows. 
The Jacobian of the nonlinear system is assembled. 
The energy balance and flux equations are linearly eliminated, forming a Schur complement for the energy density only. 
These linear eliminations are made practical by the block-diagonal nature of the LDG flux and energy balance equations. 
The globally coupled Schur complement is inverted using algebraic multigrid-preconditioned conjugate gradient producing a new energy density. 
Given the new energy density, the temperature is nonlinearly solved.
These steps are repeated until the temperature and energy density converge. 

We now fill in details of the high-level algorithm. 
For the nonlinear equation $f(u) = b$, Newton's iteration is of the form: 
	\begin{equation}
		f(u_0) + \pderiv{f}{u}\biggr\rvert_{u=u_0}(u - u_0) = b \Rightarrow u = u_0 + \paren{\pderiv{f}{u}\biggr\rvert_{u=u_0}}^{-1}(b - f(u_0)) \,, 
	\end{equation}
where $(\cdot)_0$ denotes the previous solution iterate or the initial guess.  
Newton's method applied to Eq.~\ref{eq:absem} yields: 
	\begin{equation}
	\begin{bmatrix}
		\mat{M}_F & -\frac{1}{3}\mat{D}^T & \\
		\mat{D} & \mat{M}_E + \mat{P} & -\mat{B}(\cdot) \\ 
		& -\mat{M}_a & \tilde{\mat{B}}(\cdot) 
	\end{bmatrix} 
	\begin{bmatrix} 
		\vec{F}_0 \\ E_0 \\ T_0 
	\end{bmatrix}
	+ 
	\begin{bmatrix} 
		\mat{M}_F & -\frac{1}{3}\mat{D}^T \\ 
		\mat{D} & \mat{M}_E + \mat{P} & -\delta \mat{B} \\ 
		& -\mat{M}_a & \delta\tilde{\mat{B}}
	\end{bmatrix}
	\begin{bmatrix} 
		\vec{F} - \vec{F}_0 \\
		E - E_0 \\ 
		T - T_0 
	\end{bmatrix}
	= 
	\begin{bmatrix} 
		q_F \\ q_E \\ q_T 	
	\end{bmatrix} \,,
	\end{equation}
where $\delta \mat{B}$ and $\delta\tilde{\mat{B}}$ are linearized emission operators assembled from the bilinear forms
	\begin{equation}
		u^T \delta \mat{B} v = \int \sigma_P \pderiv{B}{T}\biggr\rvert_{T=T_0} \, u v \ud \x = \int 4 \sigma_Pa c T_0^3 \, u v \ud \x \,, \quad u,v \in Y_1\,,
	\end{equation}
	\begin{equation}
		u^T \delta \tilde{\mat{B}} v = \int \paren{\frac{C_v}{\Delta t} + 4 \sigma_P ac T_0^3} uv \ud \x \,, \quad u,v \in Y_1 \,. 
	\end{equation}
Note that the gray opacities are only updated in the outer iteration when a new high-order energy density is available. 
While the temperature-dependence of the collapsing spectra associated with $\sigma_P$ and $\sigma_F$ could be updated inside the inner iteration, we have found that fixing the gray opacities within the inner iteration improves robustness without degrading convergence of the outer iteration. 
Within the linearization, this means that only the gray emission term, $B(T) = ac T^4$, has temperature dependence, yielding a sparser Jacobian with functional form equivalent to a gray, diffusion TRT system.
With this choice for the gray opacities, all operators associated with the energy density and flux are linear.
The Newton approximation simplifies to: 
	\begin{equation}
	\begin{bmatrix} 
		\mat{M}_F & -\frac{1}{3}\mat{D}^T \\ 
		\mat{D} & \mat{M}_E + \mat{P} & -\delta \mat{B} \\ 
		& -\mat{M}_a & \delta\tilde{\mat{B}} 
	\end{bmatrix}
	\begin{bmatrix} 
		\vec{F} \\
		E \\ 
		T - T_0 
	\end{bmatrix}
	= 
	\begin{bmatrix} 
		q_F \\ q_E + \mat{B}(T_0) \\ q_T - \tilde{\mat{B}}(T_0)
	\end{bmatrix} \,. 
	\end{equation}
Note that $\mat{M}_F$ and $\delta\tilde{\mat{B}}$ are both block-diagonal-by-element and can be locally inverted on each element. 
For $\mat{M}_F$, this property is by design of the LDG numerical fluxes that do not couple $\vec{F}$ along interior mesh interfaces while $\delta\tilde{\mat{B}}$ has this property since there are no spatial gradients that induce inter-element coupling in the energy balance equation. 
Using these properties, the Jacobian can be factorized as follows:
	\begin{equation}
	\begin{bmatrix} 
		\mat{I} & \mat{D}\mat{M}_F^{-1} & -\delta\mat{B}(\delta\tilde{\mat{B}})^{-1}\\ 
		& \mat{I} &  \\
		& & \mat{I} 
	\end{bmatrix}
	\begin{bmatrix} 
		\mat{S} \\
		-\frac{1}{3}\mat{D}^T & \mat{M}_F \\
		-\mat{M}_a & & \delta\tilde{\mat{B}} 
	\end{bmatrix}
	\begin{bmatrix} 
		E \\ \vec{F} \\ T - T_0 
	\end{bmatrix}
	= \begin{bmatrix} 
		q_E + \mat{B}(T_0) \\ 
		q_F \\ 
		q_T - \tilde{\mat{B}}(T_0) 
	\end{bmatrix}
	\end{equation}
where 
	\begin{equation} \label{eq:schur}
		\mat{S} = \frac{1}{3}\mat{D}\mat{M}_F^{-1}\mat{D}^T + \mat{M}_E + \mat{P} - \delta\mat{B}(\delta \tilde{\mat{B}})^{-1} \mat{M}_a
	\end{equation}
is the Schur complement. 
Note that the rows and columns of $E$ and $\vec{F}$ have been swapped to expose the upper and lower triangular structure of the factorization. 
This linear system can be solved by forming the source for the Schur complement through back substitution of the first matrix of the factorization, solving the globally coupled Schur complement with algebraic~multigrid-preconditioned conjugate~gradient, and forward substitution of the second matrix of the factorization to solve for the updated flux. 
A standard Newton method would continue the forward substitution to the temperature, resulting in a linear update of the form 
	\begin{equation} \label{eq:linear_T_update}
		T = T_0 + \delta\tilde{\mat{B}}^{-1} (q_T - \tilde{\mat{B}}(T_0) + \mat{M}_a E) \,. 
	\end{equation}
However, such an update is prone to overshooting the temperature and can lead to non-physical negativities at intermediate stages of the Newton iteration that crash the simulation.
Instead, we apply nonlinear elimination where the energy balance equation is solved nonlinearly given an energy density. 
That is, 
	\begin{equation} \label{eq:nonlinear_T_update}
		\tilde{\mat{B}}(T) = \mat{M}_a E + q_T 
	\end{equation}
is solved with a nested Newton iteration.
Note that this nonlinear solve corresponds to the original nonlinear equation (i.e.~the third row of Eq.~\ref{eq:absem}) and not the linearized Jacobian. 
However, in effect, the Newton iteration is equivalent to repeatedly solving Eq.~\ref{eq:linear_T_update} swapping $T\rightarrow T_0$ until the temperature converges and satisfies Eq.~\ref{eq:nonlinear_T_update} to iterative tolerances. 
Since the energy balance equation is not coupled in space, the temperature is eliminated pointwise. 
It has been observed that nonlinear elimination reduces the likelihood of linearization-induced negativities by advancing the temperature in lockstep with the energy density.
\begin{algorithm}
\caption{Full discrete SM algorithm for TRT.}
\label{alg:alg}
\begin{algorithmic}
\WHILE{time $<$ final time}
	\STATE Compute time-edge sources and multigroup opacity
	\WHILE{outer not converged}
		\STATE Compute emission source given previous temperature 
		\STATE Solve high-order problem with transport sweep on emission source 
		\STATE Compute gray opacities, multigroup opacity correction, and SM sources 
		\WHILE{inner not converged}
			\STATE Assemble Jacobian 
			\STATE Linearly eliminate flux and linearized energy balance equation 
			\STATE Solve linearized diffusion problem with algebraic multigrid-preconditioned conjugate gradient 
			\STATE Apply element-by-element back substitution to solve for flux 
			\STATE Given energy density, \emph{nonlinearly} solve energy balance for new temperature 
		\ENDWHILE
	\ENDWHILE
\ENDWHILE
\end{algorithmic}
\end{algorithm}
The full algorithm including time integration, the outer iteration between the high and low-order systems, and the inner, Newton-type iteration are outlined in Algorithm \ref{alg:alg}. 

We conclude this section by showing that the Schur complement in Eq.~\ref{eq:schur} remains positive definite in general and is symmetric and positive definite when the spatial discretization is lumped. 
The sum of the operators $\frac{1}{3}\mat{D}\mat{M}_F^{-1}\mat{D}^T$, the penalty matrix $\mat{P}$, and the mass matrix $\mat{M}_E$ discretize $\int \nabla u \cdot \frac{1}{3\tilde{\sigma_F}} \nabla E \ud \x + \int \tilde{\sigma}_E\, u E \ud \x$ and are thus symmetric and positive definite. 
Note that the stabilization procedure used in the LDG method ensures that $\mat{D}\mat{M}_F^{-1}\mat{D}^T$ is symmetric and positive definite on any mesh regardless of $\mat{P}$. 
Thus, we must only investigate the remaining term, $\delta\mat{B}(\delta\tilde{\mat{B}})^{-1}\mat{M}_a$, that arises from the linear elimination of the energy balance equation associated with the Jacobian. 
Using the equivalence between $L^2$ projection and inverse mass matrices along with the symmetry of the mass matrices $\delta \mat{B}$ and $\delta\tilde{\mat{B}}$, observe that 
	\begin{equation}
	\begin{aligned}
		x^T \delta\mat{B}(\delta\tilde{\mat{B}})^{-1} = y^T &\iff \delta \mat{B} x = \delta\tilde{\mat{B}} y \\
		&\iff \int u\paren{4 \sigma_P ac T_0^3 x - \paren{\frac{C_v}{\Delta t} + 4 \sigma_P a c T_0^3}y}\ud \x = 0 \\
		&\iff y = \Pi\paren{\frac{4 \sigma_P ac T_0^3}{\frac{C_v}{\Delta t} + 4 \sigma_P ac T_0^3}x} \,,
	\end{aligned}
	\end{equation}
where $\Pi$ denotes the $L^2$ projection onto $Y_1$. 
For arbitrary test function $u \in Y_1$, we then have that 
	\begin{equation} \label{eq:pseudo_fission}
		u^T \delta\mat{B}(\delta\tilde{\mat{B}})^{-1} \mat{M}_a E = \int \Pi\!\paren{\frac{4 \sigma_P ac T_0^3}{\frac{C_v}{\Delta t} + 4 \sigma_P ac T_0^3}u}\sigma_E E \ud \x \,. 
	\end{equation}
The term arising in the linearization then has a fission-like functional form. 
For positive temperature, the coefficient within the $L^2$ projection satisfies the bounds  
	\begin{equation*}
		0 \leq \frac{4 \sigma_P ac T_0^3}{\frac{C_v}{\Delta t} + 4 \sigma_P ac T_0^3} < 1 
	\end{equation*}
since $C_v/\Delta t > 0$ makes the denominator larger than the numerator. 
The pseudo-fission term then corresponds to a non-multiplying medium and thus $\mat{S}$ is non-singular.  
As seen in Eq.~\ref{eq:pseudo_fission}, this coefficient scales $\sigma_E$ and thus $\mat{M}_E - \delta\mat{B}(\delta\tilde{\mat{B}})^{-1}\mat{M}_a$ remains positive definite. 

The linearized perturbation to the diffusion system is in general not symmetric since, in Eq.~\ref{eq:pseudo_fission}, the test function, $u$, is inside the $L^2$ projection but the trial function, $E$, is not. 
However, when the spatial discretization is lumped, such as through the use of a nodal quadrature rule, the $L^2$ projection reduces to a pointwise equivalence and the triple mass matrix product term $\delta \mat{B}(\delta\tilde{\mat{B}})^{-1}\mat{M}_a$ reduces to a diagonal matrix. 
Therefore, the Schur complement generically remains positive definite and in the special case of a lumped spatial discretization is also symmetric. 

\section{Numerical Results} \label{sec:numerical_results}
We now present numerical results comparing the accuracy, solution quality, and performance of the consistent and independent SM TRT schemes. 
The above methods were implemented using the MFEM finite element framework \cite{mfem-web,mfem-2024} in the \emph{allium} code package \cite{allium}. 
The low-order systems are preconditioned by \texttt{BoomerAMG} from the \textit{hypre} sparse linear algebra package \cite{hypre}. 
All bilinear and linear forms are lumped through numerical integration with $\dim$-dimensional tensor products of the two-point Gauss-Lobatto quadrature rule (i.e.~trapezoidal rule). 
Quadrature-based lumping ensures all terms in both the low and high-order systems are integrated in a discretely consistent manner. 
The outer and inner SM iterations use a residual-based, relative stopping criterion of the form: 
	\begin{equation} \label{eq:relative_tol}
		\| \vec{u}^{k+1} - \vec{u}^k \|_2 < \epsilon \| \vec{u}^1 - \vec{u}^0 \|_2 \,,
	\end{equation}
where $\vec{u} = \vector{T & E}$, superscripts indicate the iteration index, $\vec{u}^0$ is the initial guess on entry to the iteration, and $\epsilon = 10^{-3}$ is the tolerance. 
The solutions at the previous time step and previous outer iteration are used as the initial guesses to the outer and inner iterations, respectively. 
For the inner iteration, this relative stopping criterion ensures the low-order system is always solved to a factor of $\epsilon$ more accurately than the error associated with each iterate of the outer iteration. 
For the outer iteration, the residual-based relative stopping criterion was observed to help avoid stalls in accuracy arising from insufficiently tight iterative tolerances in temporal convergence studies. 
The linearized low-order diffusion system is solved iteratively with preconditioned conjugate gradient to a residual-based, relative stopping criterion of $10^{-10}$. 
The local, nonlinear elimination step is solved to a residual-based, relative tolerance of $10^{-10}$. 
Note that the nonlinear elimination step is applied independently on each cell and thus convergence is measured locally on each cell. 
We report the maximum iterations required to converge the local nonlinear elimination in each cell over the entire mesh. 
We apply maximum iteration limits of 100, 30, 50, and 40 for the outer fixed-point, inner nonlinear, conjugate gradient, and pointwise elimination iterations, respectively. 

Performance is compared to an unaccelerated transport scheme to measure the effectiveness of the additional work associated with the low-order diffusion system. 
The unaccelerated scheme iteratively solves the TRT system by: 1) applying a sweep on an emission source computed from the previous iterate's temperature or the initial guess and 2) local, nonlinear elimination of the temperature. 
This iteration is terminated with the stopping criterion in Eq.~\ref{eq:relative_tol} with $\epsilon = 10^{-4}$. 
The unaccelerated scheme is solved to a tighter tolerance since, without acceleration, $\|\vec{u}^{1} - \vec{u}^0\|_2$ is larger for the unaccelerated scheme than for an SM method. 
We have found that a $10\times$ tighter tolerance for the unaccelerated scheme produces similar accuracy and solution quality between the two methods.
As with SM, the nonlinear elimination step is solved to a residual-based, relative tolerance of $10^{-10}$. 
For the unaccelerated method, the maximum iteration limits are 500 outer and 40 pointwise elimination iterations, respectively. 

When reported, error is computed with respect to a discrete reference solution that is temporally and spatially refined. 
This requires computing the norm of the difference between solutions defined on distinct meshes. 
The error is computed in the $L^2(\D)$ norm using a $\dim$-fold tensor product of the three-point Gauss-Legendre quadrature rule on each element in the coarser mesh. 
For each quadrature point on the coarse mesh, a search is performed to find the element in the reference mesh that contains the quadrature point. 
The reference solution can then be interpolated using the finite element expansion on each element to evaluate the reference solution at each quadrature point. 
All errors are reported as relative to the magnitude of the reference solution. 
All plots of discrete solutions show the within-element linear variance of the solution by interpolating the solution across a range of points within the element. 

Reported runtimes of algorithms or kernels within an algorithm are aggregated across ranks and represent the maximum time required by a rank to complete. 
Wall time is computed as the sum of the time spent in each cycle minus the time spent in IO and disk operations such as writing simulation data to disk for visualization or simulation restart purposes. 
We will compare the costs of the four most expensive components of the SM algorithm: the high-order transport sweep, the low-order Newton iteration, the linear preconditioned conjugate gradient iteration, and the pointwise Newton iteration for eliminating the temperature. 
The linear solve cost includes AMG setup, forming the Schur complement in the energy density, solving the globally coupled system with preconditioned conjugate gradient, and back-solving for the first moment. 
Note that the cost of the low-order Newton solve includes the nested costs of the linear solve and the pointwise elimination of the temperature. 
We report costs inclusively such that the total cost includes the sweep and low-order solve costs and the low-order solve cost includes the costs of the linear preconditioned conjugate gradient and pointwise elimination. 
The remaining cost at each cycle arises from evaluating the temperature-dependence of the opacities, applying mass matrices for the time integration scheme, and computing source terms for the high-order spatial discretization such as the Planck emission and boundary condition sources. 
These costs are common among the algorithms we present. 

This section comprises four test problems: two one-dimensional problems, a gray Marshak wave and the multigroup, multi-material Larsen problem, and two two-dimensional problems, the gray crooked pipe and a multigroup lattice problem.
The Marshak wave problem exhibits strong diffusive behavior and is thus a stress test for the iterative robustness of the methods. 
Larsen's problem stresses the effectiveness of gray acceleration on a strongly frequency-dependent problem. 
The low cost of these one-dimensional problems facilitates the generation of space-time resolved discrete reference solutions with which the accuracy can be evaluated in space-time convergence studies, enabling investigation into efficiency of the consistent and independent approaches on these problems. 
The gray crooked pipe problem has a challenging geometry that induces ray effects at early times and the temperature-independent opacities produce a strongly diffusive regime at late times, both of which combine to test the overall robustness of the algorithm in terms of the amount of algorithmic hardening procedures, such as negative flux fixups and positivity-preserving flooring, required to avoid simulation crashes. 
Lastly, the multigroup lattice problem uses analytic opacities inspired by real opacity data relevant to inertial confinement fusion (ICF) and thus further evaluates the overall robustness and effectiveness of our approach on an ICF-like proxy problem. 
Unfortunately, the expense of obtaining space-time resolved discrete reference solutions for the two-dimensional problems was prohibitive. 
Thus, we focus on evaluating solution quality and performance by comparing solutions generated on representative ``coarse'' and ``fine'' resolutions. 
The two-dimensional problems are run in parallel with 11 and 44 ranks on the coarse and fine resolutions, respectively, with 4 OpenMP threads per rank. 
These problems are summarized in Table \ref{tab:problem_description}.
\begin{table}
\centering
\caption{Description of the four numerical test problems.}
\label{tab:problem_description}
\begin{tabular}{cccccc}
\toprule
Problem & Dimension & Groups & Number of Angles & Final Time & Discrete Reference \\
\midrule
Marshak Wave & 1 & Gray & 6 & \SI{2.5}{\ns} & Yes \\
Larsen & 1 & 33 & 6 & \SI{10}{\ns} & Yes \\ 
Crooked Pipe & 2 & Gray & 36, 144 & \SI{50}{\sh} & No \\ 
Lattice & 2 & 17 & 64 & \SI{50}{\ns} & No \\
\bottomrule
\end{tabular}
\end{table}

\subsection{Marshak Wave} \label{sub:marshak_wave}
The Marshak wave problem is a one-dimensional, gray benchmark problem designed to induce very stiff, localized absorption-emission physics. 
This problem is very challenging to resolve in space and time due to the large span of the opacity's magnitude as the problem heats from cold. 
Here, we investigate the accuracy of two choices for the low-order system's time-edge source term as discussed in Remark \ref{rem:time}, conduct a space-time accuracy study against a discrete reference solution to verify our implementation and compare accuracy between the consistent and independent schemes, and compare the solution quality and performance of the consistent and independent methods. 

\subsubsection{Definition}
The Marshak wave problem models radiation impinging on a one-dimensional slab where the material properties are strongly temperature dependent. 
The computational domain is $[0,\SI{0.05}{\cm}]$ and the simulation time is \SI{2.5}{\ns}. 
The opacity is $\sigma(T) = 10^{12}/T^3\si{\per\cm}$ and the heat capacity is \SI{3e12}{\erg\per\cm\cubed\per\eV}. 
The material temperature and radiation fields are initially in equilibrium at a temperature of \SI{1}{\eV}. 
At time $t=0$, a strong radiation source at $x=0$ corresponding to a temperature of \SI{1}{\keV} is turned on. 
Initially, the material is cold and extremely optically thick with $\sigma(\SI{1}{\eV}) = 10^{12}\si{\per\cm}$. 
As the material equilibrates with the boundary condition, the opacity dips to $\sigma(\SI{1000}{\eV}) = 10^3\si{\per\cm}$, allowing radiation to penetrate further into the domain, creating a wave with a strong gradient in temperature that moves from left to right through the domain. 
Thus, the opacity spans nine orders of magnitude during the evolution of this problem. 

We consider meshes with 32, 64, 128, 256, and 512 elements resulting in mesh sizes between \SI{1.5625e-3}{\cm} and \SI{9.765625e-5}{\cm} and five time step refinements between \SI{5e-4}{\ns} and \SI{8e-3}{\ns}. 
Gauss-Legendre $S_6$ angular quadrature is used. 
A spatially and temporally refined discrete reference solution is used to assess the accuracy of the methods. 
The reference solution is computed with the consistent SM method on a mesh of \num{8192} elements and a time step of size \SI{1e-4}{\ns}. 
The reference solution is then $16\times$ and $5\times$ more refined in space and time, respectively, than the highest-resolution simulation in this study. 
Use of a discrete reference enables assessment of the accuracy and efficiency of the methods on a challenging, fully nonlinear problem. 

Note that this is an exceptionally challenging problem: even the discrete reference mesh has a cold optical depth of $\sigma(\SI{1}{\eV}) h \approx 10^6$ and a warm optical depth of $\sigma(\SI{1000}{\eV}) h \approx 10^{-3}$. 
The range of mesh sizes considered in the spatial convergence study cover the range of $0.1 \leq \sigma(\SI{1000}{\eV}) h \leq 1.5$ and thus represent convergence in mesh size as the warm opacity becomes resolved in space. 
However, none of the meshes we consider adequately resolve the cold opacity. 
Despite this under-resolution in the thickest parts of the domain, we are still able to show asymptotic convergence in space and time so we still consider the discrete reference spatially resolved enough. 
Furthermore, the large magnitude of the opacity makes this problem exceptionally stiff. 
In particular, the unaccelerated method could not converge rapidly enough to be practical, requiring $\mathcal{O}(10^4)$ iterations to converge on the first time step. 
Thus, for this problem, we only compare the SM methods. 
For these reasons, this problem represents a challenging iterative and robustness stress test for the SM algorithms. 

\subsubsection{Time-Edge Source}
\begin{figure}
	\centering
	\begin{subfigure}{0.4\textwidth}
		\centering
		\includegraphics[width=\textwidth]{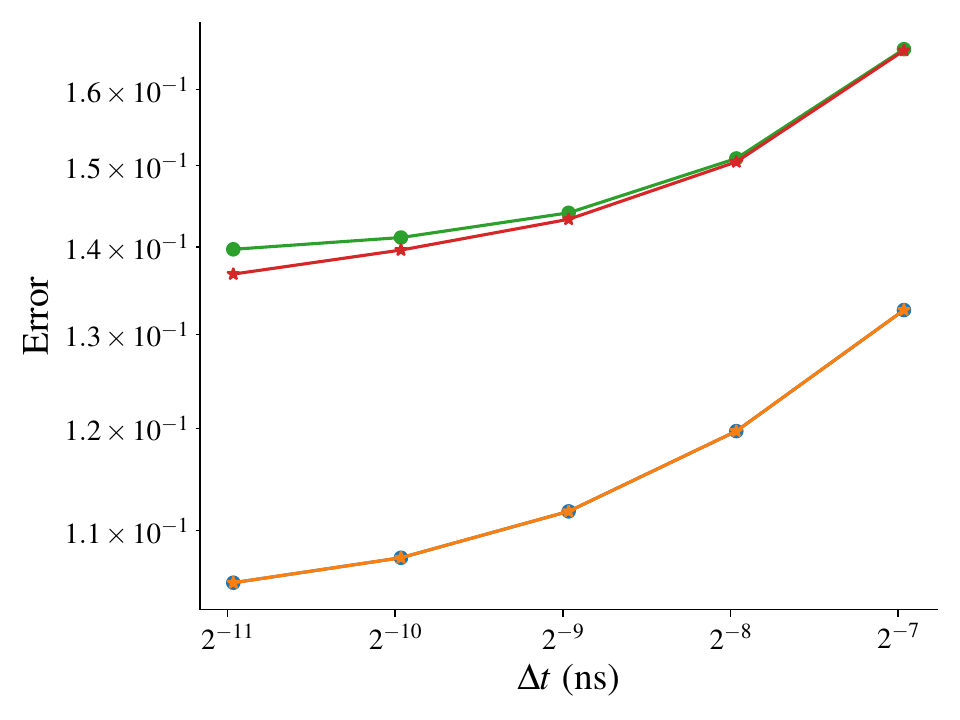}
		\caption{$N_e = 32$}
	\end{subfigure}
	\begin{subfigure}{0.4\textwidth}
		\centering
		\includegraphics[width=\textwidth]{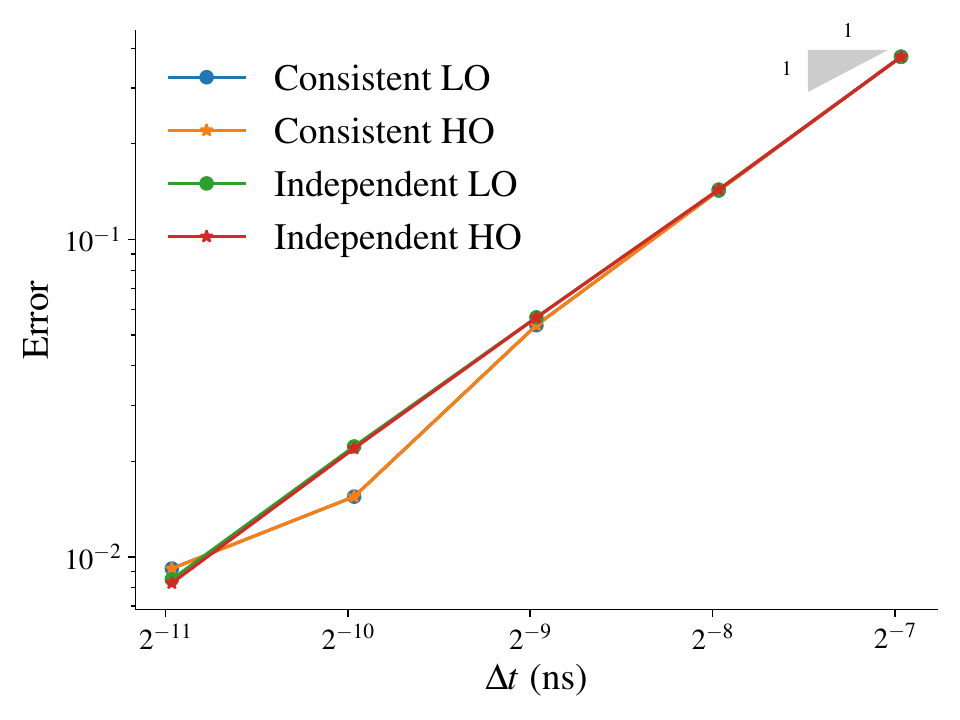}
		\caption{$N_e = 512$}
	\end{subfigure}
	\caption{A comparison of the error with respect to the discrete reference solution when the SM methods use the low and high-order solution as the previous time step's solution.}
	\label{fig:marshak_reset}
\end{figure}
First, we compare the effect of using the high or low-order solution in the time-edge source term for time integration of the low-order system as discussed in Remark \ref{rem:time}. 
Figure \ref{fig:marshak_reset} shows error with respect to the discrete reference solution as the time step is reduced on the coarsest and most refined meshes with 32 and 512 elements, respectively. 
The independent scheme shows a small difference in accuracy for small time step sizes on the coarse mesh. 
On the refined mesh, both schemes produce the same accuracy. 
The consistent scheme is designed so that the difference between the high and low-order solutions is dependent on iterative tolerances and thus the consistent scheme shows identical error behavior for both the coarse and fine meshes. 
These results indicate that optimal convergence in time can be achieved with either method. 
Thus, we elect to use the high-order solution as the time-edge solution since the high-order solution can be guaranteed to be positive through the use of local, negative flux fixups. 
We postulate that this is a marginally more robust approach on under-resolved problems as it results in a more positive source term for the low-order system. 

\subsubsection{Space-Time Accuracy and Efficiency}
\begin{figure}
	\centering
	\begin{subfigure}{0.49\textwidth}
		\centering
		\includegraphics[width=\textwidth]{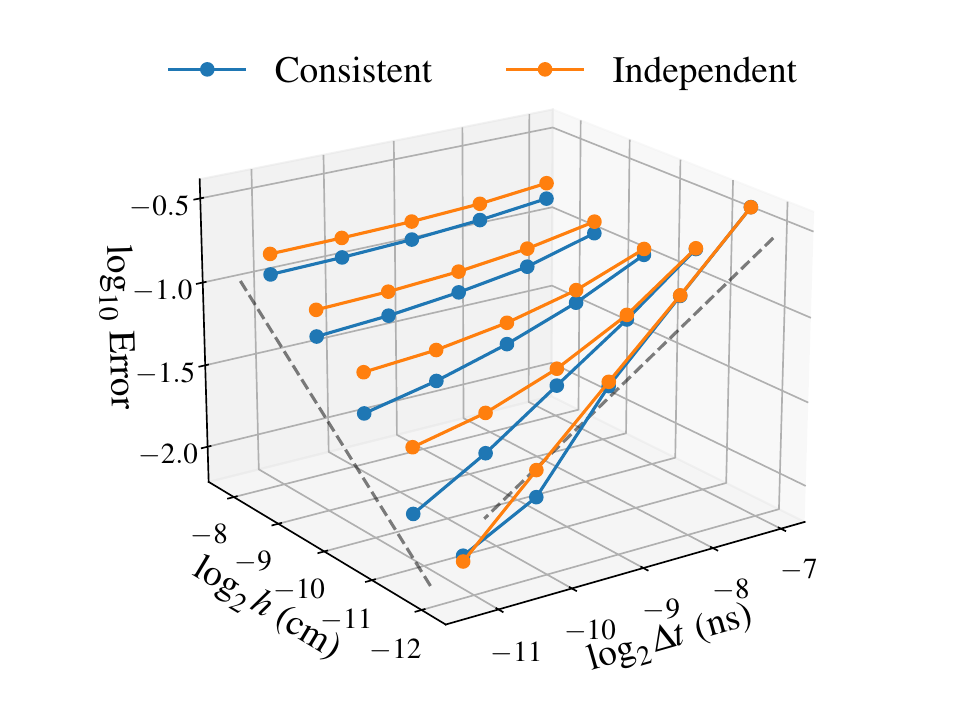}
		\caption{}
		\label{fig:marshak_err}
	\end{subfigure}
	\begin{subfigure}{0.49\textwidth}
		\centering
		\includegraphics[width=\textwidth]{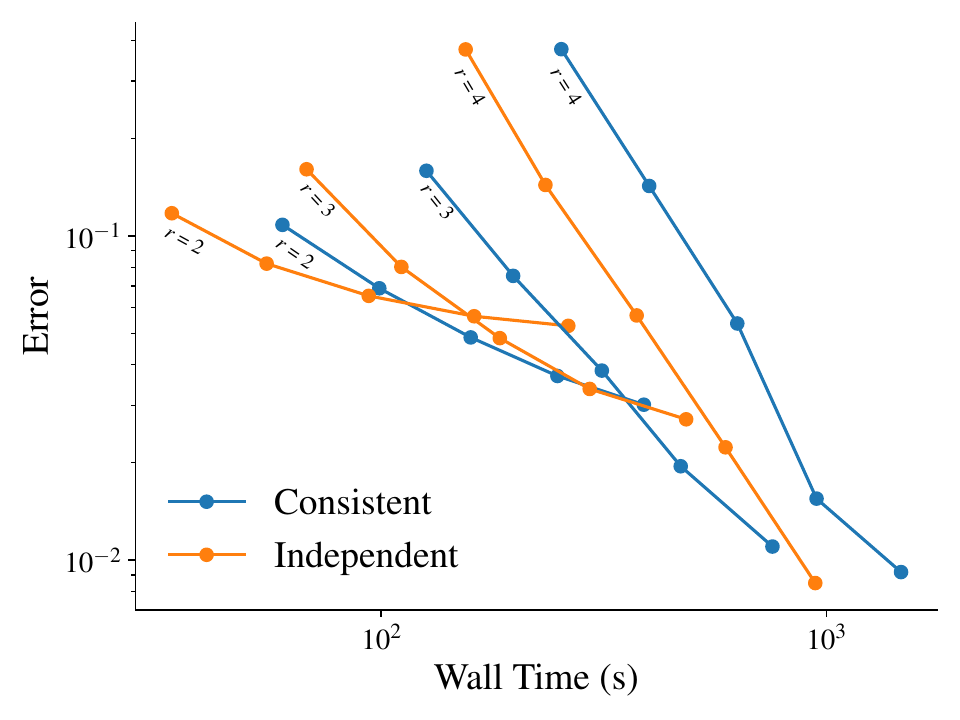}
		\caption{}
		\label{fig:marshak_eff}
	\end{subfigure}
	\caption{The (a) error and (b) efficiency as the mesh size, $h$, and time step size, $\Delta t$, are reduced. In (a), reference first-order lines in space and time are provided. In (b), efficiency is plotted as accuracy versus runtime. Each continuous line represents the efficiency associated with a fixed spatial mesh as the time step is reduced. The line labels indicate the number of mesh refinements with $r=2$, $r=3$, and $r=4$ representing meshes with 128, 256, and 512 elements, respectively.}
\end{figure}
The error in temperature with respect to the discrete reference solution across all spatial and temporal refinements is shown in Fig.~\ref{fig:marshak_err}. 
Reference first-order lines in both space and time are provided. 
For the smallest time step, both the consistent and independent methods converge with first-order accuracy in space. 
While the spatial discretizations are expected to converge with second-order accuracy, convergence may be reduced to first-order on this problem due to the lack of regularity in the solution -- the temperature profile is continuous with a nearly discontinuous derivative at the wave front -- or due to the use of piecewise-constant opacities on this strongly opacity-driven problem.
On the most refined mesh, both schemes converge at the expected first-order accuracy in time. 
Aside from the most resolved spatial mesh, the consistent method is more accurate than the independent method and shows asymptotic convergence in time on coarser meshes than the independent method. 
While consistent is more accurate, it is also more costly. 
The efficiency, in accuracy per simulation wall time, is shown in Fig.~\ref{fig:marshak_eff}. 
Here, continuous lines represent simulations with a fixed mesh size as the time step is reduced. 
The line labels denote the associated mesh has $32 \times 2^r$ elements. 
The fastest method to achieve a given error (i.e.~lowest cost across a fixed horizontal line) alternates between consistent and independent as the desired error is lowered; in certain cases it is more efficient to solve for more unknowns with the independent method than for fewer with the consistent scheme. 
For $r=4$, where the mesh is sufficiently resolved such that consistent and independent produce similar errors, independent is uniformly more efficient than consistent. 

\subsubsection{Solution Quality and Performance}
\begin{figure}
	\centering
	\begin{subfigure}{0.32\textwidth}
		\centering
		\includegraphics[width=\textwidth]{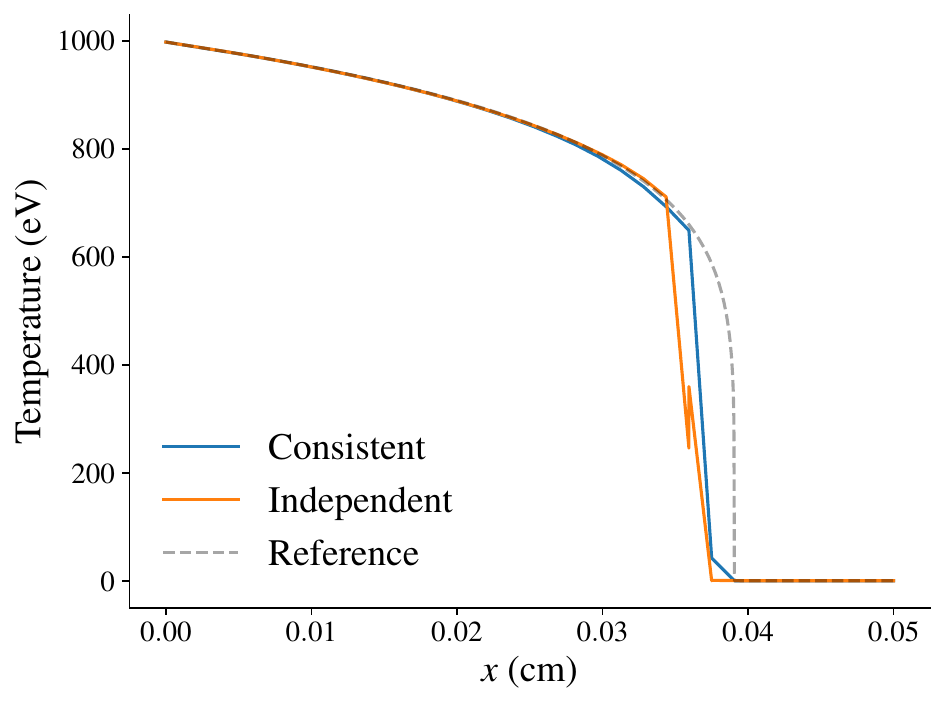}
		\caption{}
	\end{subfigure}
	\begin{subfigure}{0.32\textwidth}
		\centering
		\includegraphics[width=\textwidth]{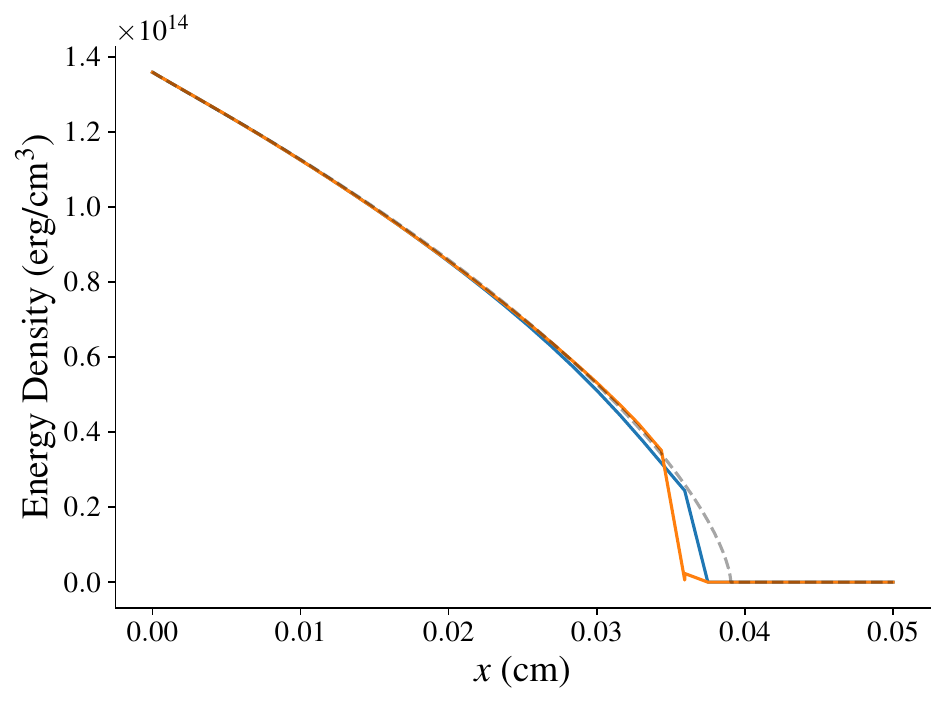}
		\caption{}
	\end{subfigure}
	\begin{subfigure}{0.32\textwidth}
		\centering
		\includegraphics[width=\textwidth]{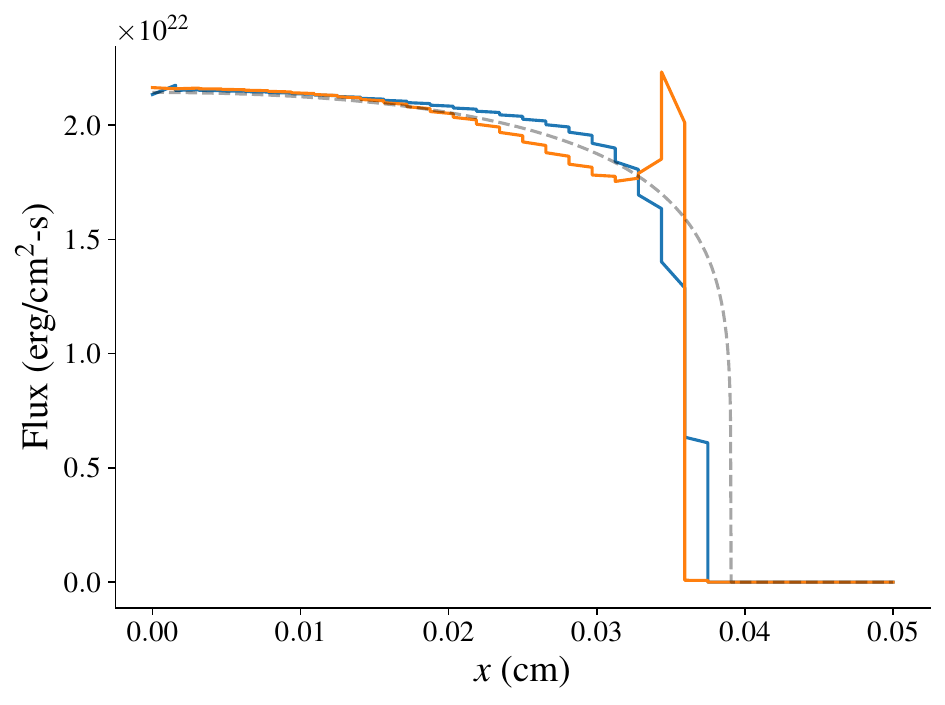}
		\caption{}
	\end{subfigure}
	\caption{The (a) temperature, (b) energy density, and (c) flux on the coarse Marshak wave problem with 32 elements in space and a time step of \SI{4e-3}{\ns}. The discrete reference solution is included. }
	\label{fig:marshak_soln_coarse}
\end{figure}
\begin{figure}
	\centering
	\begin{subfigure}{0.32\textwidth}
		\centering
		\includegraphics[width=\textwidth]{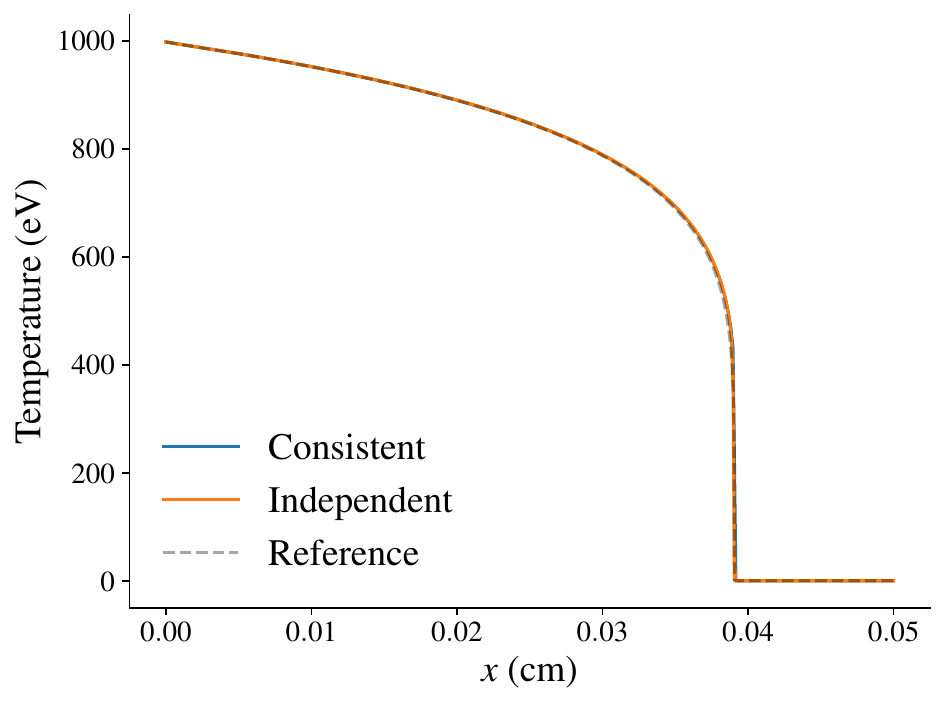}
		\caption{}
	\end{subfigure}
	\begin{subfigure}{0.32\textwidth}
		\centering
		\includegraphics[width=\textwidth]{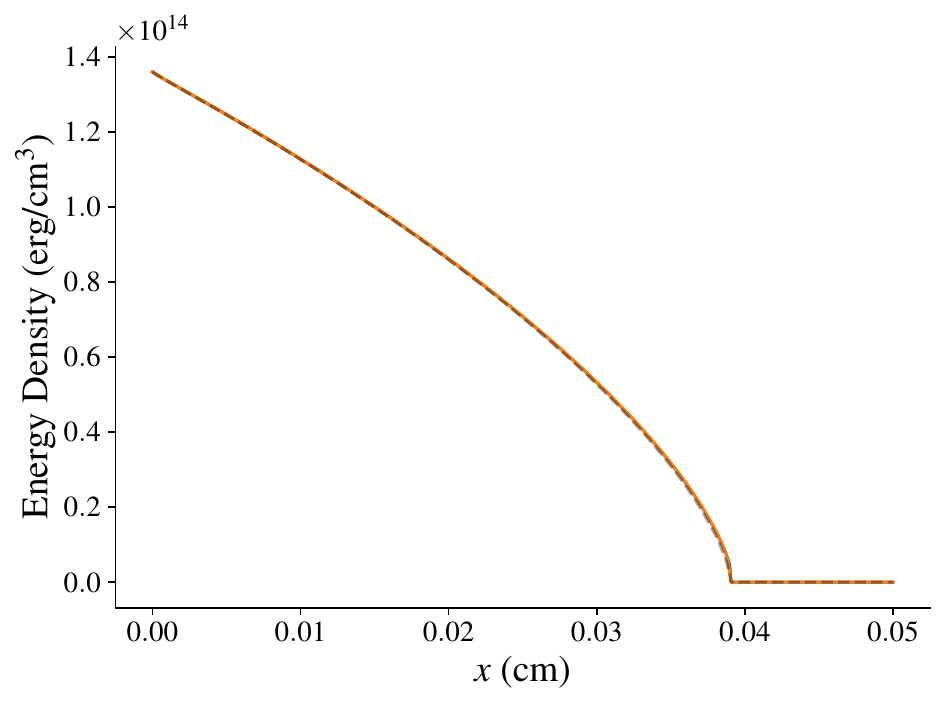}
		\caption{}
	\end{subfigure}
	\begin{subfigure}{0.32\textwidth}
		\centering
		\includegraphics[width=\textwidth]{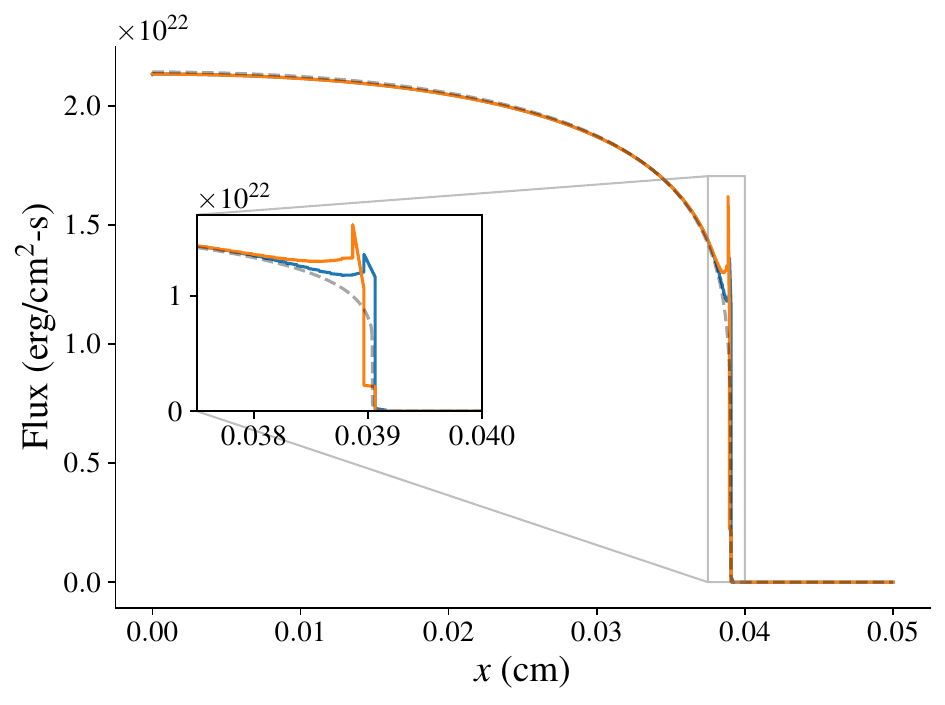}
		\caption{}
	\end{subfigure}
	\caption{The (a) temperature, (b) energy density, and (c) flux on the refined Marshak wave problem with 512 elements in space and a time step of \SI{5e-4}{\ns}. A spatially and temporally resolved reference solution is included. }
	\label{fig:marshak_soln_fine}
\end{figure}

We now investigate the solution quality of the methods. 
Figure \ref{fig:marshak_soln_coarse} shows the temperature, energy density, and flux on the least spatially and temporally refined mesh and time step at the final simulation time compared to the spatially and temporally resolved reference solution. 
The independent method has oscillations in the temperature, energy density, and flux at the wave front whereas consistent maintains a monotonically decreasing profile in all variables. 
Note that we plot the linear variance of the solution in each element by evaluating the solution at a range of points within each element and that, despite the oscillations seen in the element-wise linear solution, the average temperature in each cell remains monotonic. 
Lack of resolution slows the wave speed in both methods with consistent slightly closer to the wave speed of the discrete reference solution. 
It is possible that slope limiting, as applied in \cite{me}, could be applied to reduce the impact of non-monotonic energy densities on the temperature. 
These figures are repeated for the finest mesh and time step in Fig.~\ref{fig:marshak_soln_fine}. 
The consistent and independent solutions are visually indistinguishable in temperature and energy density. 
However, both methods still have a spurious oscillation in the flux with the independent method producing a larger deviation from the discrete reference.

\begin{table}
	\centering
	\caption{Summary of iterative performance and cost of the sweep, low-order nonlinear solve, low-order linear solve, and nonlinear temperature elimination for the consistent and independent SM methods on the Marshak wave problem. Here, ``coarse'' and ``fine'' refer to simulations corresponding to the least and most resolved mesh and time step sizes in the space-time convergence study, respectively.}
	\label{tab:summary_marshak}
	\begin{tabular}{cccccccc}
\toprule
 &  & \multicolumn{2}{c}{Coarse}  &  & \multicolumn{2}{c}{Fine} \\
\cmidrule{3-4}\cmidrule{6-7}
 & Metric & Consistent & Independent & & Consistent & Independent \\
\midrule
\multirow{6}{*}{\rotatebox{90}{Sweeps}} & Min & 4 & 4 & & 3 & 3 \\
 & Max & 11 & 11 & & 11 & 11 \\
 & Avg. & 5.21 & 5.21 & & 5.37 & 5.37 \\
 & Total & \num{3266} & \num{3266} & & \num{53732} & \num{53732} \\
 & Cost (s) & \num{0.6} & \num{0.6} & & \num{113.0} & \num{113.7} \\
 & Avg. Cost (ms) & \num{0.2} & \num{0.2} & & \num{2.1} & \num{2.1} \\
\addlinespace
\multirow{6}{*}{\rotatebox{90}{LO Nonlinear}} & Min & 2 & 2 & & 2 & 2 \\
 & Max & 10 & 10 & & 9 & 9 \\
 & Avg. & 2.57 & 2.57 & & 2.49 & 2.49 \\
 & Total & \num{8456} & \num{8456} & & \num{136619} & \num{136617} \\
 & Cost (s) & \num{15.3} & \num{15.8} & & \num{1246.6} & \num{1271.1} \\
 & Avg. Cost (ms) & \num{4.7} & \num{4.8} & & \num{23.2} & \num{23.7} \\
\addlinespace
\multirow{6}{*}{\rotatebox{90}{LO Linear}} & Min & 1 & 1 & & 1 & 1 \\
 & Max & 14 & 14 & & 15 & 15 \\
 & Avg. & 10.62 & 10.62 & & 12.59 & 12.59 \\
 & Total & \num{89792} & \num{89791} & & \num{1720502} & \num{1720559} \\
 & Cost (s) & \num{8.2} & \num{8.2} & & \num{417.9} & \num{423.9} \\
 & Avg. Cost (ms) & \num{1.0} & \num{1.0} & & \num{3.1} & \num{3.1} \\
\addlinespace
\multirow{6}{*}{\rotatebox{90}{Energy Balance}} & Min & 1 & 1 & & 1 & 1 \\
 & Max & $40^*$ & $40^*$ & & $40^*$ & $40^*$ \\
 & Avg. & 4.16 & 4.15 & & 4.12 & 4.12 \\
 & Total & \num{35139} & \num{35088} & & \num{562751} & \num{562568} \\
 & Cost (s) & \num{1.8} & \num{1.8} & & \num{172.6} & \num{173.8} \\
 & Avg. Cost (ms) & \num{0.2} & \num{0.2} & & \num{1.3} & \num{1.3} \\
\addlinespace
 & Wall Time (s) & \num{16.3} & \num{16.9} & & \num{1440.6} & \num{1467.5} \\
\bottomrule
\end{tabular}
	$^*$ indicates maximum iterations reached
\end{table}
The iterative performance and costs of the algorithms on the coarsest and finest mesh and time step sizes are presented in Table \ref{tab:summary_marshak}. 
The consistent and independent methods on both the coarsest and finest problems exhibit robust iterative performance: the average number of iterations required for each component of the algorithm matches closely. 
We stress that the unaccelerated method would have required significantly more sweeps and total runtime meaning the low-order operations still represent a large improvement in overall efficiency compared to the unaccelerated scheme. 
This will be explored in subsequent sections where it is practical to obtain a solution with the unaccelerated method. 

\subsection{Larsen}
The Larsen problem is a multi-material, frequency-dependent, one-dimensional benchmark. 
The frequency dependence of the opacities is designed to spatially separate the fast and slow dynamics associated with high and low-frequency radiation, respectively, resulting in a problem with a truly frequency-dependent solution that cannot be captured with a gray transport model. 
Thus, this problem stresses our use of a gray low-order system to accelerate the multigroup high-order system. 
This subsection proceeds as follows: we define the problem, compare the use of piecewise-constant and linear gray opacities in the low-order system as discussed in Remark \ref{rem:opacity}, conduct a space-time accuracy study with respect to a discrete reference solution, and compare performance of the consistent, independent, and unaccelerated schemes. 

\subsubsection{Definition}
\begin{figure}
	\centering
	\includegraphics[width=0.4\textwidth]{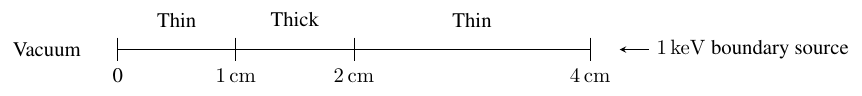}
	\caption{Materials and boundary conditions for Larsen's problem.}
	\label{fig:larsen_diag}
\end{figure}
The domain is $[0,\SI{4}{\cm}]$. 
The problem consists of a thick region between two thin regions with a \SI{1}{\keV} boundary source on the right edge of the domain at $x=\SI{4}{\cm}$. 
The material layout is depicted in Fig.~\ref{fig:larsen_diag}. 
The thick and thin opacities have the form 
	\begin{equation}
		\sigma(T) = \frac{\alpha}{\nu^3}\paren{1 - e^\frac{-\nu}{T}} \,,
	\end{equation} 
where $\alpha_\text{thick} = 10^{12}\si{\eV\cubed\per\cm}$ and $\alpha_\text{thin} = 10^9\si{\eV\cubed\per\cm}$, respectively. 
33 logarithmically spaced groups between \SI{1e-2}{\eV} and \SI{300}{\keV} are used. 
On this problem, the optically thin material is transparent to the high-frequency portion of the spectrum emitted by the \SI{1}{\keV} boundary source, resulting in the inner optically thick material heating much more rapidly than the thin region adjacent to the boundary source. 
This spatially separates the dynamics of the high and low-frequency radiation, creating an inherently frequency-dependent solution. 
The final time is \SI{10}{\ns} to allow for the temperature wavefront to pass through the thick region into the final thin region. 
A discrete reference solution computed with \num{16408} elements and a time step of size \SI{2.5e-4}{\ns} is used to compare solution quality.
The reference mesh is built from a base of 4096 elements per centimeter and has non-uniform refinement at the domain boundaries and on both sides of all material interfaces. 
The reference mesh has a maximum element width of \SI{2.4e-4}{\cm} and minimum width $20\times$ smaller than the maximum. 
Gauss-Legendre $S_6$ angular quadrature is used. 

\subsubsection{Gray Opacity Representations}
\begin{figure}
	\begin{subfigure}{0.32\textwidth}
		\centering
		\includegraphics[width=\textwidth]{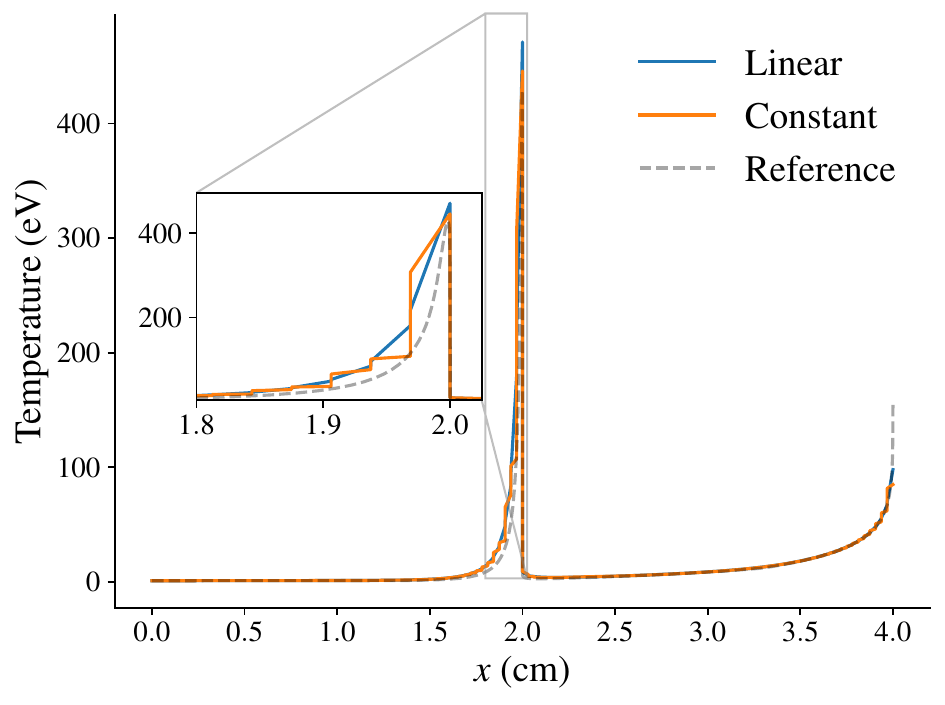}
		\caption{$t = \SI{1e-1}{\ns}$}
        \label{fig:larsen_opacity_a}
	\end{subfigure}
	\begin{subfigure}{0.32\textwidth}
		\centering
		\includegraphics[width=\textwidth]{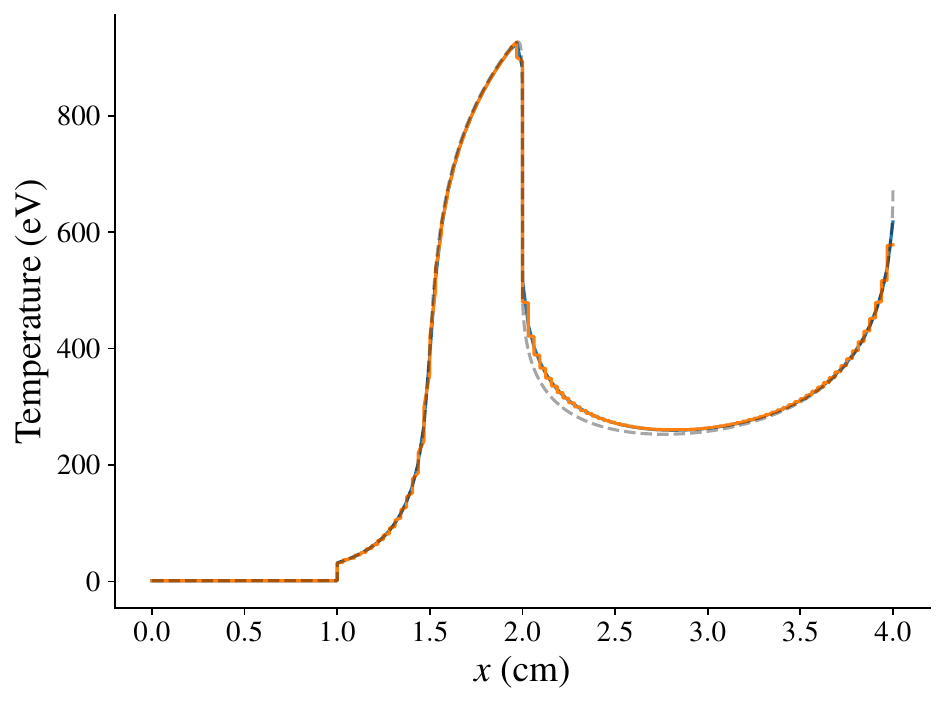}
		\caption{$t = \SI{1}{\ns}$}
	\end{subfigure}
	\begin{subfigure}{0.32\textwidth}
		\centering
		\includegraphics[width=\textwidth]{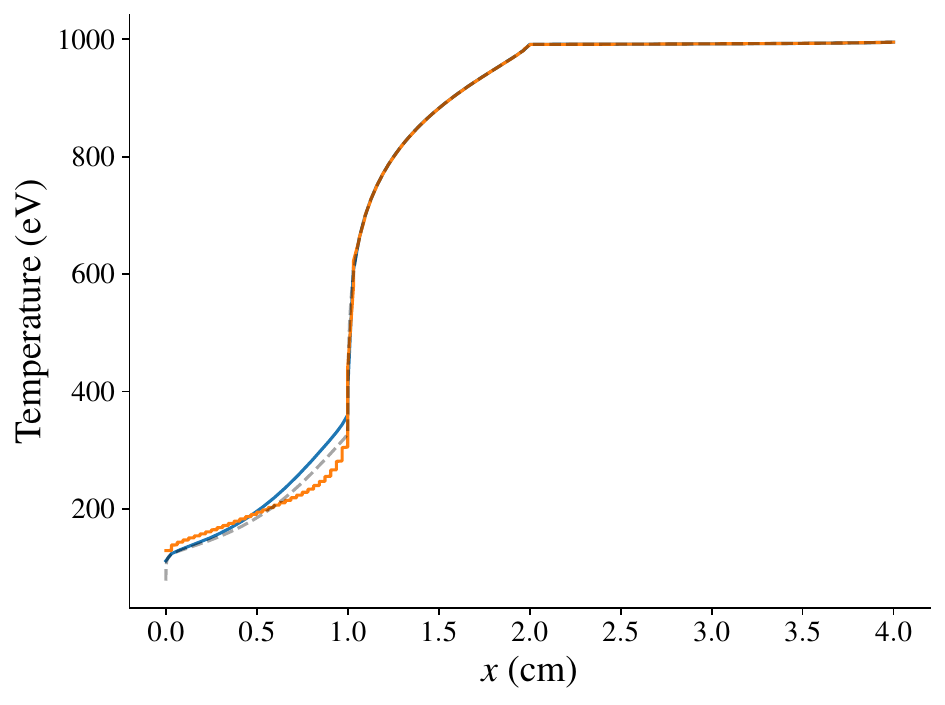}
		\caption{$t = \SI{10}{\ns}$}
		\label{fig::laren_opacity_c}
	\end{subfigure}
	\caption{Temperature solutions generated with the consistent SM method using piecewise-linear and piecewise-constant representations for the gray opacities on Larsen's problem.}
	\label{fig:larsen_opacity}
\end{figure}
We first compare the use of piecewise-linear versus piecewise-constant representations of the \emph{gray} opacities.
As discussed in Remark \ref{rem:opacity}, while the multigroup opacities are piecewise-constant, the collapsing spectra for the gray opacities have spatial dependence. 
Here, we show that solution quality is significantly improved when the spatial dependence of the collapsing spectra is incorporated into the representation of the gray opacities such that $\sigma_E$, $\sigma_F$, and $\sigma_P$ are approximated in $Y_1$. 
The consistent method is used on a mesh with 128 elements and a time step size of \SI{1e-2}{\ns}. 
The temperature at three snapshots in time when piecewise-linear and piecewise-constant gray opacities are used is shown in Fig.~\ref{fig:larsen_opacity}. 
The piecewise-constant gray opacity imprints a piecewise-constant spatial dependence onto the temperature, resulting in the staircase behavior seen in Fig.~\ref{fig:larsen_opacity}. 
In contrast, use of piecewise-linear gray opacities results in a smoothly varying solution. 
In addition, compared to the discrete reference, piecewise-constant has a slower wave speed and less accurate wave shape than piecewise-linear. 
Thus, we choose to use piecewise-linear gray opacities for the remainder of the document. 

\subsubsection{Space-Time Accuracy and Efficiency}
\begin{figure}
	\centering
	\begin{subfigure}{0.45\textwidth}
		\centering
		\includegraphics[width=\textwidth]{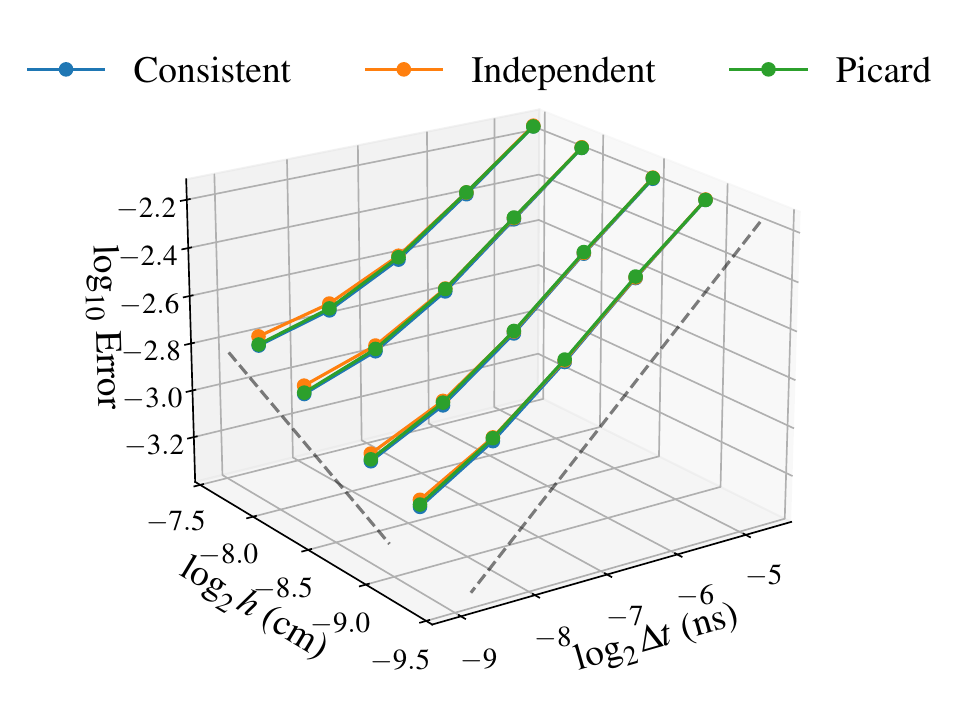}
		\caption{Temperature, $t=\SI{1}{\ns}$}
	\end{subfigure}
	\begin{subfigure}{0.45\textwidth}
		\centering
		\includegraphics[width=\textwidth]{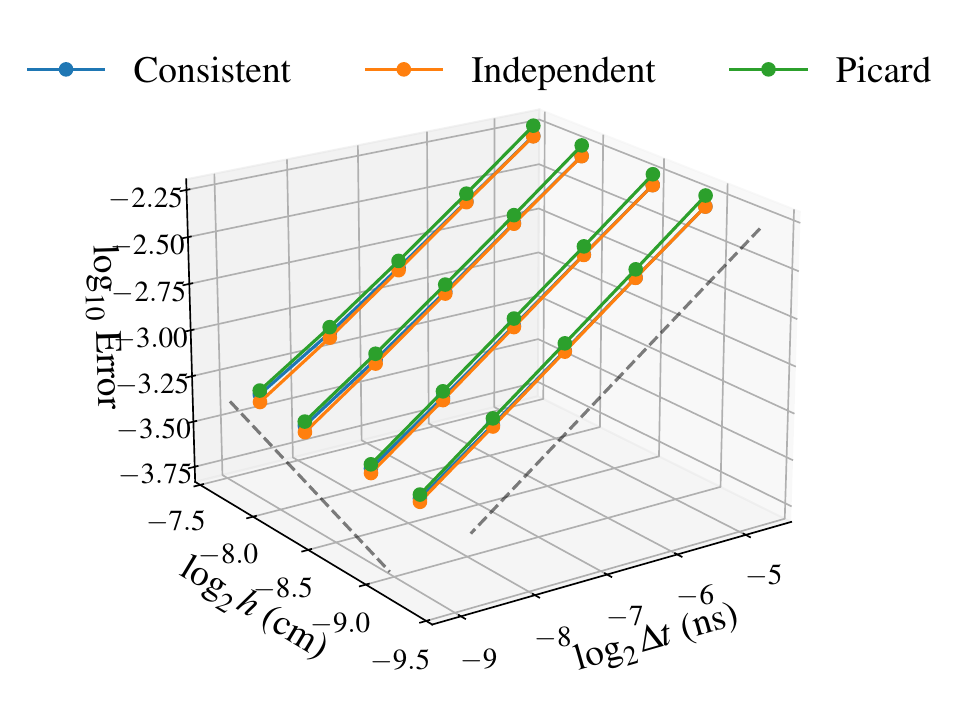}
		\caption{Energy Density, $t=\SI{1}{\ns}$}
	\end{subfigure}
	\begin{subfigure}{0.45\textwidth}
		\centering
		\includegraphics[width=\textwidth]{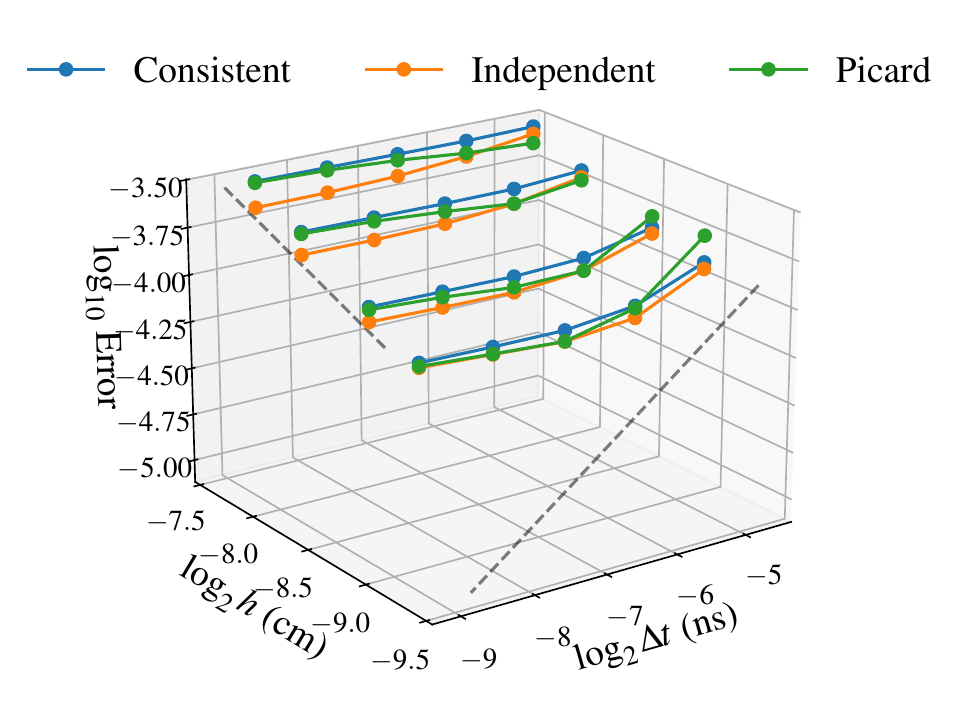}
		\caption{Temperature, $t=\SI{10}{\ns}$}
	\end{subfigure}
	\begin{subfigure}{0.45\textwidth}
		\centering
		\includegraphics[width=\textwidth]{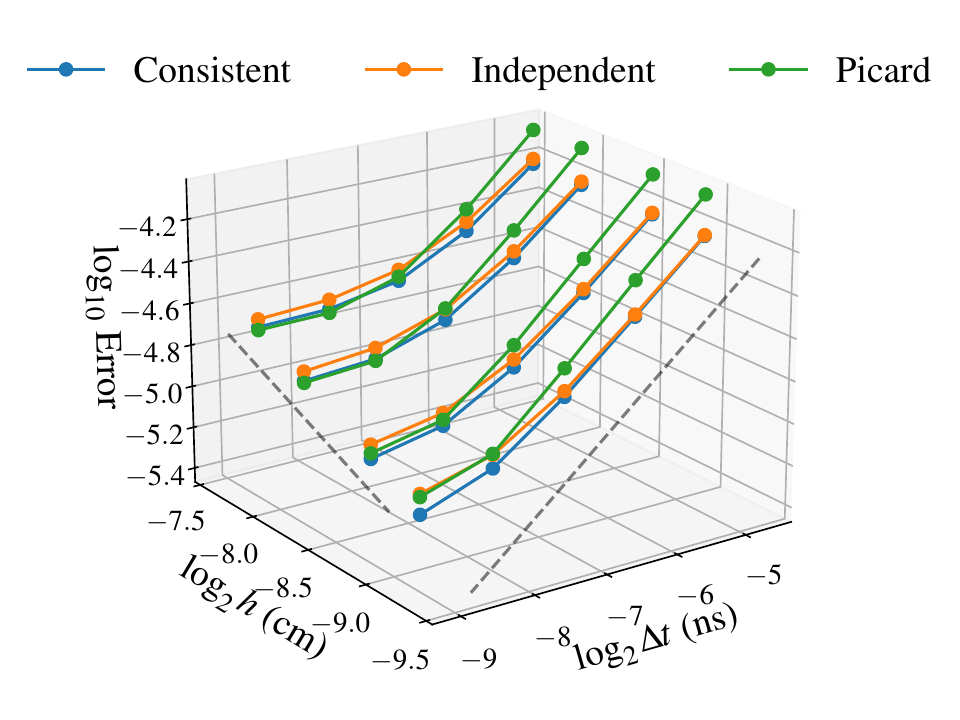}
		\caption{Energy Density, $t=\SI{10}{\ns}$}
	\end{subfigure}
	\caption{Space-time accuracy with respect to a space-time resolved discrete reference solution on Larsen's problem in the temperature and energy density at $t=\SI{1}{\ns}$ and $t=\SI{10}{\ns}$. Reference first-order-accurate lines are provided in both space and time. }
	\label{fig:larsen_error}
\end{figure}
Figure \ref{fig:larsen_error} shows a space-time convergence study for the temperature and energy density at $t=\SI{1}{\ns}$ and $t = \SI{10}{\ns}$. 
We consider four spatial meshes that have a base of 192, 256, 384, and 512 elements per centimeter with non-uniform refinement on both sides of material interfaces. 
The resulting meshes have 786, \num{1042}, \num{1554}, and \num{2066} elements with maximum element widths of \SI{5.2e-3}{\cm}, \SI{3.9e-3}{\cm}, \SI{2.6e-3}{\cm}, and \SI{1.9e-3}{\cm}, respectively. 
The minimum element widths are $13\times$ smaller than the maximum for each mesh. 
Four time step sizes between \SI{2.5e-3}{\ns} and \SI{4e-2}{\ns} are used. 
Accuracy is computed with respect to the discrete reference. 
Convergence in space is assessed with respect to the \emph{maximum} element width in the mesh. 
Different convergence behavior is seen for the temperature and energy density at the earlier and later times. 
The consistent, independent, and unaccelerated methods produce similar error behavior for all meshes and time step sizes. 
At $t=\SI{1}{\ns}$, the temperature and energy density converge with first-order accuracy in time. 
In space, the temperature is also first-order while the energy density has reduced convergence. 
At $t=\SI{10}{\ns}$, the spatial error dominates temporal error in both variables, resulting in stagnating errors. 
For the energy density, temporal convergence stagnates only for the smallest time step with good first-order convergence in both space and time observed on all other time steps. 
The temperature stagnates for all time steps but shows first-order accuracy in space. 
As with the Marshak wave problem, the observed first-order convergence in space may be due to the lack of regularity in the solutions (e.g.~the temperature in Fig.~\ref{fig:larsen_opacity_a}) or representing the multigroup absorption opacity as piecewise-constant. 
While the non-uniform meshes do improve convergence, we suspect that this problem is difficult to resolve in space and time without dynamic refinement. 
Notably, the independent method produces a more accurate temperature than the consistent method at $t=\SI{10}{\ns}$. 
However, this may be due to cancellation of errors within the more numerically diffusive independent method as neither method shows asymptotic convergence at $t=\SI{10}{\ns}$. 

The unaccelerated method exhibits accuracy between that of the independent and consistent SM methods. 
While unaccelerated and consistent SM are formally equivalent up to iterative tolerances, recall that different stopping criteria are used for the unaccelerated and SM methods. 
The fact that the unaccelerated method produces accuracy similar to that of the consistent SM method indicates that our choice of tolerance is fair. 

\begin{figure}
	\centering
	\begin{subfigure}{0.4\textwidth}
		\centering
		\includegraphics[width=\textwidth]{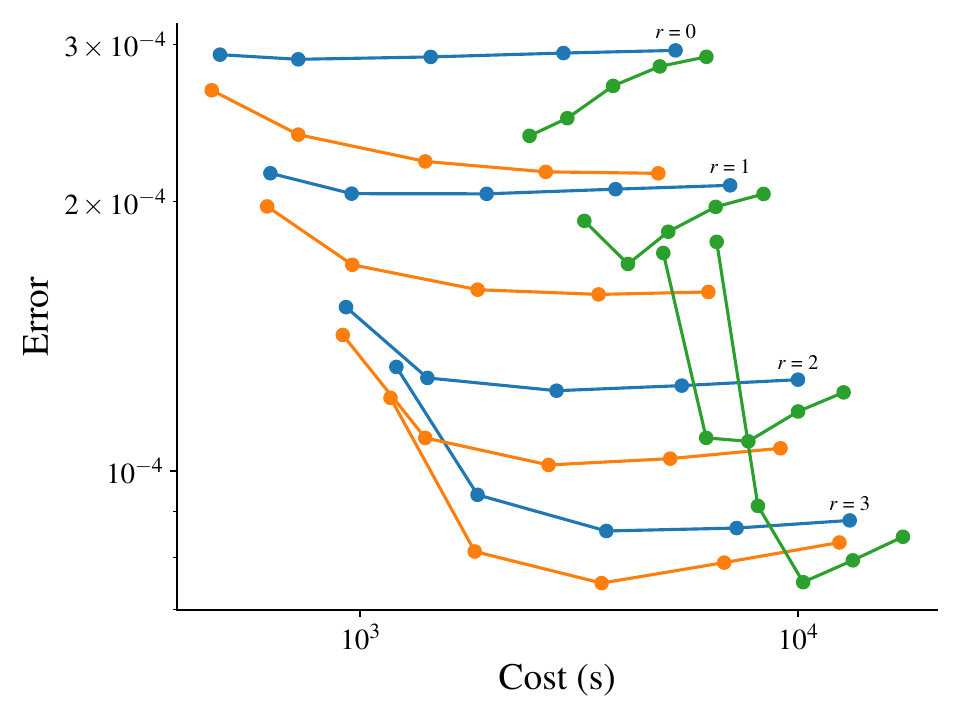}
		\caption{Temperature, $t=\SI{10}{\ns}$}
	\end{subfigure}
	\begin{subfigure}{0.4\textwidth}
		\centering
		\includegraphics[width=\textwidth]{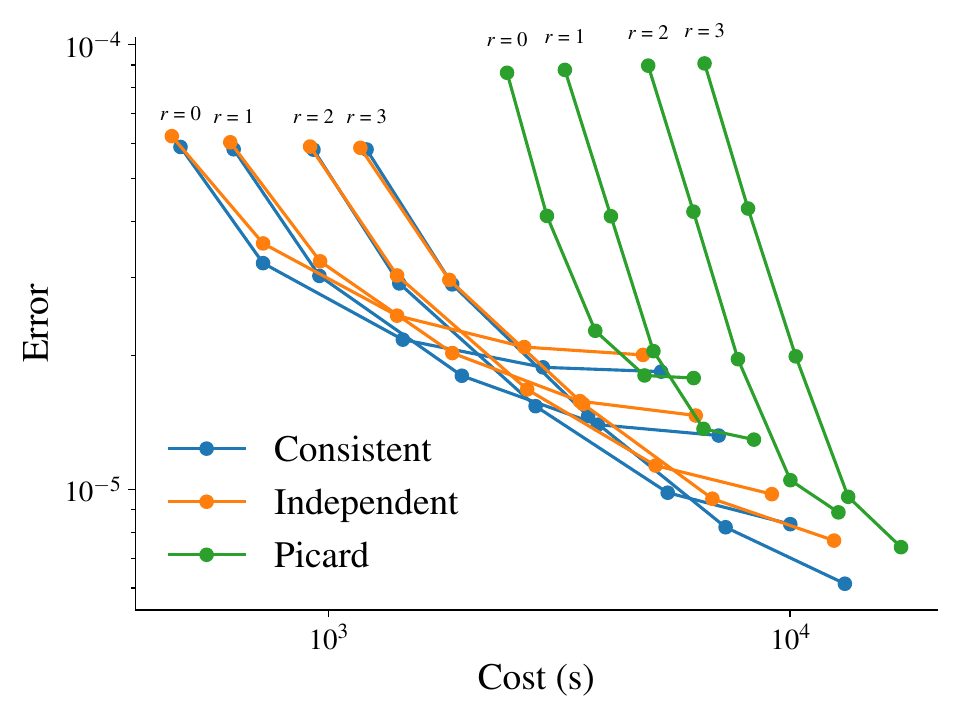}
		\caption{Energy Density, $t=\SI{10}{\ns}$}
	\end{subfigure}
	\caption{Space-time efficiency in terms of accuracy versus runtime on Larsen's problem. Each continuous line represents the efficiency associated with a fixed spatial mesh as the time step is reduced. The line labels indicate successively refined spatial meshes.}
	\label{fig:larsen_eff}
\end{figure}
The efficiency of the methods in temperature and energy density at $t=\SI{10}{\ns}$ is shown in Fig.~\ref{fig:larsen_eff}. 
As with the Marshak problem, continuous lines represent simulations with a fixed mesh as the time step is reduced. 
Each pair of continuous lines is labeled such that $r=0$ is the coarsest mesh and $r=3$ the finest. 
On each mesh and time step size, consistent and independent are similar in cost. 
Thus, the slight differences in accuracy shown in Fig.~\ref{fig:larsen_error} drive the efficiency of the methods. 
Looking at vertical lines of fixed cost, the independent method is more efficient than consistent in the temperature at $t=\SI{10}{\ns}$ while consistent has a slight advantage in efficiency in terms of the energy density. 
In contrast, the unaccelerated method is similar in accuracy to the SM methods but much more expensive. 
However, as the time step reduces, convergence for the unaccelerated method improves. 
For small enough time steps, unaccelerated can converge in a similar number of iterations as the SM methods without performing the extra low-order operations. 
Thus, it is expected that the unaccelerated method will be more efficient for very small time steps. 
Here, we see that the SM methods are slightly more efficient than the unaccelerated method for the smallest time steps but much more efficient as the time step increases, a testament to the SM algorithm's improved robustness with respect to cost as the time step increases. 

\subsubsection{Performance}
\begin{figure} 
	\centering
	\begin{subfigure}{0.45\textwidth}
		\centering
		\includegraphics[width=\textwidth]{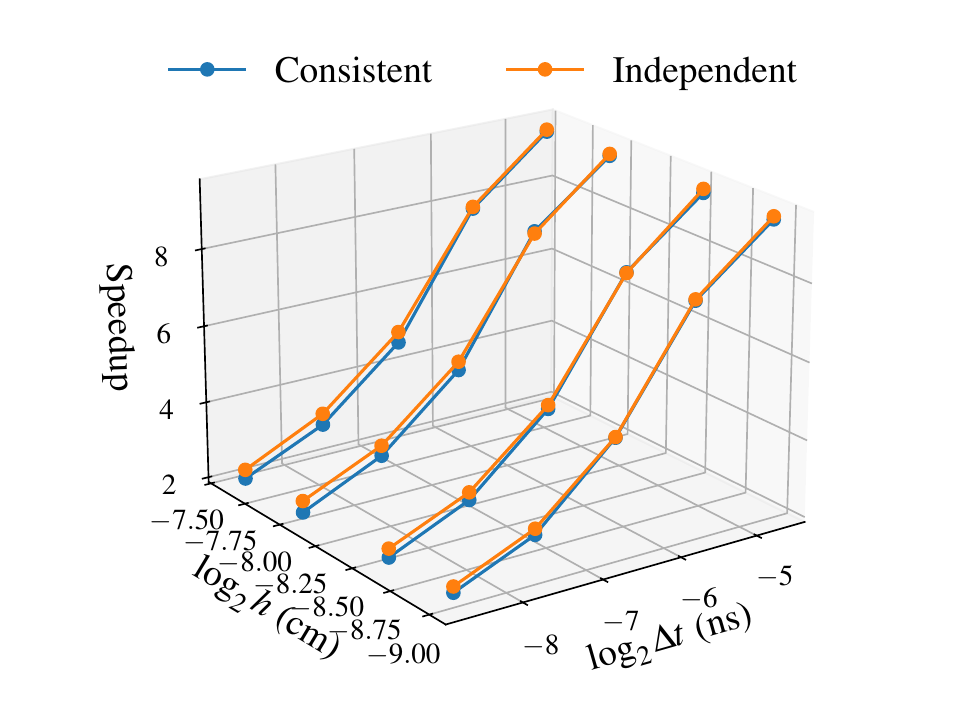}
		\caption{Sweeps}
	\end{subfigure}
	\begin{subfigure}{0.45\textwidth}
		\centering
		\includegraphics[width=\textwidth]{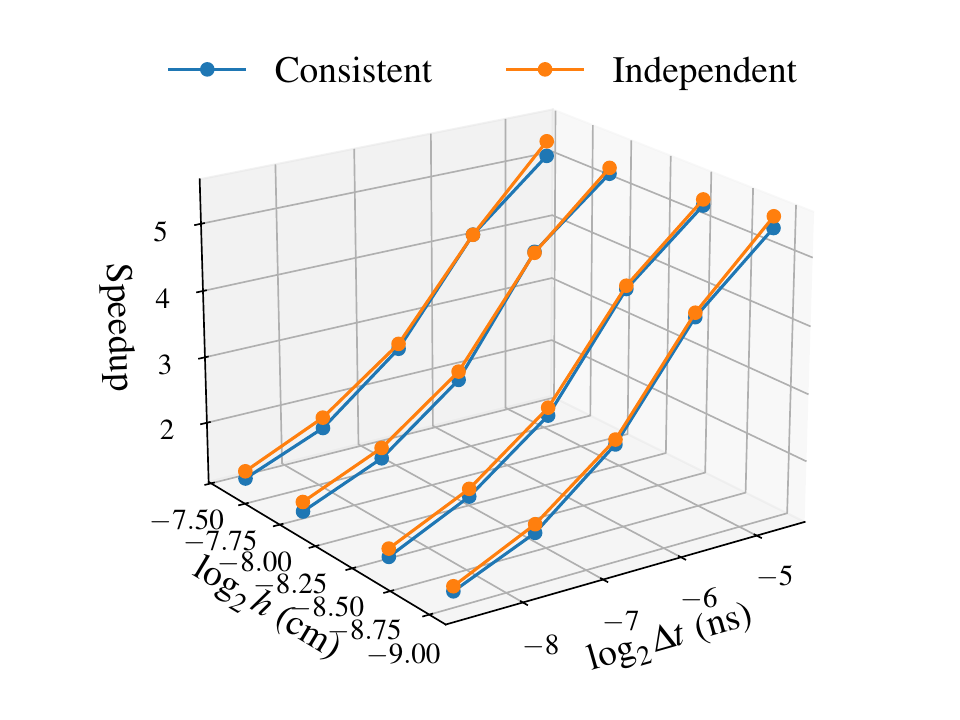}
		\caption{Wall Time}
	\end{subfigure}
	\caption{Speedup in terms of (a) sweeps and (b) wall time of the two SM methods relative to the unaccelerated method as a function of mesh and time step size. }
	\label{fig:larsen_speedup}
\end{figure}
Figure \ref{fig:larsen_speedup} shows the speedup, in terms of total number of sweeps performed across the entire simulation and wall time, of the SM methods with respect to the unaccelerated method. 
On the smallest time steps, where the unaccelerated method is most efficient, consistent and independent were $1.17\times$ and $1.28\times$ faster than the unaccelerated method. 
As the time step size increases, the SM methods remain rapidly convergent while the unaccelerated method degrades, resulting in speedups up to $5.39\times$ and $5.56\times$ for consistent and independent, respectively. 
The speedup in terms of sweeps is higher with consistent and independent requiring up to $9.5\times$ and $9.6\times$ fewer sweeps than the unaccelerated method.

\begin{table}
	\centering
	\caption{Summary of iterative performance and cost of the sweep, low-order nonlinear solve, low-order linear solve, and nonlinear temperature elimination for the consistent SM, independent SM, and unaccelerated methods on Larsen's problem. Here, ``coarse'' and ``fine'' refer to simulations corresponding to the least and most resolved mesh and time step sizes in the space-time convergence study, respectively.}
	\label{tab:summary_larsen}
	\begin{tabular}{cccccccccc}
\toprule
 &  & \multicolumn{3}{c}{Coarse}  &  & \multicolumn{3}{c}{Fine} \\
\cmidrule{3-5}\cmidrule{7-9}
 & Metric & Consistent & Independent & Unaccelerated & & Consistent & Independent & Unaccelerated \\
\midrule
\multirow{6}{*}{\rotatebox{90}{Sweeps}} & Min & 5 & 5 & 41 & & 3 & 3 & 3 \\
 & Max & $100^*$ & $100^*$ & $500^*$ & & 12 & 12 & 37 \\
 & Avg. & 15.65 & 15.62 & 116.34 & & 6.89 & 6.44 & 15.49 \\
 & Total & \num{4929} & \num{4921} & \num{36646} & & \num{68879} & \num{64399} & \num{154918} \\
 & Cost (s) & \num{142.1} & \num{141.1} & \num{925.6} & & \num{4704.3} & \num{4453.2} & \num{9186.7} \\
 & Avg. Cost (ms) & \num{28.8} & \num{28.7} & \num{25.3} & & \num{68.3} & \num{69.2} & \num{59.3} \\
\addlinespace
\multirow{6}{*}{\rotatebox{90}{LO Nonlinear}} & Min & 2 & 2 & -- & & 2 & 2 & -- \\
 & Max & 12 & 11 & -- & & 4 & 4 & -- \\
 & Avg. & 2.24 & 2.25 & -- & & 2.03 & 2.03 & -- \\
 & Total & \num{13119} & \num{13432} & -- & & \num{139398} & \num{130439} & -- \\
 & Cost (s) & \num{209.2} & \num{201.2} & -- & & \num{4274.0} & \num{4065.9} & -- \\
 & Avg. Cost (ms) & \num{42.4} & \num{40.9} & -- & & \num{62.1} & \num{63.1} & -- \\
\addlinespace
\multirow{6}{*}{\rotatebox{90}{LO Linear}} & Min & 4 & 4 & -- & & 5 & 5 & -- \\
 & Max & 21 & 20 & -- & & 17 & 17 & -- \\
 & Avg. & 15.66 & 15.67 & -- & & 14.59 & 14.69 & -- \\
 & Total & \num{205406} & \num{210473} & -- & & \num{2033811} & \num{1915974} & -- \\
 & Cost (s) & \num{51.8} & \num{53.0} & -- & & \num{985.3} & \num{932.2} & -- \\
 & Avg. Cost (ms) & \num{3.9} & \num{3.9} & -- & & \num{7.1} & \num{7.1} & -- \\
\addlinespace
\multirow{6}{*}{\rotatebox{90}{Energy Balance}} & Min & 1 & 1 & 2 & & 1 & 1 & 1 \\
 & Max & $40^*$ & $40^*$ & 7 & & 5 & 5 & 4 \\
 & Avg. & 7.01 & 6.97 & 3.31 & & 2.99 & 3.02 & 3.48 \\
 & Total & \num{92001} & \num{93593} & \num{121293} & & \num{417063} & \num{393862} & \num{538855} \\
 & Cost (s) & \num{51.9} & \num{49.4} & \num{745.3} & & \num{733.1} & \num{737.5} & \num{3828.6} \\
 & Avg. Cost (ms) & \num{4.0} & \num{3.7} & \num{20.3} & & \num{5.3} & \num{5.7} & \num{24.7} \\
\addlinespace
 & Wall Time (s) & \num{472.4} & \num{462.6} & \num{1940.5} & & \num{13131.2} & \num{12440.0} & \num{17378.1} \\
\bottomrule
\end{tabular}
	$^*$ indicates maximum iterations reached
\end{table}
Table \ref{tab:summary_larsen} provides a summary of the iterative performance and cost of the primary components of the algorithm for the coarsest and finest resolutions in space and time for the SM and unaccelerated methods. 
All three methods exhibit difficulties converging the outermost iteration on the coarsest resolution. 
In particular, all three methods fail to converge within their respective maximum allowed iterations for the first 5 time steps and another 10 time steps near where high frequency radiation heats the inner thick material causing the temperature spike seen in Fig.~\ref{fig:larsen_opacity_a}. 
Outside of these 15 outlying time steps, the reported average values represent the behavior of the algorithms well. 
Unlike the unaccelerated method, the SM methods also struggle to converge the nonlinear temperature elimination iteration in the first five time steps. 
This may be due to negativities in the low-order solution causing non-convergence in the local Newton solver.
Note that the outer nonlinear low-order iteration still converges in at most twelve iterations even when the temperature elimination solve fails to converge. 
Convergence of the nonlinear elimination solve is measured locally on each element and thus it is likely that only a few elements are struggling to converge. 
These iterative issues are resolved on the fine problem with all methods converging the outer, inner, and temperature elimination iterations across all time steps. 
On both resolutions, the low-order nonlinear and linear solves are robust converging in 2-3 and 14-16 iterations on average, respectively. 

\subsection{Crooked Pipe}
The crooked pipe problem is a two-dimensional, multi-material, gray problem with constant opacities. 
The early-time dynamics are challenging to resolve in space and angle, resulting in negativity-inducing ray effects in the thin pipe region. 
At long times, the domain has heated but the temperature-independent opacities remain large, inducing the equilibrium diffusion limit, a numerically challenging regime to resolve. 
Overall, this problem stresses the robustness of the algorithms, particularly their response to negativity, and their ability to remain iteratively effective in an absorption-emission dominated, strongly diffusive equilibrium limit. 
This problem is challenging to resolve in space and requires a long final time to reach the numerically interesting diffusive regime. 
Both of these factors make obtaining a sufficiently resolved, space-time accurate discrete reference challenging. 
Thus, we focus on comparing the methods in terms of their solution quality and performance on representative spatially coarse and fine resolutions.
This section concludes with a comparison of the SM methods' sensitivity to angular resolution in terms of solution quality and performance. 

\subsubsection{Definition}
\begin{figure}
	\centering
	\includegraphics[width=0.65\textwidth]{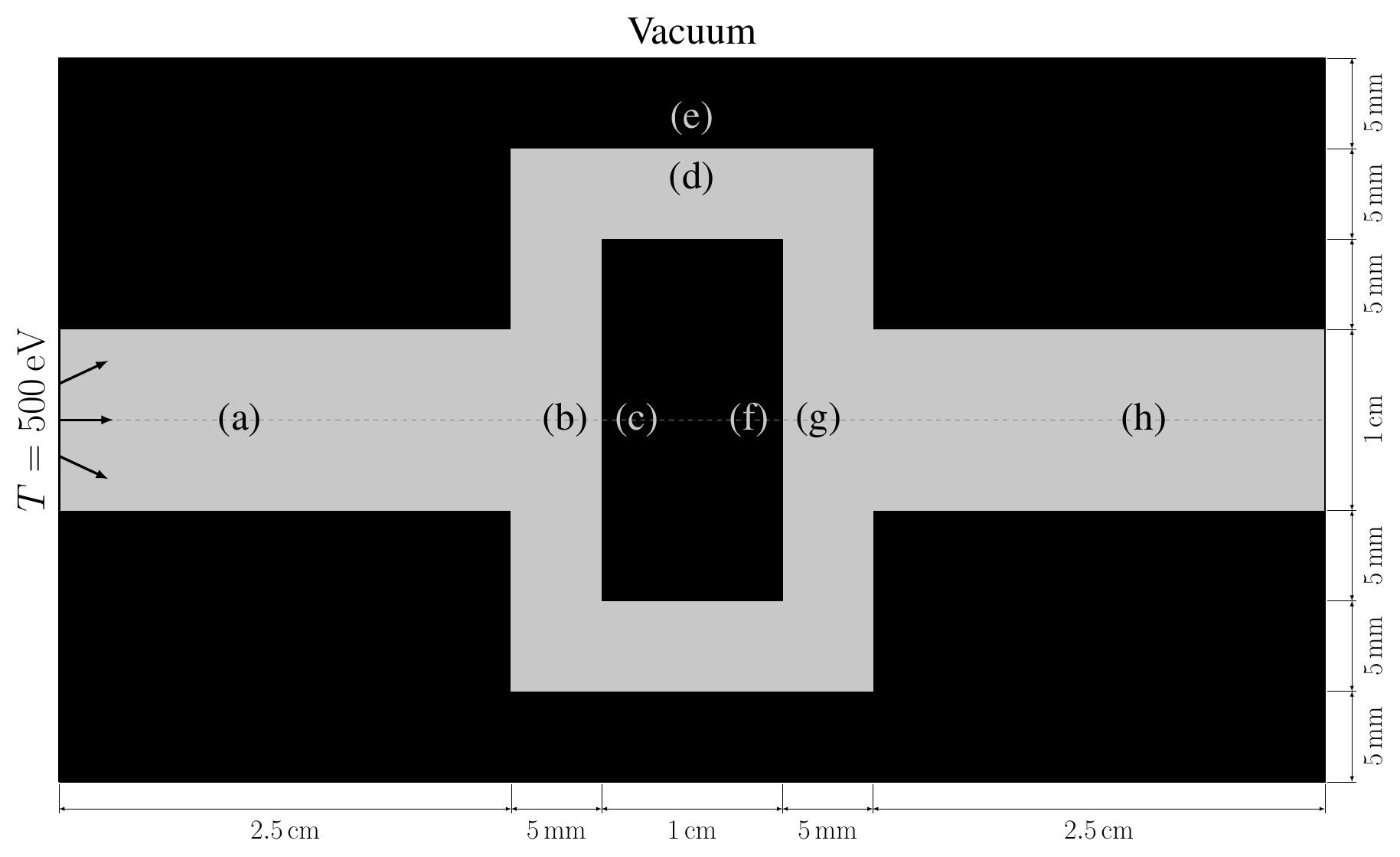}
	\caption{Problem description for the gray crooked pipe problem. The gray and black regions represent the pipe and wall, respectively. A reflection plane along $y=0$ is used to halve the computational domain. The solutions at the spatial locations labeled (a)--(h) are recorded at each time step. }
	\label{fig:cp_diag}
\end{figure}
\begin{figure}
	\centering
	\includegraphics[width=0.6\textwidth,trim={0 0 0 0}, clip]{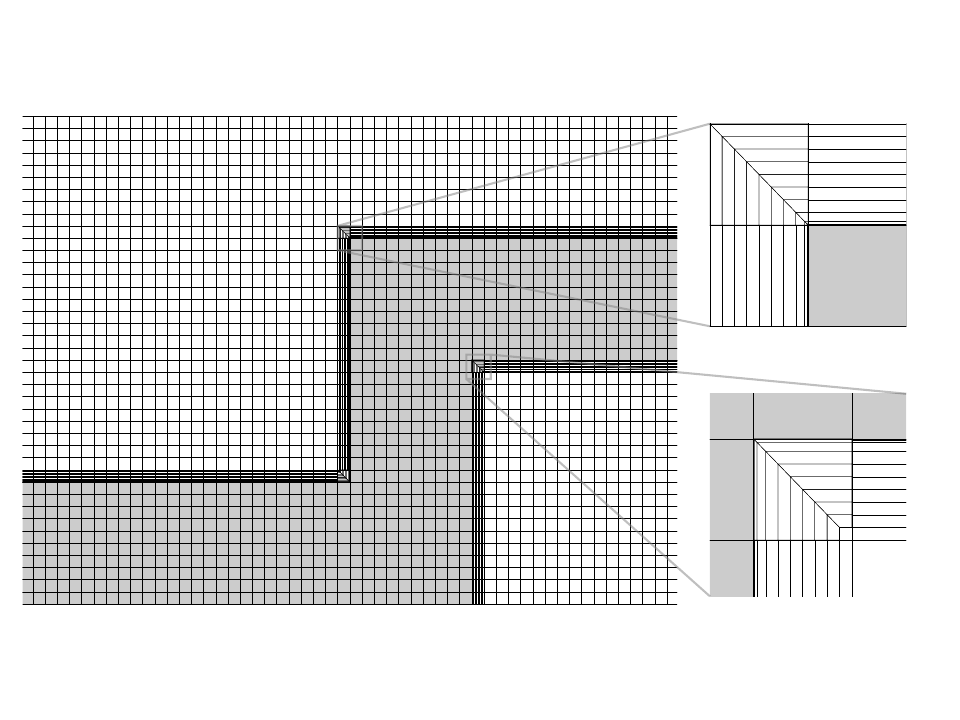}
	\caption{A depiction of the refinement pattern used to resolve the crooked pipe problem in space. Conformal refinement is applied along the wall side of all pipe-wall interfaces.}
	\label{fig:cp_mesh}
\end{figure}

The geometry and materials are depicted in Fig.~\ref{fig:cp_diag}. 
The problem consists of a thick and a thin material, denoted the ``wall'' and ``pipe'', respectively, and is driven by a \SI{500}{\eV} boundary source at the left inlet of the pipe. 
The thick and thin opacities are $\sigma_\text{thick} = \SI{2000}{\per\cm}$ and $\sigma_\text{thin} = \SI{0.2}{\per\cm}$, respectively. 
We use a constant time step of size \SI{1e-2}{\sh} where $\SI{1}{\sh} = \SI{10}{\ns}$. 
The final time is \SI{50}{\sh}.  
We consider two spatial meshes: a uniformly spaced mesh with $140 \times 40$ elements and a fine mesh with a base of $252 \times 72$ elements and conforming refinement along the wall side of all material interfaces. 
The fine mesh uses logarithmic refinement that feathers into uniform refinement in the direction perpendicular to material interfaces. 
This mix of logarithmic and uniform refinement helps resolve interface-induced boundary layers as radiation penetrates into the wall region. 
A pattern of trapezoidal and square elements is used at corners to turn the refinement in tandem with the geometry of the pipe without creating hanging nodes. 
A portion of the refined mesh is shown in Fig.~\ref{fig:cp_mesh}. 
The coarse mesh has a characteristic mesh length of \SI{5e-2}{\cm}. 
The fine mesh has minimum and maximum characteristic mesh lengths of \SI{8.7841e-5}{\cm} and \SI{2.777e-2}{\cm}, respectively, and \num{22504} elements. 
The fine mesh then has between $1.8\times$ and $570\times$ more resolution than the coarse mesh and roughly $4\times$ the number of elements. 
Unless otherwise noted, we use an $S_{6}$ square Chebyshev-Legendre angular quadrature rule with 36 angles.
The coarse and fine problems are solved in parallel on one and four nodes, respectively, where each node has 44 processors. 
The zero and scale negative flux fixup \cite{hamilton2009negative} is used inside the sweep to ensure the high-order solution is positive while preserving particle balance. 
When present, negative temperatures are floored to a minimum temperature of \SI{1e-8}{\eV}. 

\subsubsection{Solution Quality} \label{sec:gcp_solution_quality}
\begin{figure}
	\centering
	\begin{subfigure}{0.4\textwidth}
		\centering
		\begin{overpic}[width=\textwidth, trim={0 0 875 0}, clip]{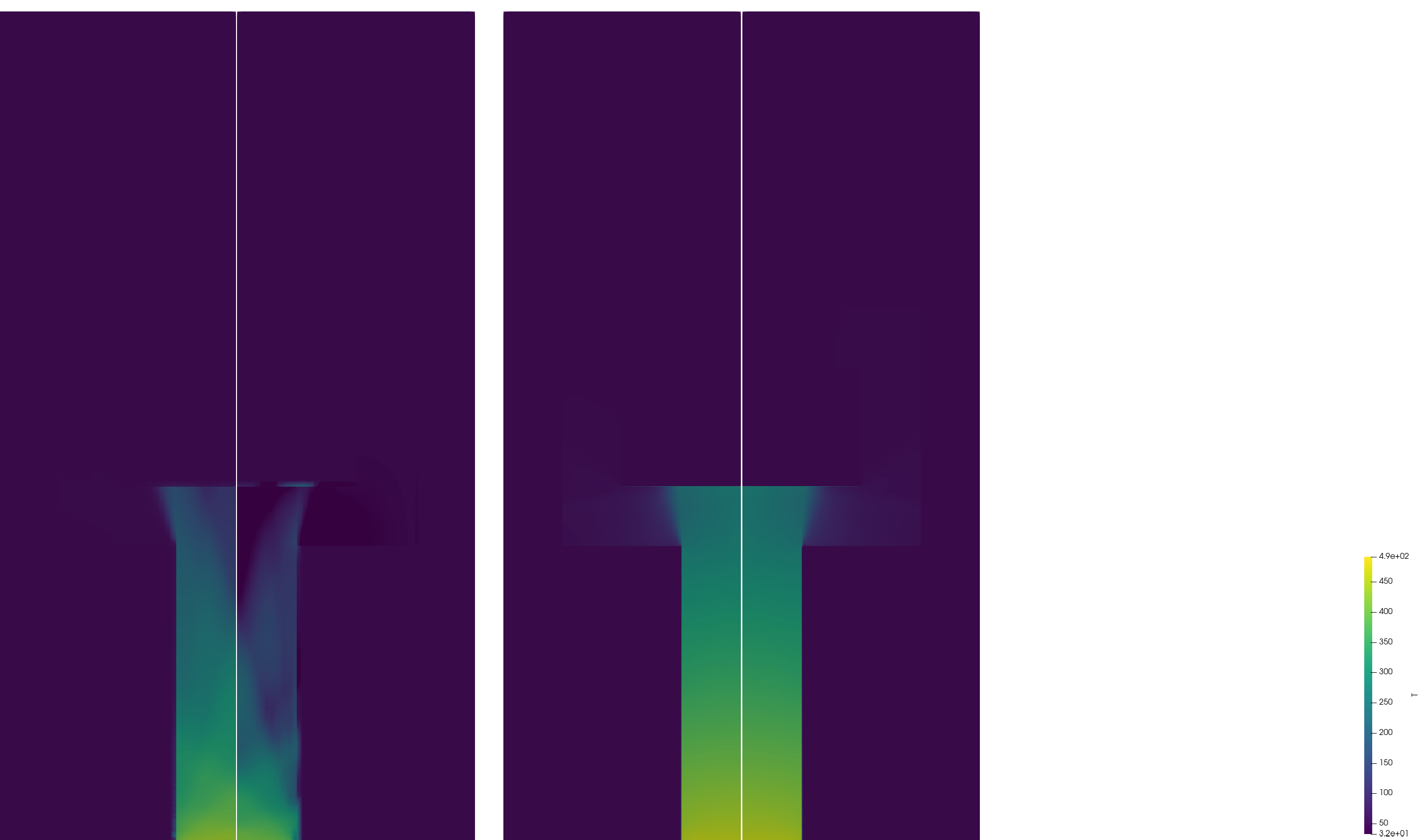}
			\put(4,80){\textcolor{white}{\scalebox{0.6}{Consistent}}}
			\put(28,80){\textcolor{white}{\scalebox{0.6}{Independent}}}
			\put(24,86){\makebox(0,0){\scalebox{0.6}{Coarse}}}

			\put(55,80){\textcolor{white}{\scalebox{0.6}{Consistent}}}
			\put(79,80){\textcolor{white}{\scalebox{0.6}{Independent}}}
			\put(75,86){\makebox(0,0){\scalebox{0.6}{Fine}}}
		\end{overpic}
		\caption{$t = \SI{0.1}{\sh}$}
	\end{subfigure}
	\quad
	\begin{subfigure}{0.4\textwidth}
		\centering
		\begin{overpic}[width=\textwidth, trim={0 0 875 0}, clip]{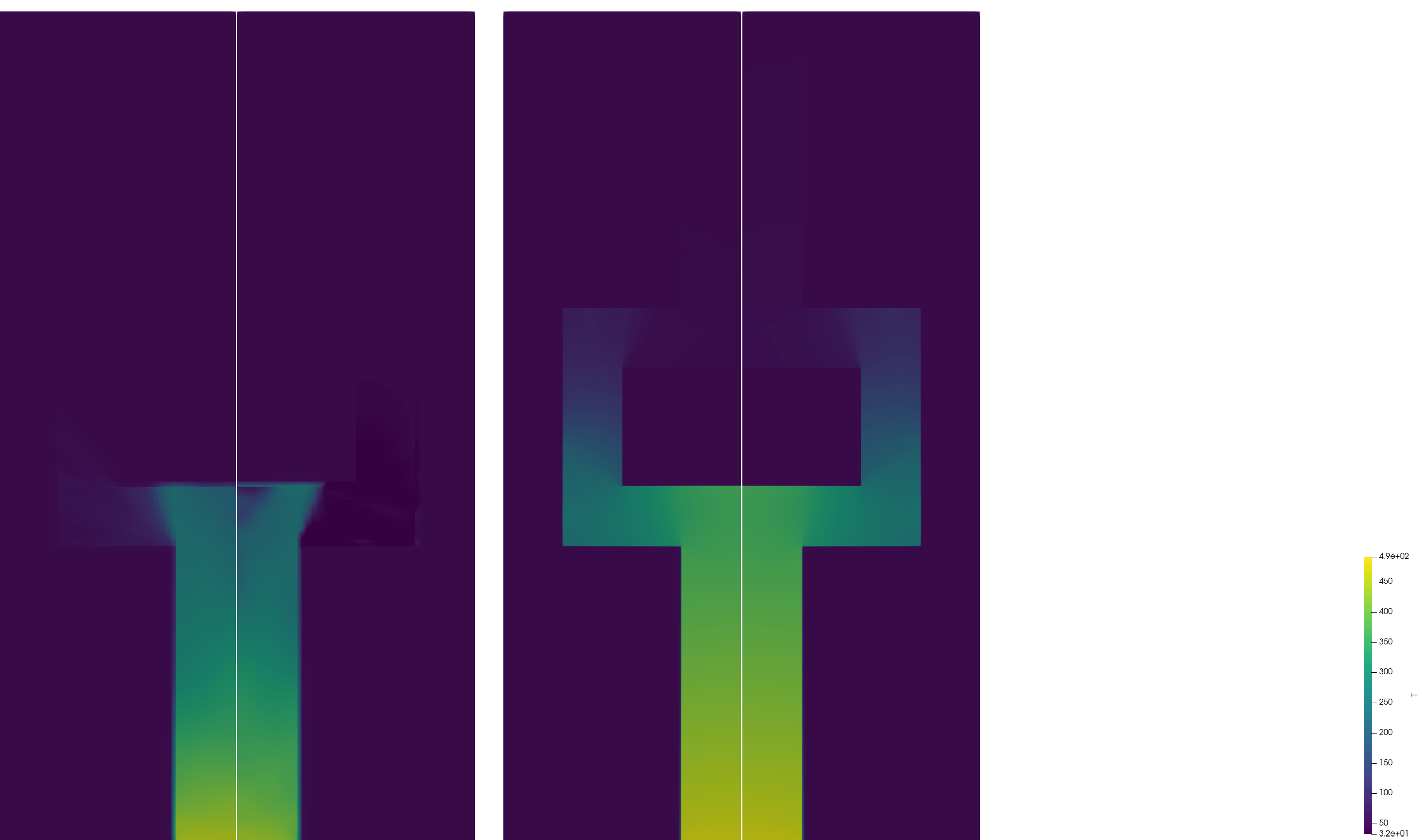}
			\put(4,80){\textcolor{white}{\scalebox{0.6}{Consistent}}}
			\put(28,80){\textcolor{white}{\scalebox{0.6}{Independent}}}
			\put(24,86){\makebox(0,0){\scalebox{0.6}{Coarse}}}

			\put(55,80){\textcolor{white}{\scalebox{0.6}{Consistent}}}
			\put(79,80){\textcolor{white}{\scalebox{0.6}{Independent}}}
			\put(75,86){\makebox(0,0){\scalebox{0.6}{Fine}}}
		\end{overpic}
		\caption{$t = \SI{0.5}{\sh}$}
	\end{subfigure}
	\begin{subfigure}{0.4\textwidth}
		\centering
		\begin{overpic}[width=\textwidth, trim={0 0 875 0}, clip]{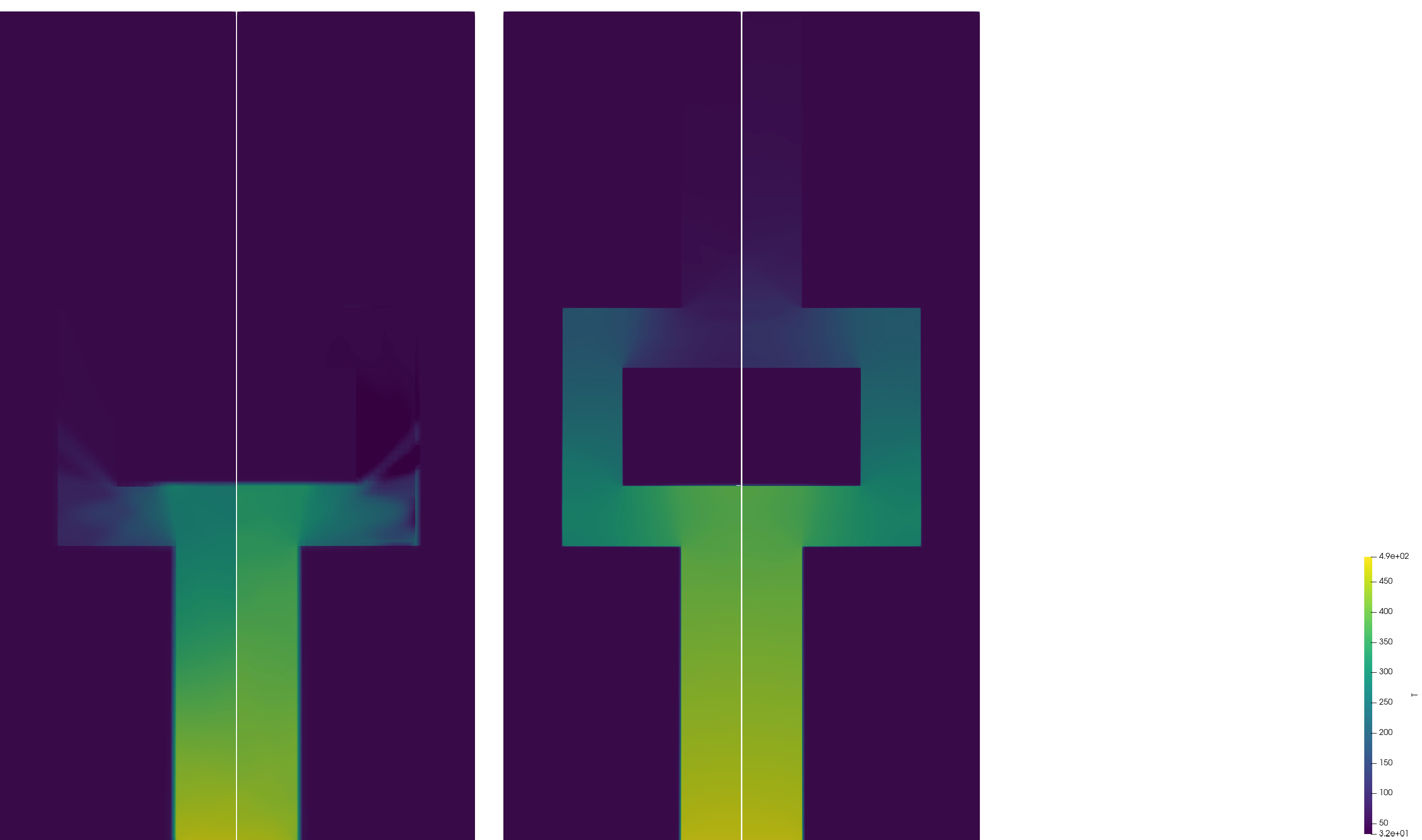}
			\put(4,80){\textcolor{white}{\scalebox{0.6}{Consistent}}}
			\put(28,80){\textcolor{white}{\scalebox{0.6}{Independent}}}
			\put(24,86){\makebox(0,0){\scalebox{0.6}{Coarse}}}

			\put(55,80){\textcolor{white}{\scalebox{0.6}{Consistent}}}
			\put(79,80){\textcolor{white}{\scalebox{0.6}{Independent}}}
			\put(75,86){\makebox(0,0){\scalebox{0.6}{Fine}}}
		\end{overpic}
		\caption{$t = \SI{1}{\sh}$}
	\end{subfigure}
	\quad
	\begin{subfigure}{0.4\textwidth}
		\centering
		\begin{overpic}[width=\textwidth, trim={0 0 875 0}, clip]{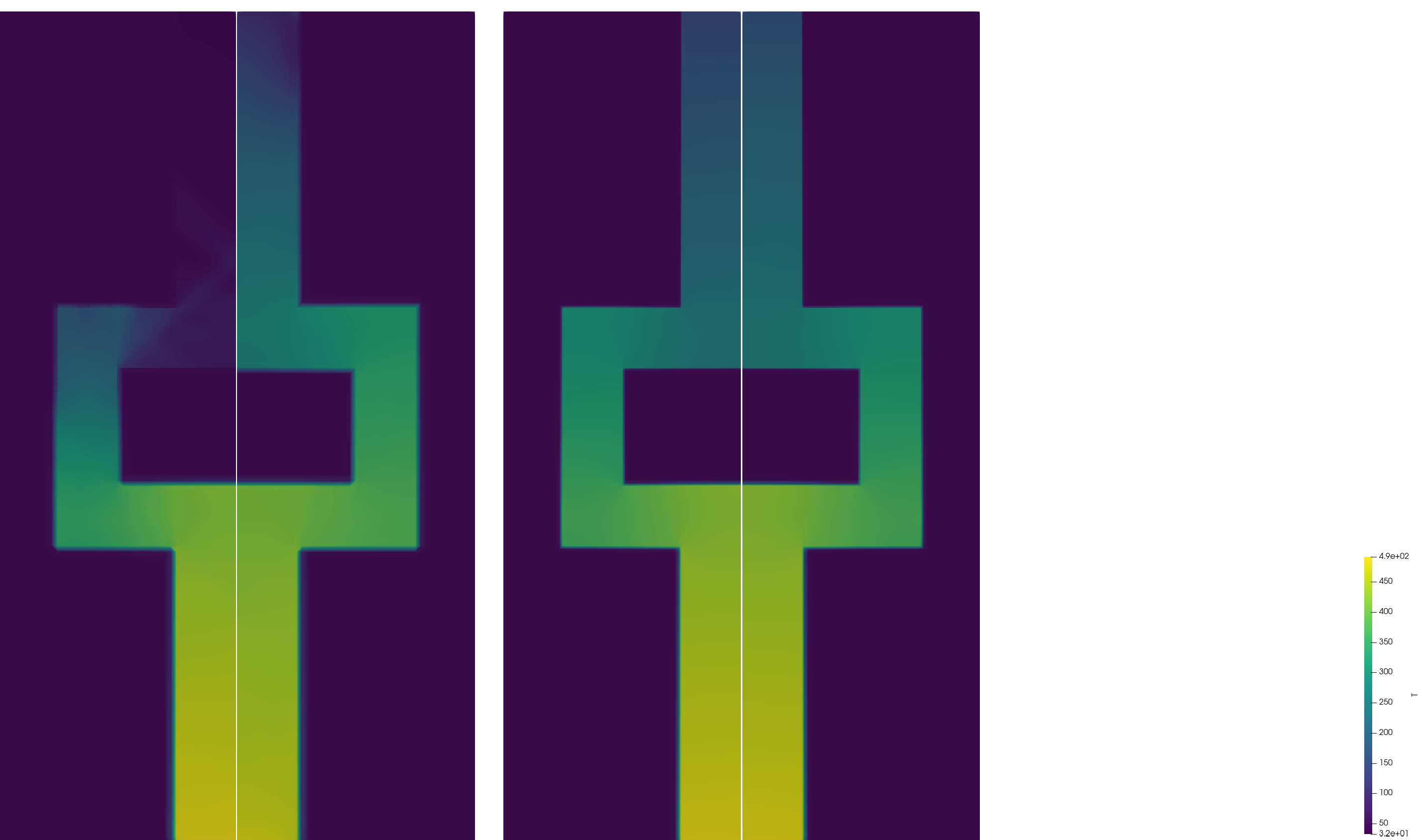}
			\put(4,80){\textcolor{white}{\scalebox{0.6}{Consistent}}}
			\put(28,80){\textcolor{white}{\scalebox{0.6}{Independent}}}
			\put(24,86){\makebox(0,0){\scalebox{0.6}{Coarse}}}

			\put(55,80){\textcolor{white}{\scalebox{0.6}{Consistent}}}
			\put(79,80){\textcolor{white}{\scalebox{0.6}{Independent}}}
			\put(75,86){\makebox(0,0){\scalebox{0.6}{Fine}}}
		\end{overpic}
		\caption{$t = \SI{5}{\sh}$}
	\end{subfigure}
	\begin{subfigure}{0.4\textwidth}
		\centering
		\begin{overpic}[width=\textwidth, trim={0 0 875 0}, clip]{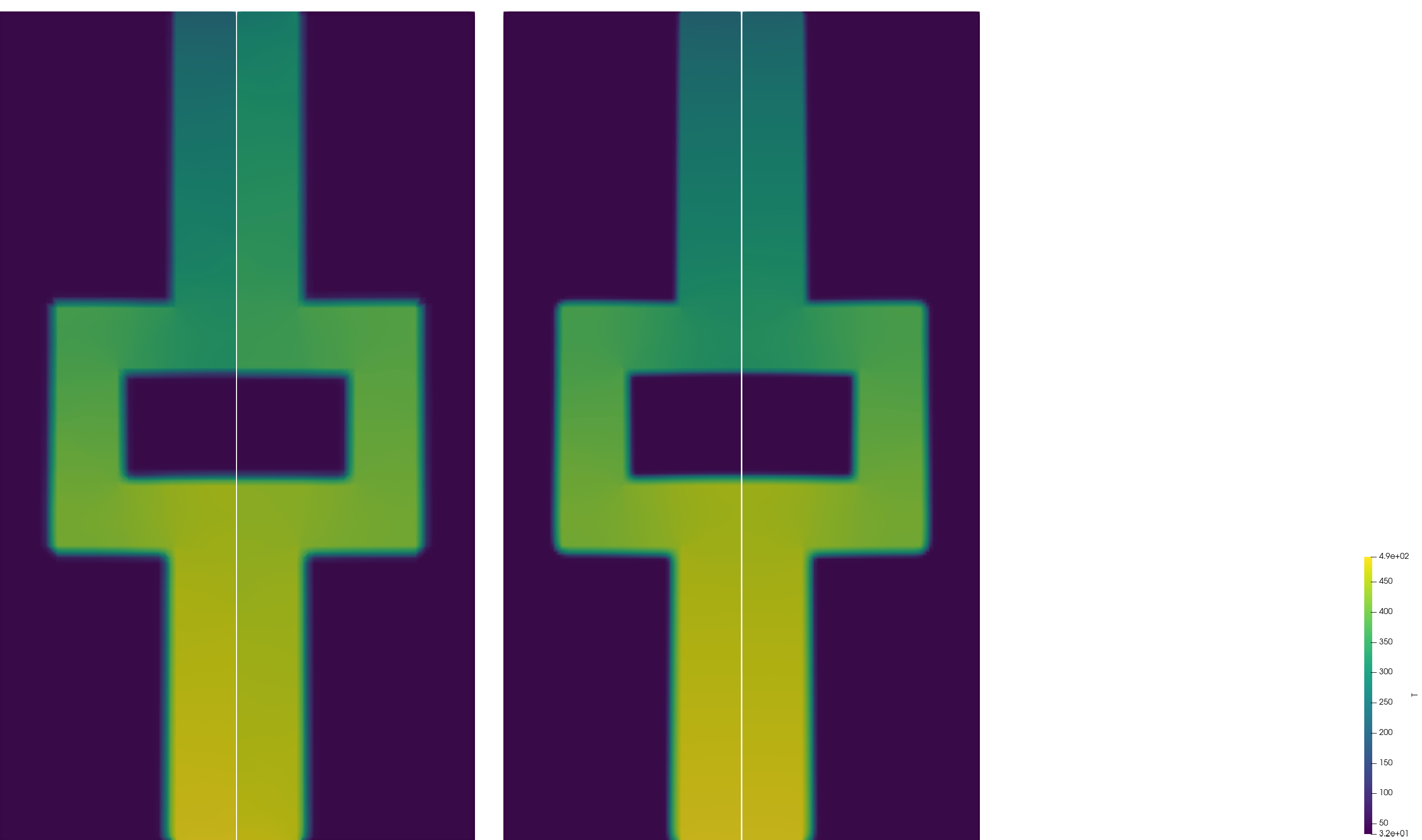}
			\put(4,80){\textcolor{white}{\scalebox{0.6}{Consistent}}}
			\put(28,80){\textcolor{white}{\scalebox{0.6}{Independent}}}
			\put(24,86){\makebox(0,0){\scalebox{0.6}{Coarse}}}

			\put(55,80){\textcolor{white}{\scalebox{0.6}{Consistent}}}
			\put(79,80){\textcolor{white}{\scalebox{0.6}{Independent}}}
			\put(75,86){\makebox(0,0){\scalebox{0.6}{Fine}}}
		\end{overpic}
		\caption{$t = \SI{50}{\sh}$}
	\end{subfigure}
	\caption{Temperature profiles over time from the consistent and independent methods on the coarse and fine crooked pipe problems. }
	\label{fig:gcp_temperature}
\end{figure}
Figure \ref{fig:gcp_temperature} shows the evolution of the temperature computed with the consistent and independent methods on the coarse and fine meshes. 
Note that, to save space, we omit comparison to the unaccelerated method as the consistent SM and unaccelerated methods produced solutions that were visually equivalent; detailed comparison is delayed until discussion of the tracer data. 
On the coarse mesh, the independent method produces negative temperatures at early times resulting in non-physical cold spots at the wavefront. 
While a negative flux fixup guarantees the high-order solution is positive, the low-order solution can be negative on under-resolved problems. 
For the independent method, the negativities in the low-order energy density are severe enough to induce negative temperatures, requiring flooring to keep the simulation running. 
In contrast, the consistent method has positive temperatures on the same mesh without the need for flooring. 
This is likely due to the consistent low-order system being more accurate than the independent low-order system. 
The independent method also has a significantly faster wave speed; at $t = \SI{5}{\sh}$, the wavefront associated with the independent method has traveled the length of the pipe while the consistent method's solution is only at the final turn. 
On the fine mesh, temperature flooring was not required for either method. 
In addition, the consistent and independent solutions are qualitatively much closer. 
Note that the wave speeds computed on the coarse and fine meshes are very different suggesting that the uniform $140\times 40$ mesh is too coarse to produce an accurate solution. 
In other words, while independent did struggle with solution quality, requiring energy-conservation-violating floors to continue running, this reduced solution quality occurred on an unrealistically coarse mesh where the more accurate consistent method is also incorrect. 

\begin{figure}
	\centering
	\foreach \x in {0,4,5,1,2,6,7,3}{
		\begin{subfigure}{0.24\textwidth}
			\centering
			\includegraphics[width=\textwidth]{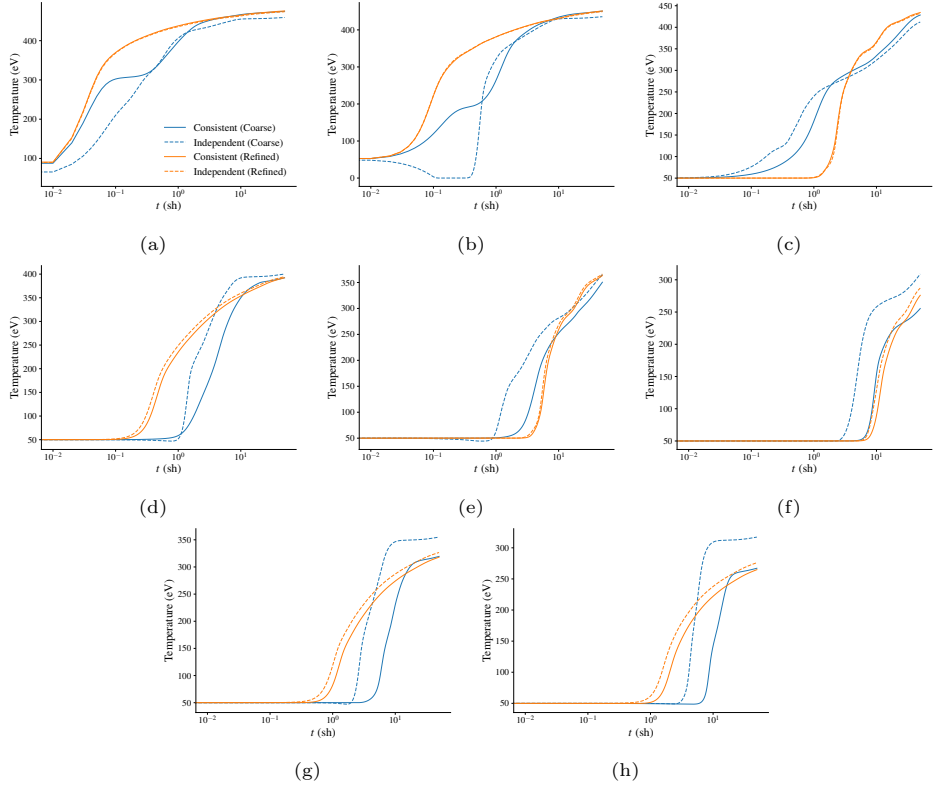}
			\caption{}
		\end{subfigure}	
	}
	\caption{Evolution of the temperature over time for the consistent and independent methods on the coarse and fine crooked pipe problems. Subfigures (a)--(h) correspond to the spatial locations labeled in Fig.~\ref{fig:cp_diag}. }
	\label{fig:gcp_tracer}
\end{figure}
The wave speeds of the consistent SM, independent SM, and unaccelerated methods are compared in detail using tracer locations where the solution is recorded at a single point in space at every time step. 
The evolution of the temperature at the eight spatial locations labeled in Fig.~\ref{fig:cp_diag} is shown in Fig.~\ref{fig:gcp_tracer}. 
Location (b), in particular, shows the effect of temperature flooring. 
The coarse independent simulation dips down to the minimum temperature for $t\in [\SI{0.1}{\sh}, \SI{0.5}{\sh}]$. 
On the other hand, the coarse consistent and unaccelerated simulations exhibit more physical behavior, resulting in smoother filling of the pipe. 
At tracer locations (d)--(h), coarse independent heats 5-\SI{10}{\sh} earlier than coarse consistent. 
At all locations, coarse consistent and coarse independent do not agree with each other nor the fine-mesh solutions, suggesting the coarse simulations are severely under-resolved in space and cannot accurately resolve the steep spatial gradients in the solution. 
The coarse independent simulation also equilibrates to incorrect final temperatures, failing to obtain the steady-state solution achieved by the fine-mesh simulations. 
On the coarse mesh, the largest deviations between the consistent and independent temperatures at the final time occur at tracer locations (f) and (h) with deviations of \SI{26.7}{\eV} and \SI{25}{\eV}, respectively. 
On the fine mesh, such deviations are reduced at all tracer locations with deviations of \SI{5.5}{\eV} and \SI{5.6}{\eV} at locations (f) and (h), respectively. 
This $\approx\!5\times$ reduction in deviation between consistent and independent suggests that the independent method can produce high quality solutions despite producing poor solution quality on the coarse mesh. 
In all locations, the unaccelerated method is closer to the consistent method than the independent method with closer agreement on the refined problem. 
These noticeable differences may be due to differences in where and when negative flux fixups are applied in addition to the relative stopping criteria differences explored in the previous section. 

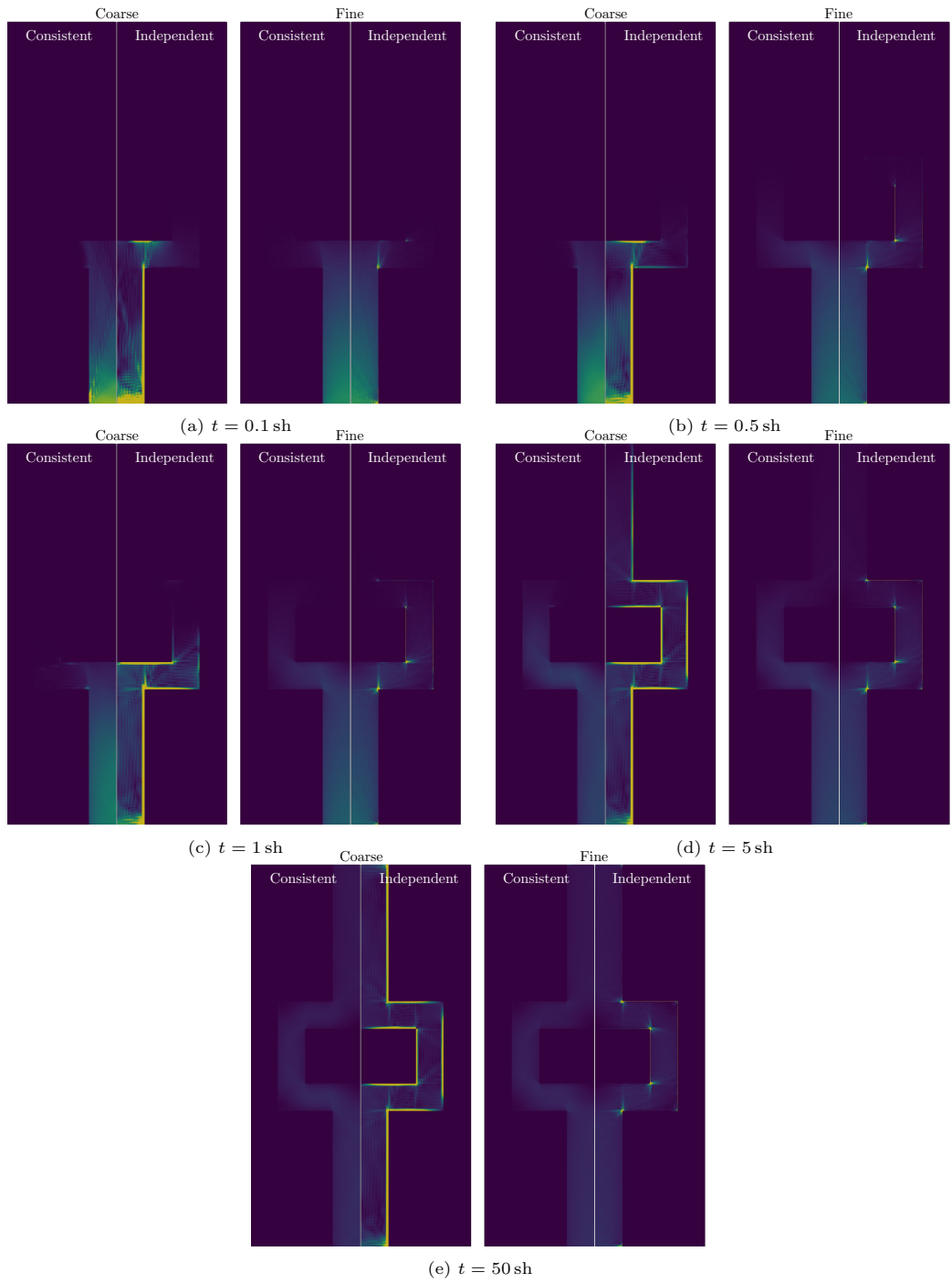
\begin{figure}
	\centering
	\begin{subfigure}{0.4\textwidth}
		\centering
		\begin{overpic}[width=\textwidth, trim={0 0 875 0}, clip]{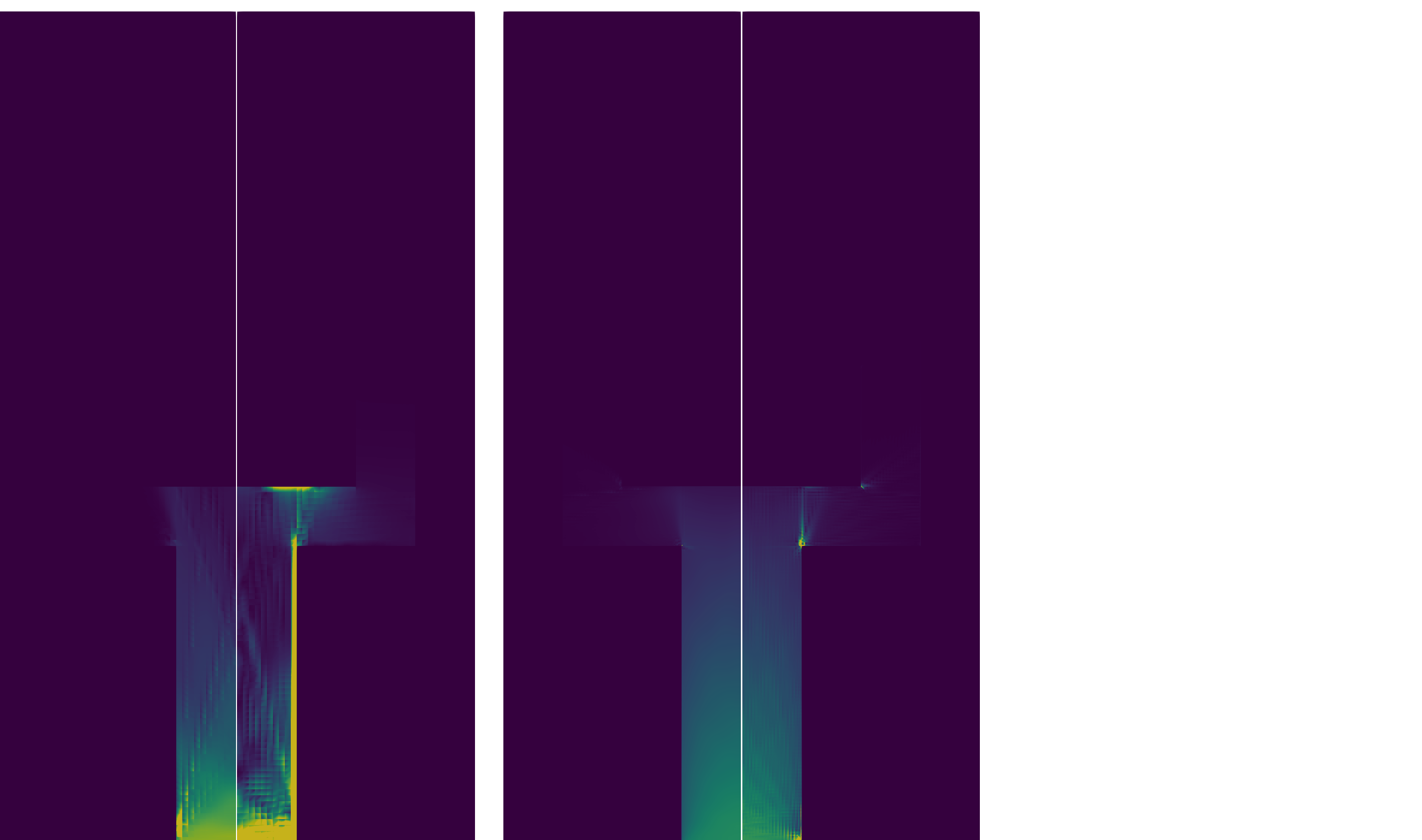}
			\put(4,80){\textcolor{white}{\scalebox{0.6}{Consistent}}}
			\put(28,80){\textcolor{white}{\scalebox{0.6}{Independent}}}
			\put(24,86){\makebox(0,0){\scalebox{0.6}{Coarse}}}

			\put(55,80){\textcolor{white}{\scalebox{0.6}{Consistent}}}
			\put(79,80){\textcolor{white}{\scalebox{0.6}{Independent}}}
			\put(75,86){\makebox(0,0){\scalebox{0.6}{Fine}}}
		\end{overpic}
		\caption{$t = \SI{0.1}{\sh}$}
	\end{subfigure}
	\quad
	\begin{subfigure}{0.4\textwidth}
		\centering
		\begin{overpic}[width=\textwidth, trim={0 0 875 0}, clip]{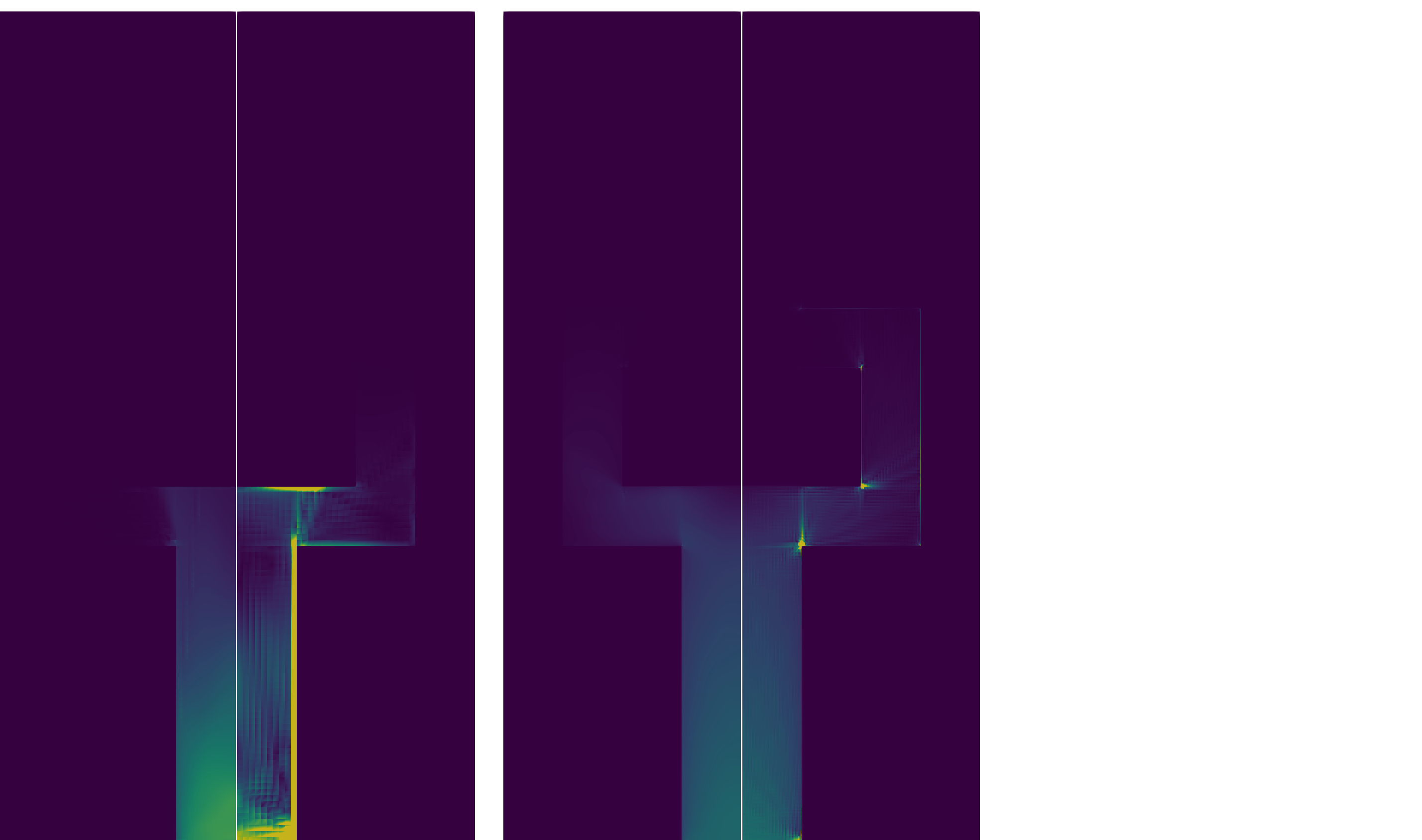}
			\put(4,80){\textcolor{white}{\scalebox{0.6}{Consistent}}}
			\put(28,80){\textcolor{white}{\scalebox{0.6}{Independent}}}
			\put(24,86){\makebox(0,0){\scalebox{0.6}{Coarse}}}

			\put(55,80){\textcolor{white}{\scalebox{0.6}{Consistent}}}
			\put(79,80){\textcolor{white}{\scalebox{0.6}{Independent}}}
			\put(75,86){\makebox(0,0){\scalebox{0.6}{Fine}}}
		\end{overpic}
		\caption{$t = \SI{0.5}{\sh}$}
	\end{subfigure}
	\begin{subfigure}{0.4\textwidth}
		\centering
		\begin{overpic}[width=\textwidth, trim={0 0 875 0}, clip]{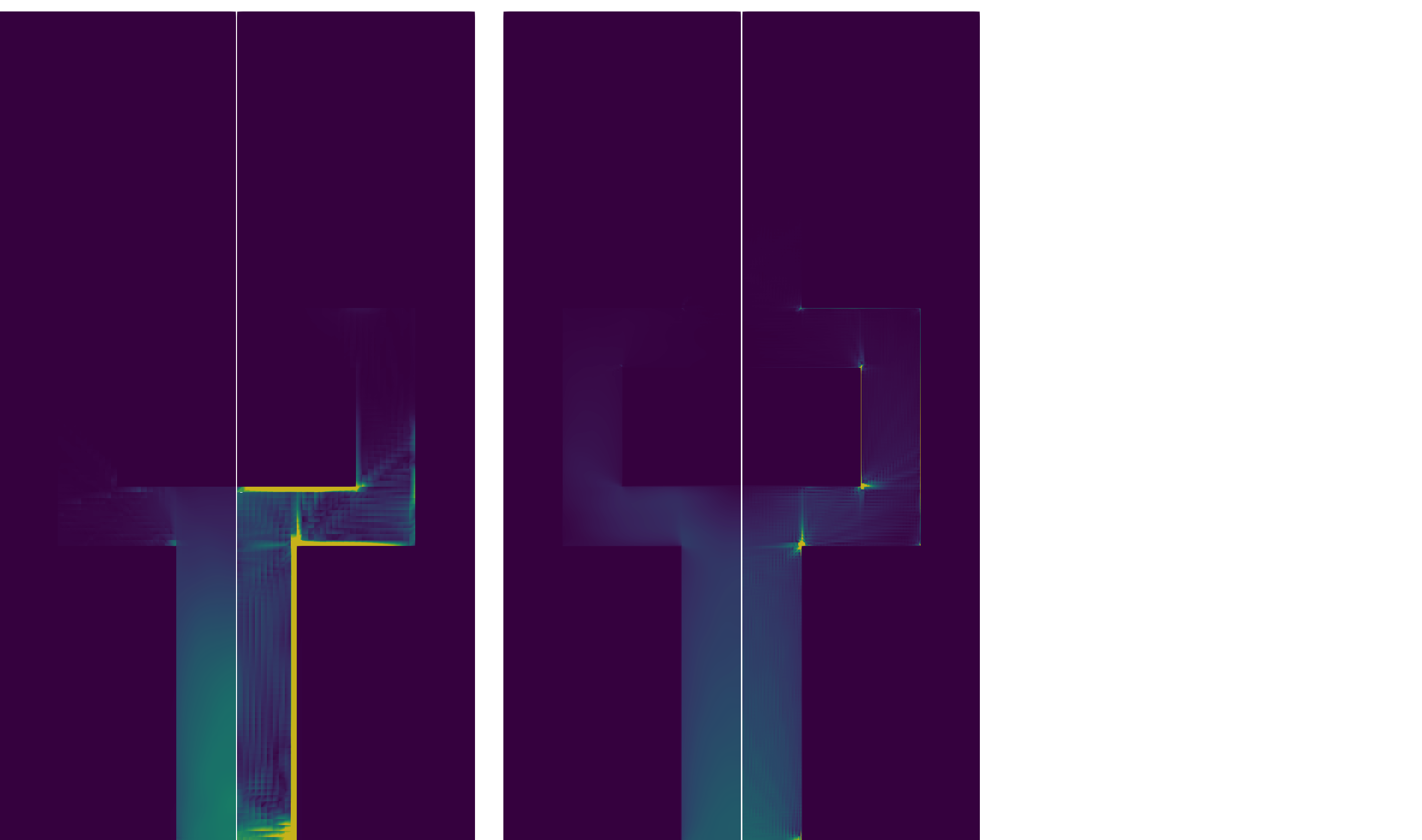}
			\put(4,80){\textcolor{white}{\scalebox{0.6}{Consistent}}}
			\put(28,80){\textcolor{white}{\scalebox{0.6}{Independent}}}
			\put(24,86){\makebox(0,0){\scalebox{0.6}{Coarse}}}

			\put(55,80){\textcolor{white}{\scalebox{0.6}{Consistent}}}
			\put(79,80){\textcolor{white}{\scalebox{0.6}{Independent}}}
			\put(75,86){\makebox(0,0){\scalebox{0.6}{Fine}}}
		\end{overpic}
		\caption{$t = \SI{1}{\sh}$}
	\end{subfigure}
	\quad
	\begin{subfigure}{0.4\textwidth}
		\centering
		\begin{overpic}[width=\textwidth, trim={0 0 875 0}, clip]{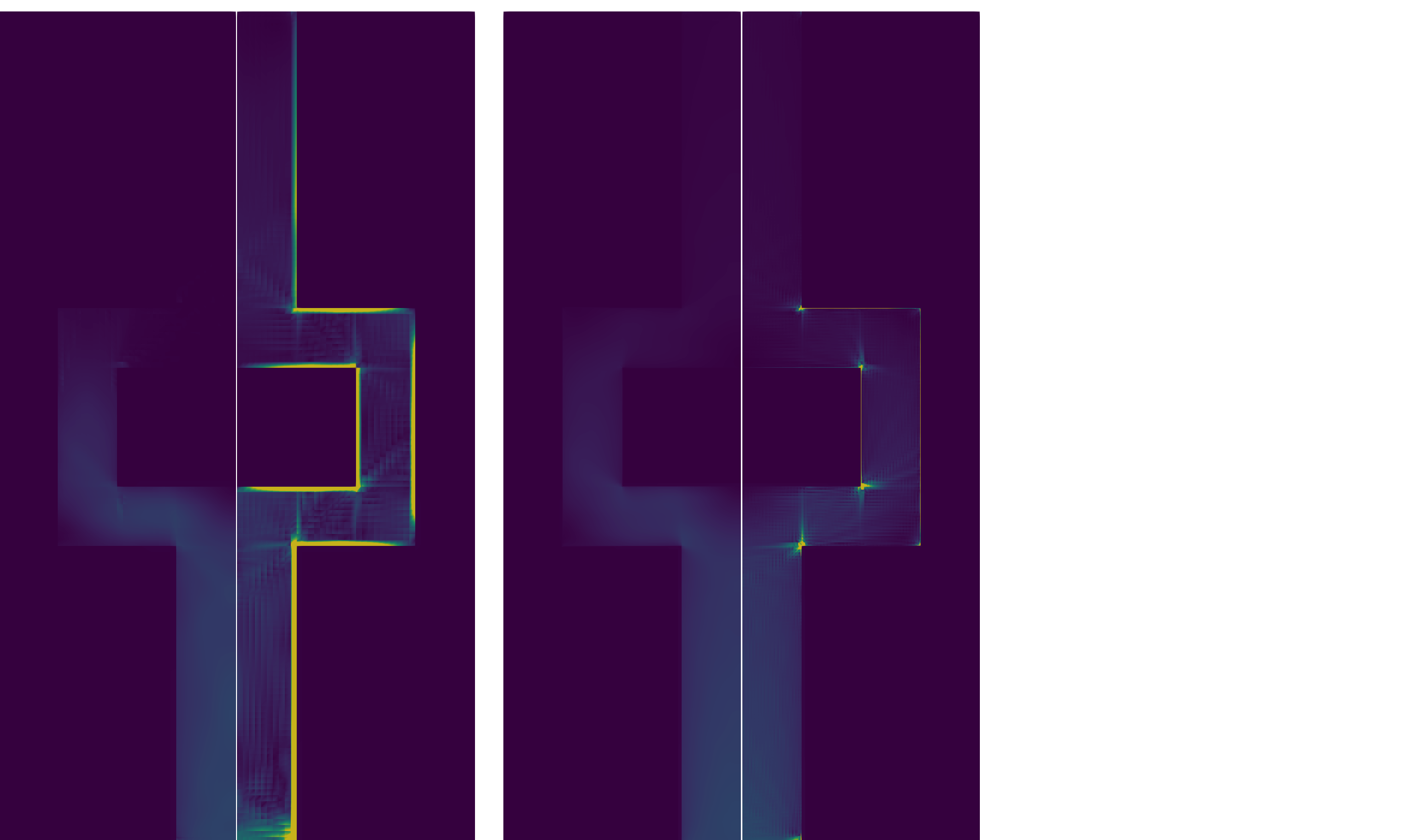}
			\put(4,80){\textcolor{white}{\scalebox{0.6}{Consistent}}}
			\put(28,80){\textcolor{white}{\scalebox{0.6}{Independent}}}
			\put(24,86){\makebox(0,0){\scalebox{0.6}{Coarse}}}

			\put(55,80){\textcolor{white}{\scalebox{0.6}{Consistent}}}
			\put(79,80){\textcolor{white}{\scalebox{0.6}{Independent}}}
			\put(75,86){\makebox(0,0){\scalebox{0.6}{Fine}}}
		\end{overpic}
		\caption{$t = \SI{5}{\sh}$}
	\end{subfigure}
	\begin{subfigure}{0.4\textwidth}
		\centering
		\begin{overpic}[width=\textwidth, trim={0 0 875 0}, clip]{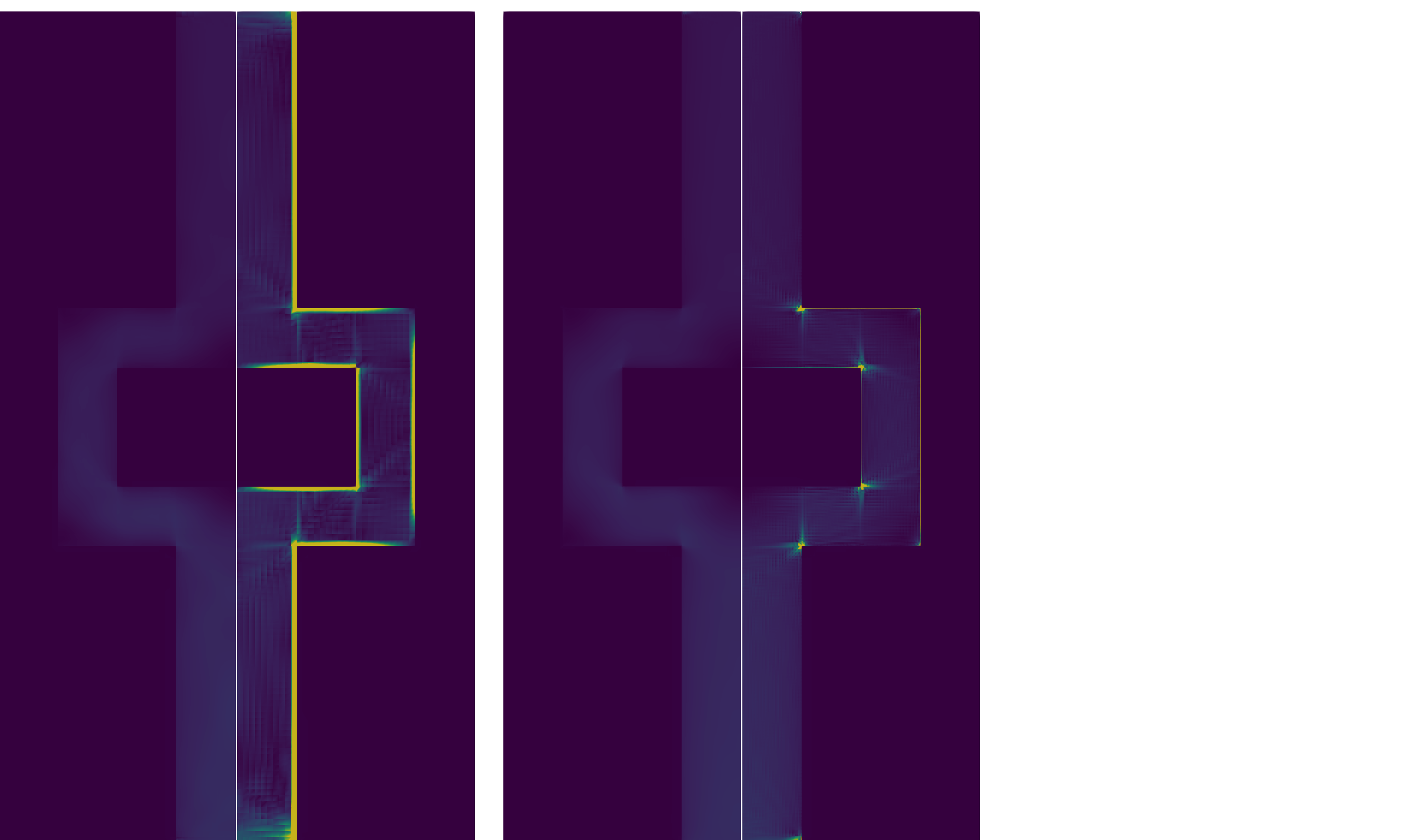}
			\put(4,80){\textcolor{white}{\scalebox{0.6}{Consistent}}}
			\put(28,80){\textcolor{white}{\scalebox{0.6}{Independent}}}
			\put(24,86){\makebox(0,0){\scalebox{0.6}{Coarse}}}

			\put(55,80){\textcolor{white}{\scalebox{0.6}{Consistent}}}
			\put(79,80){\textcolor{white}{\scalebox{0.6}{Independent}}}
			\put(75,86){\makebox(0,0){\scalebox{0.6}{Fine}}}
		\end{overpic}
		\caption{$t = \SI{50}{\sh}$}
	\end{subfigure}
	\caption{The magnitude of the flux over time from the consistent and independent methods on the coarse and fine crooked pipe problems. }
	\label{fig:gcp_flux}
\end{figure}
Next, we compare the solution quality of the low-order flux. 
Figure \ref{fig:gcp_flux} shows the magnitudes of the flux, $\|\vec{F}\|_2$, produced by the consistent and independent methods on the coarse and fine meshes. 
On the coarse mesh and at early times, the consistent method shows imprinting of the negative flux fixup causing non-smooth oscillations in the flux. 
However, these oscillations quickly diffuse out, resulting in a smooth flux by $t=\SI{0.5}{\sh}$. 
On the other hand, the independent method shows poor solution quality in the flux throughout the simulation. 
The flux appears to non-physically accumulate along the wall side of material interfaces. 
This might be an explanation for the much faster temperature wave speed on the coarse mesh. 
This behavior is reduced on the fine mesh but hot spots are still present at all corners of material interfaces. 
The consistent method on the fine mesh produces a physically reasonable flux. 
As observed on steady-state problems \cite{csmm}, the indepenent method converges the flux (first moment) very slowly in space and is prone to poor solution quality on multi-dimensional problems with heterogeneous materials. 

\begin{figure}
	\centering
	\begin{subfigure}{0.4\textwidth}
		\centering 
		\includegraphics[width=\textwidth]{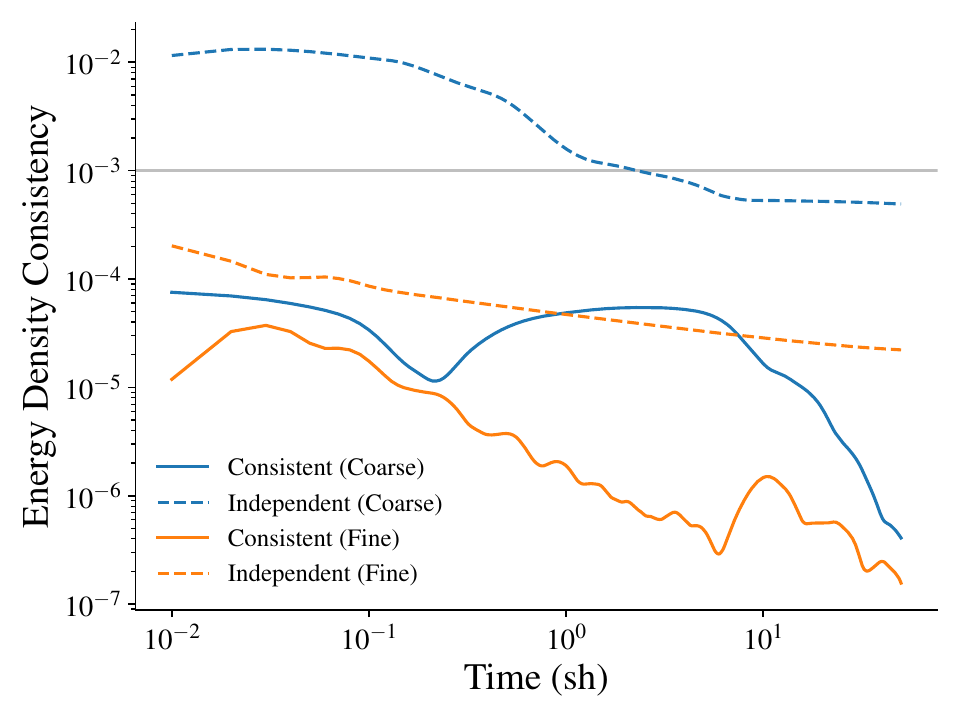}
		\caption{}
	\end{subfigure}
	\begin{subfigure}{0.4\textwidth}
		\centering 
		\includegraphics[width=\textwidth]{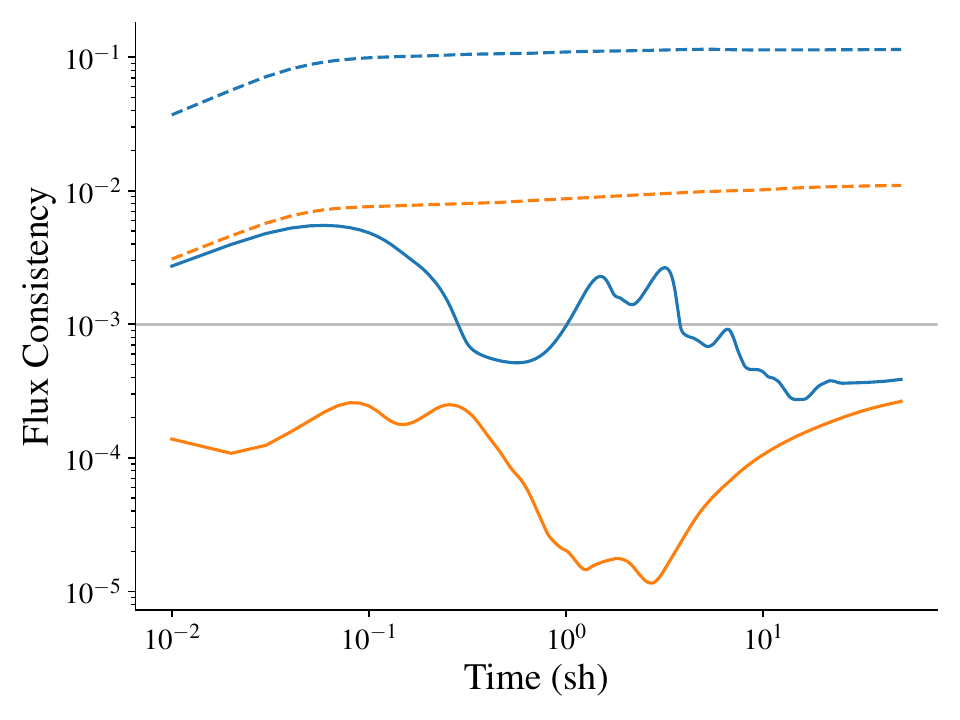}
		\caption{}
	\end{subfigure}
	\caption{The difference between the high and low-order energy density and flux as a function of simulation time for the consistent and independent SM methods on the coarse and fine crooked pipe problems.}
	\label{fig:gcp_consistency}
\end{figure}
The difference between the high and low-order energy density and flux is plotted over simulation time in Fig.~\ref{fig:gcp_consistency}. 
The horizontal line indicates the relative, residual-based tolerance of $10^{-3}$ used on this problem. 
For the energy density, the consistent method achieves agreement with the high-order energy density below the iterative tolerance at all simulation times on both the coarse and fine meshes. 
Consistency in the flux is below the tolerance only on the fine mesh likely due to the increased reliance on the negative flux fixup on the coarse mesh. 
The independent method has between $5\times$-$100\times$ less consistency in the energy density than the consistent method on both the coarse and fine meshes. 
This disparity is larger in the flux where the consistent flux was up to $1000\times$ more consistent than the independent method. 
For independent methods, consistency between high and low-order solutions is an indicator of mesh resolution. 
Lack of consistency is then yet another indicator that this problem is very challenging to resolve in space, with even the fine independent simulation exhibiting signs of lack of resolution in the flux. 

\subsubsection{Performance}
\begin{figure} 
	\centering
	\begin{subfigure}{0.4\textwidth}
		\centering
		\includegraphics[width=\textwidth]{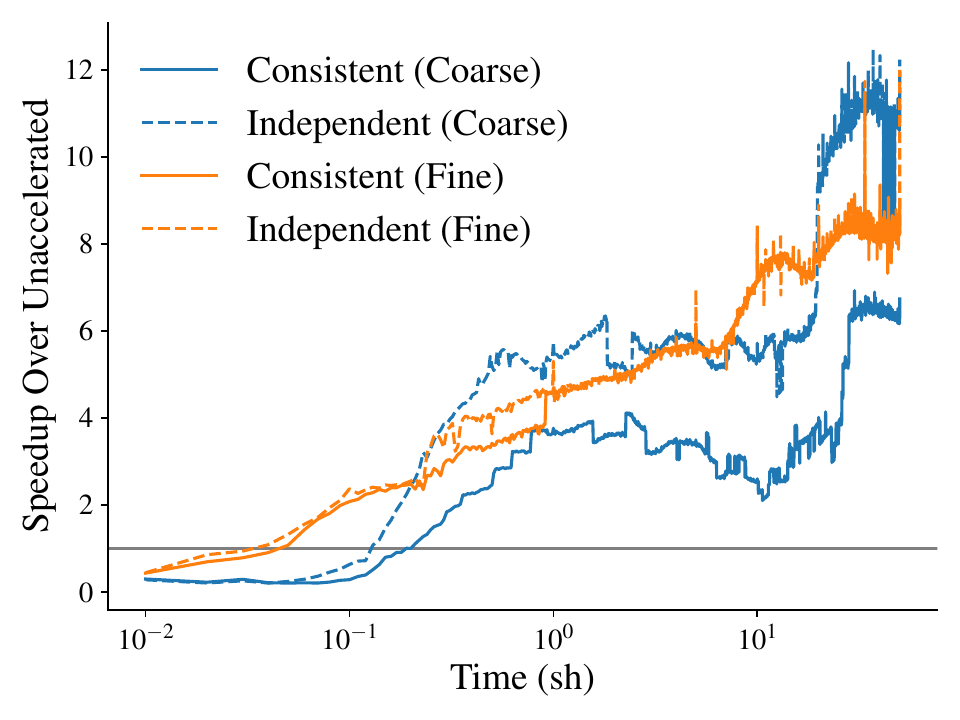}
		\caption{}
		\label{fig:gcp_performance_speedup}
	\end{subfigure}
	\begin{subfigure}{0.4\textwidth}
		\centering
		\includegraphics[width=\textwidth]{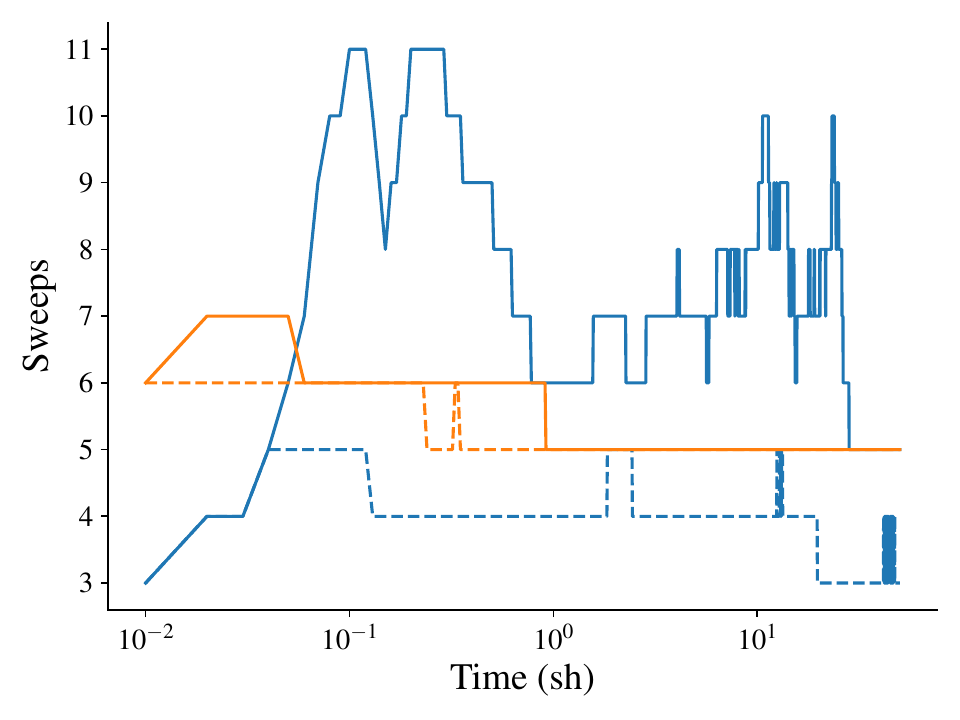}
		\caption{}
		\label{fig:gcp_performance_sweeps}
	\end{subfigure}
	\begin{subfigure}{0.4\textwidth}
		\centering
		\includegraphics[width=\textwidth]{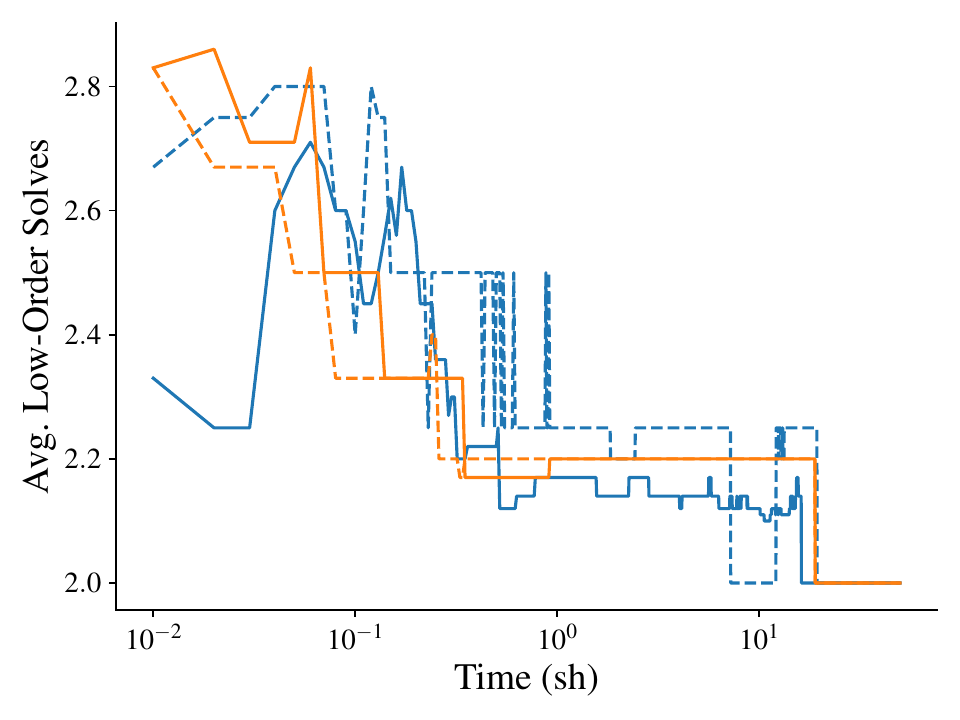}
		\caption{}
		\label{fig:gcp_performance_nonlinear}
	\end{subfigure}
	\begin{subfigure}{0.4\textwidth}
		\centering
		\includegraphics[width=\textwidth]{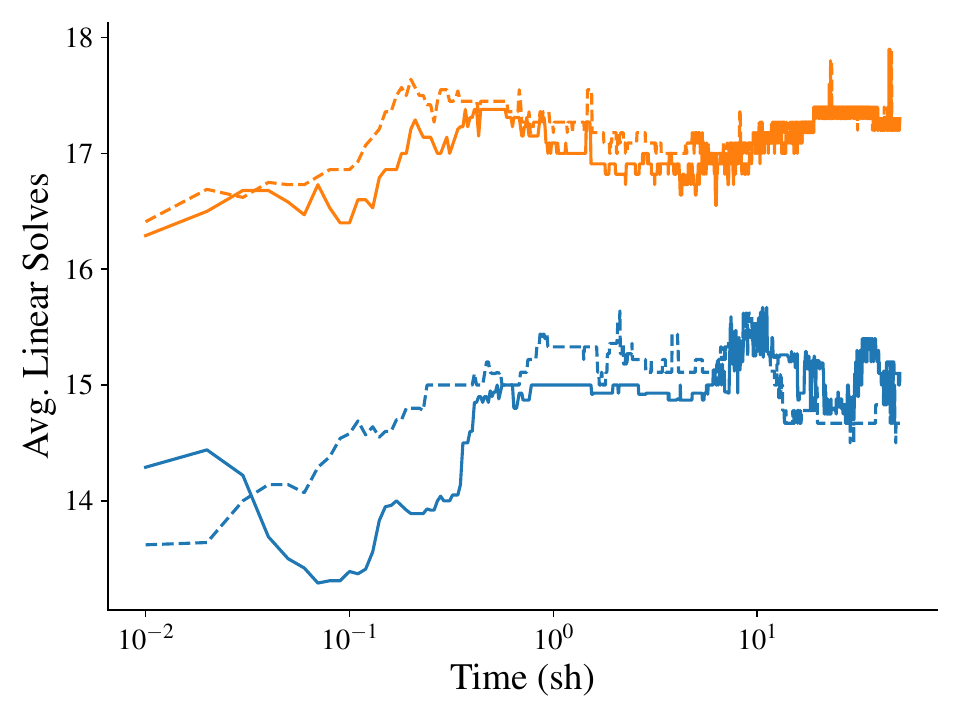}
		\caption{}
		\label{fig:gcp_performance_linear}
	\end{subfigure}
	\caption{Performance metrics as a function of simulation time for the consistent and independent SM methods on the coarse and fine crooked pipe problems.}
	\label{fig:gcp_performance} 
\end{figure}
We now discuss performance of the methods on the coarse and fine meshes. 
We consider four metrics: speedup with respect to the unaccelerated scheme, sweeps (or equivalently the number of outer iterations), average number of nonlinear low-order Newton iterations at each outer iteration, and the average number of linear, low-order solves at each outer iteration. 
These four metrics are plotted with respect to simulation time in Fig.~\ref{fig:gcp_performance}. 
A summary of iterations to convergence and cost of the sweep, nonlinear Newton iteration, linear preconditioned conjugate gradient iteration, and temperature elimination solver is provided in Table \ref{tab:summary_gcp}. 
Table \ref{tab:speedup_gcp} summarizes the speedup of the SM methods with respect to the unaccelerated method in terms of sweeps and runtime. 
At early times, absorption-emission is mild and the unaccelerated method converges as rapidly as the SM methods. 
This results in the SM methods being slower than unaccelerated in cost per time step due to the additional cost of solving the low-order system. 
However, as radiation begins to impinge on the inner blocker and turn the corner, the unaccelerated scheme's convergence degrades while the SM methods maintain the same rapid convergence. 
By the end of the simulation, where the problem has become very diffusive, the independent and consistent methods on the coarse mesh are $12\times$ and $6\times$ faster than the unaccelerated method. 
In total runtime, the coarse-mesh independent and consistent SM methods were $8.2\times$ and $4.3\times$ faster than unaccelerated. 
Recall that the independent method had very poor solution quality in both the temperature and flux on the coarse mesh. 
Thus, the nearly $2\times$ speedup of independent over consistent in total runtime on the coarse mesh comes at the cost of reduced robustness and solution quality. 

\begin{table}
	\centering
	\caption{Summary of iterative performance and cost of the sweep, low-order nonlinear solve, low-order linear solve, and nonlinear temperature elimination for the consistent SM, independent SM, and unaccelerated algorithms on the coarse and fine crooked pipe problems.}
	\label{tab:summary_gcp}
	\begin{tabular}{cccccccccc}
\toprule
 &  & \multicolumn{3}{c}{Coarse}  &  & \multicolumn{3}{c}{Fine} \\
\cmidrule{3-5}\cmidrule{7-9}
 & Metric & Consistent & Independent & Unaccelerated & & Consistent & Independent & Unaccelerated \\
\midrule
\multirow{6}{*}{\rotatebox{90}{Sweeps}} & Min & 3 & 3 & 4 & & 5 & 5 & 15 \\
 & Max & 11 & 5 & 113 & & 7 & 6 & 127 \\
 & Avg. & 6.47 & 3.47 & 100.48 & & 5.02 & 5.00 & 116.16 \\
 & Total & \num{32342} & \num{17363} & \num{502392} & & \num{25095} & \num{25025} & \num{580791} \\
 & Cost (s) & \num{572.5} & \num{311.8} & \num{7822.2} & & \num{821.8} & \num{810.9} & \num{16630.6} \\
 & Avg. Cost (ms) & \num{17.7} & \num{18.0} & \num{15.6} & & \num{32.7} & \num{32.4} & \num{28.6} \\
\addlinespace
\multirow{6}{*}{\rotatebox{90}{LO Nonlinear}} & Min & 2 & 2 & -- & & 2 & 2 & -- \\
 & Max & 3 & 4 & -- & & 4 & 4 & -- \\
 & Avg. & 2.05 & 2.07 & -- & & 2.08 & 2.08 & -- \\
 & Total & \num{66428} & \num{36240} & -- & & \num{52142} & \num{51978} & -- \\
 & Cost (s) & \num{2452.9} & \num{1287.5} & -- & & \num{2472.4} & \num{2479.8} & -- \\
 & Avg. Cost (ms) & \num{75.8} & \num{74.1} & -- & & \num{98.5} & \num{99.1} & -- \\
\addlinespace
\multirow{6}{*}{\rotatebox{90}{LO Linear}} & Min & 11 & 10 & -- & & 11 & 12 & -- \\
 & Max & 16 & 16 & -- & & 19 & 19 & -- \\
 & Avg. & 15.09 & 14.90 & -- & & 17.23 & 17.28 & -- \\
 & Total & \num{1002091} & \num{539978} & -- & & \num{898623} & \num{898090} & -- \\
 & Cost (s) & \num{685.1} & \num{370.7} & -- & & \num{916.0} & \num{941.6} & -- \\
 & Avg. Cost (ms) & \num{10.3} & \num{10.2} & -- & & \num{17.6} & \num{18.1} & -- \\
\addlinespace
\multirow{6}{*}{\rotatebox{90}{Energy Balance}} & Min & 1 & 1 & 4 & & 1 & 1 & 4 \\
 & Max & 5 & $40^*$ & 5 & & $40^*$ & $40^*$ & 10 \\
 & Avg. & 3.31 & 4.05 & 4.45 & & 8.75 & 3.92 & 4.43 \\
 & Total & \num{220071} & \num{146673} & \num{2237996} & & \num{456143} & \num{203671} & \num{2573736} \\
 & Cost (s) & \num{419.2} & \num{226.9} & \num{5118.2} & & \num{459.7} & \num{469.0} & \num{7953.5} \\
 & Avg. Cost (ms) & \num{6.3} & \num{6.3} & \num{10.2} & & \num{8.8} & \num{9.0} & \num{13.7} \\
\addlinespace
 & Wall Time (s) & \num{3146.4} & \num{1665.0} & \num{13636.5} & & \num{3373.0} & \num{3342.4} & \num{25169.8} \\
\bottomrule
\end{tabular}
	$^*$ indicates maximum iterations reached
\end{table}
\begin{table}
	\centering
	\caption{The min, max, average, and total speedups and cost in cycle time over the unaccelerated scheme on the crooked pipe problem. }
	\label{tab:speedup_gcp}
	\begin{tabular}{cccccccc}
\toprule
 &  & \multicolumn{2}{c}{Coarse}  &  & \multicolumn{2}{c}{Fine} \\
\cmidrule{3-4}\cmidrule{6-7}
 & Metric & Consistent & Independent & & Consistent & Independent \\
\midrule
\multirow{4}{*}{\rotatebox{90}{Sweeps}} & Min & 1.00 & 1.00 & & 2.50 & 2.50 \\
 & Max & 22.60 & 37.67 & & 25.40 & 25.40 \\
 & Avg. & 16.71 & 30.01 & & 23.19 & 23.22 \\
 & Total & 15.53 & 28.93 & & 23.14 & 23.21 \\
\addlinespace
\multirow{4}{*}{\rotatebox{90}{Cost}} & Min & 0.21 & 0.21 & & 0.44 & 0.45 \\
 & Max & 6.93 & 12.47 & & 11.45 & 12.00 \\
 & Avg. & 4.70 & 8.64 & & 7.50 & 7.56 \\
 & Total & 4.33 & 8.19 & & 7.46 & 7.53 \\
\bottomrule
\end{tabular}
\end{table}
The discrepancy in cost between consistent and independent is primarily due to the slower convergence of the consistent outer iteration compared to the independent method. 
Figure \ref{fig:gcp_performance_sweeps} shows the number of sweeps (and thus outer iterations) at each time step. 
On the coarse mesh, consistent requires up to 11 iterations to converge whereas independent converges in at most 5. 
This accumulates into the consistent method performing $1.8\times$ more sweeps than the independent method over the course of the simulation. 
Note that both schemes still significantly reduce the number of sweeps compared to the unaccelerated method: 
consistent and independent required $15.5\times$ and $28.9\times$ fewer sweeps than the unaccelerated method on the coarse mesh. 
On the fine mesh, consistent and independent converged in cost with $7.46\times$ and $7.53\times$ faster time-to-solution than the unaccelerated method, respectively. 
Both methods converge uniformly with 5-7 sweeps per time step, resulting in the consistent and independent methods performing $23.1\times$ and $23.2\times$ fewer sweeps across the entire simulation, respectively. 
Thus, on resolved meshes, the independent method attains good solution quality but no longer provides a significant performance advantage over consistent. 

The average number of nonlinear low-order solves per outer iteration has very small variance across the two methods and two meshes; all converged in fewer than four nonlinear iterations. 
The average number of linear iterations is also largely independent of low-order discretization. 
Both methods converge in an average of 15 iterations across all time steps on the coarse mesh. 
On the fine mesh, the average increases to 17. 
This slight decrease in iterative efficiency is likely due to the conditioning associated with a non-uniform mesh where the smallest element is $317\times$ smaller than the largest element. 

\subsubsection{Angular Quadrature Sensitivity}
\begin{figure}
\centering
\begin{subfigure}{0.24\textwidth}
	\centering
	\begin{overpic}[width=\textwidth, trim={0 0 1825 0}, clip]{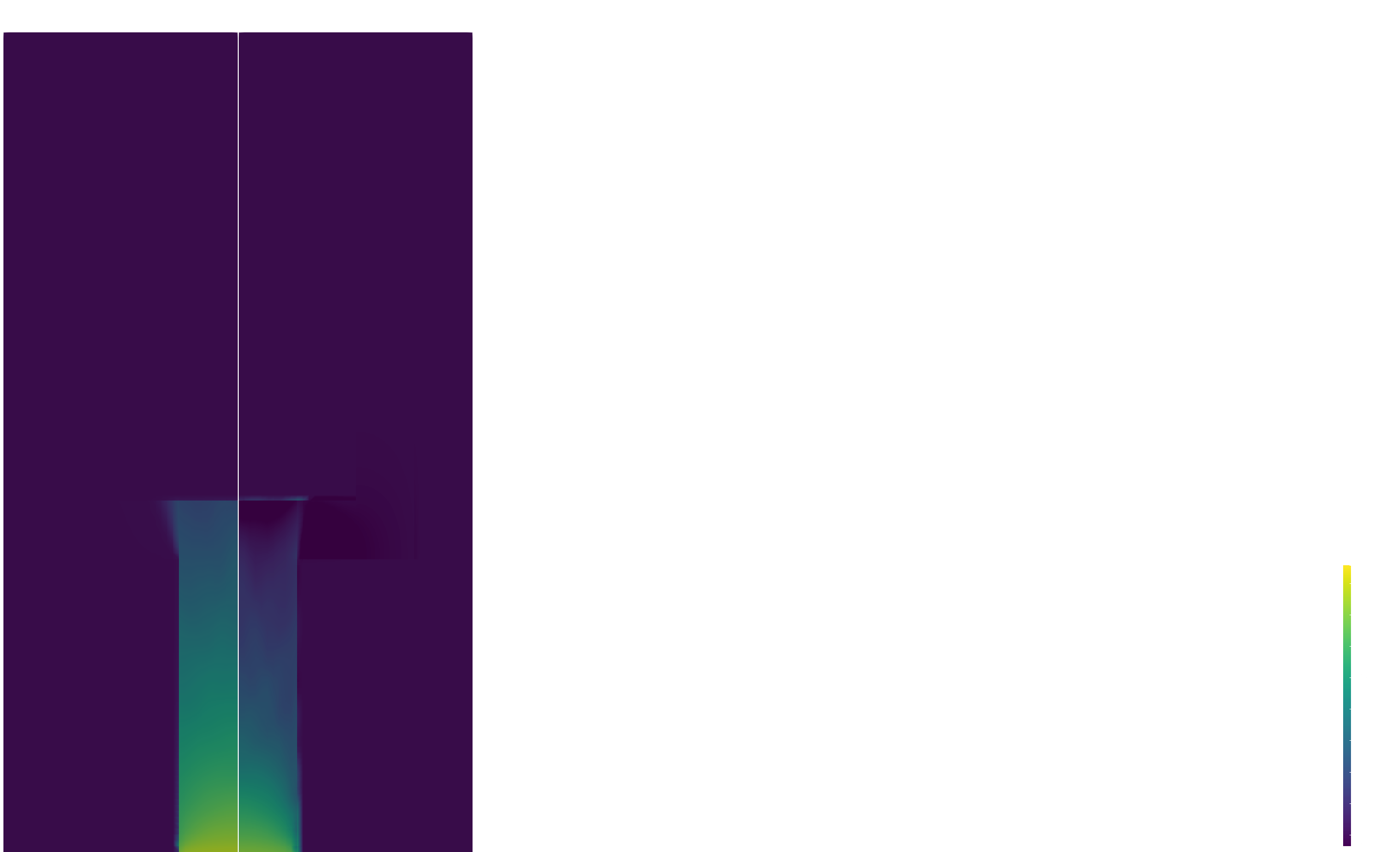}
		\put(6,92){\textcolor{white}{\scalebox{0.6}{Consistent}}}
		\put(32,92){\textcolor{white}{\scalebox{0.6}{Independent}}}
	\end{overpic}
	\caption{$t = \SI{0.1}{\sh}$}
\end{subfigure}
\begin{subfigure}{0.24\textwidth}
	\centering
	\begin{overpic}[width=\textwidth, trim={0 0 1825 0}, clip]{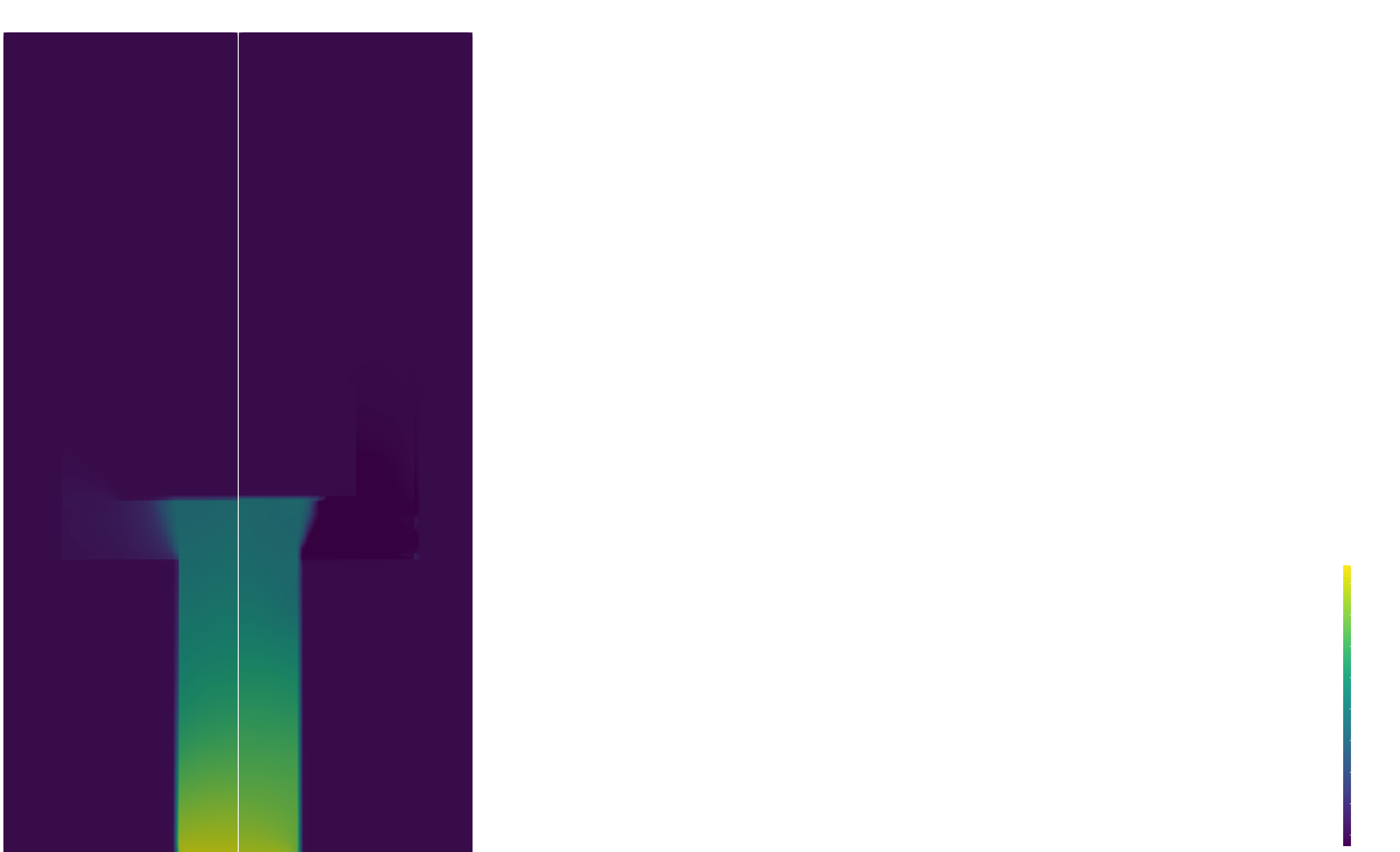}
		\put(6,92){\textcolor{white}{\scalebox{0.6}{Consistent}}}
		\put(32,92){\textcolor{white}{\scalebox{0.6}{Independent}}}
	\end{overpic}
	\caption{$t = \SI{0.5}{\sh}$}
\end{subfigure}
\begin{subfigure}{0.24\textwidth}
	\centering
	\begin{overpic}[width=\textwidth, trim={0 0 1825 0}, clip]{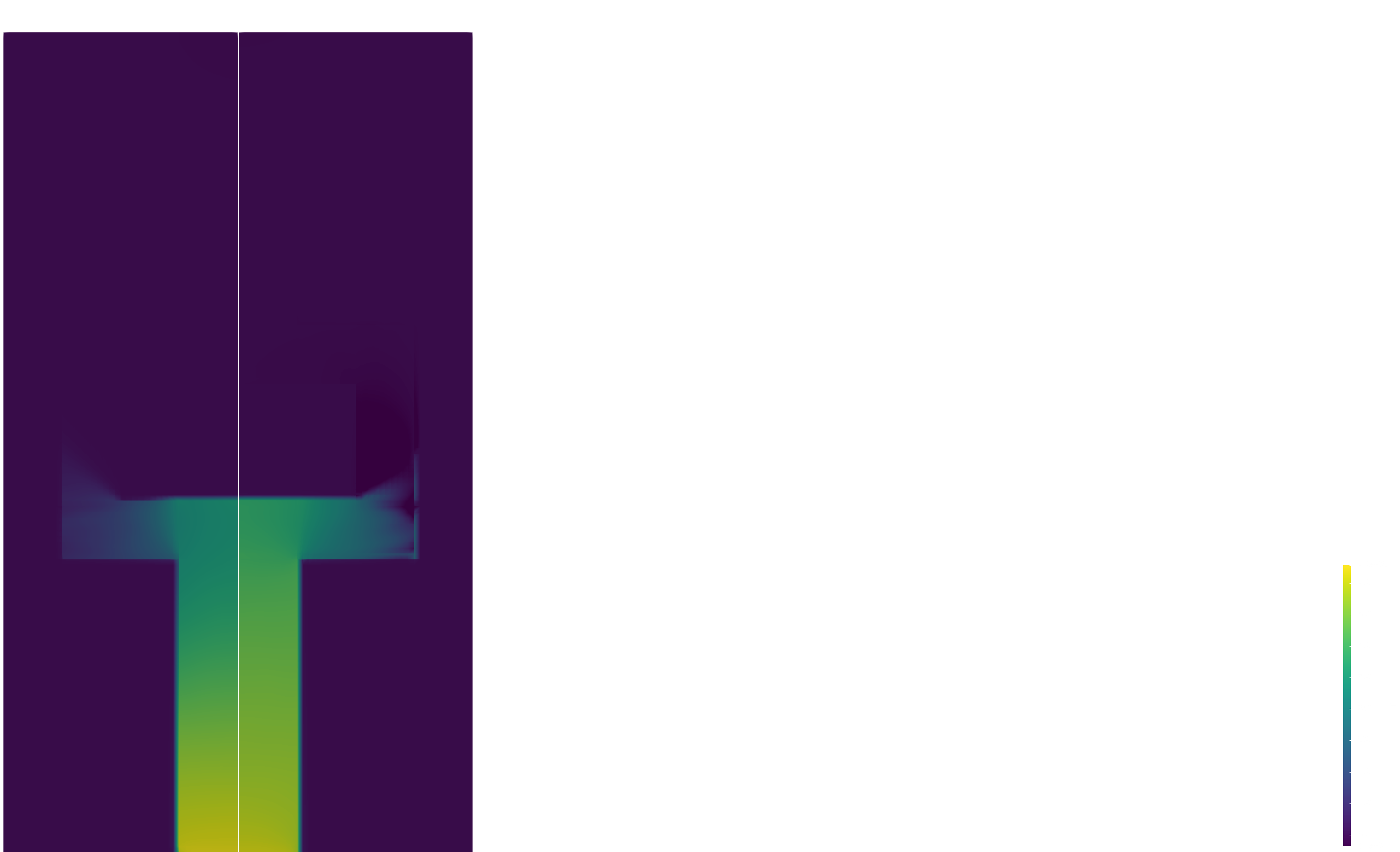}
		\put(6,92){\textcolor{white}{\scalebox{0.6}{Consistent}}}
		\put(32,92){\textcolor{white}{\scalebox{0.6}{Independent}}}
	\end{overpic}
	\caption{$t = \SI{1}{\sh}$}
\end{subfigure}
\begin{subfigure}{0.24\textwidth}
	\centering
	\begin{overpic}[width=\textwidth, trim={0 0 1825 0}, clip]{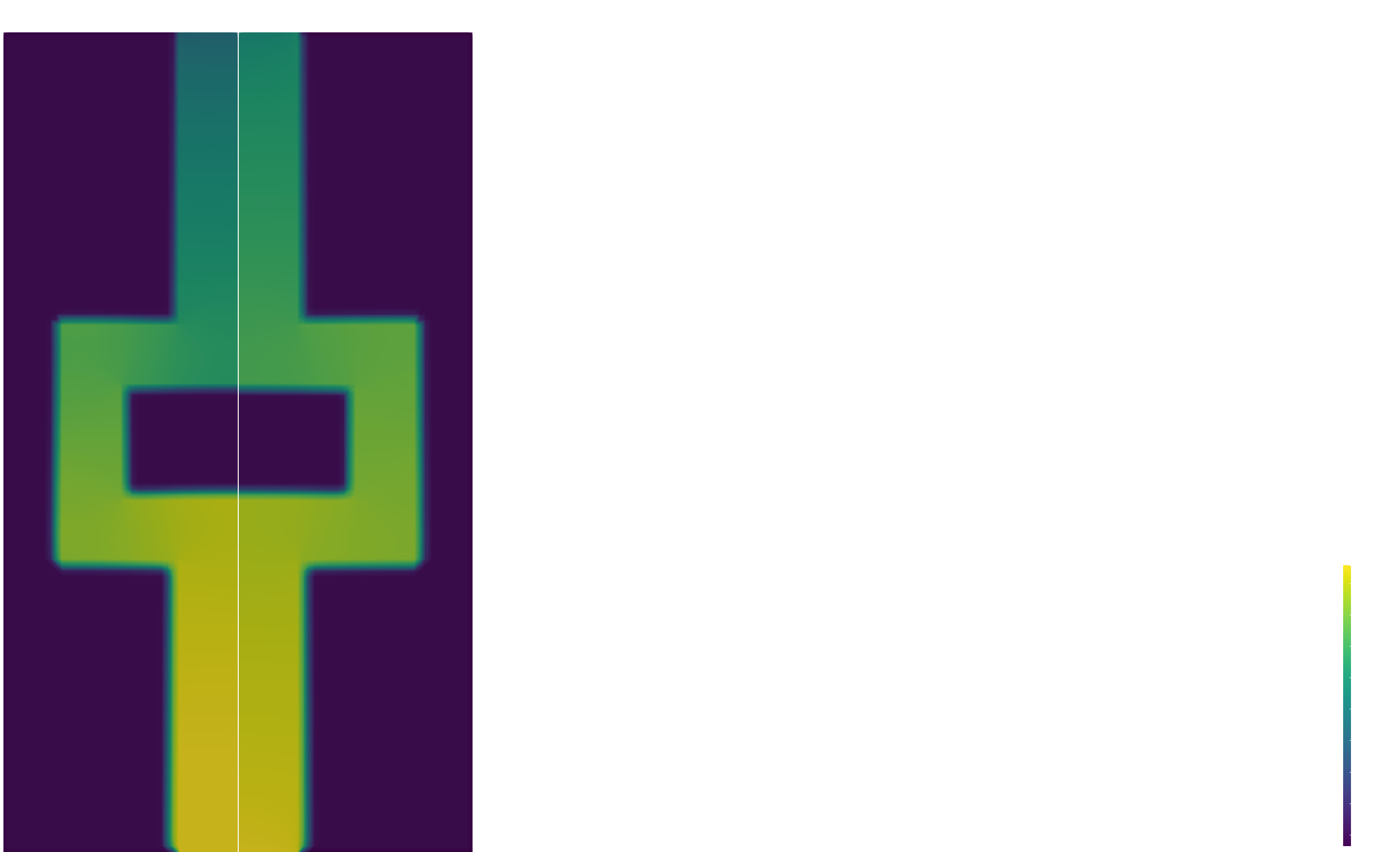}
		\put(6,92){\textcolor{white}{\scalebox{0.6}{Consistent}}}
		\put(32,92){\textcolor{white}{\scalebox{0.6}{Independent}}}
	\end{overpic}
	\caption{$t = \SI{50}{\sh}$}
\end{subfigure}
\caption{Temperature profiles over time from the consistent and independent methods on the coarse-mesh crooked pipe problem where the angular quadrature rule has been changed to $S_{12}$.}
\label{fig:gcp12_temp}
\end{figure}

\begin{figure}
\centering
\begin{subfigure}{0.24\textwidth}
	\centering
	\begin{overpic}[width=\textwidth, trim={0 0 1825 0}, clip]{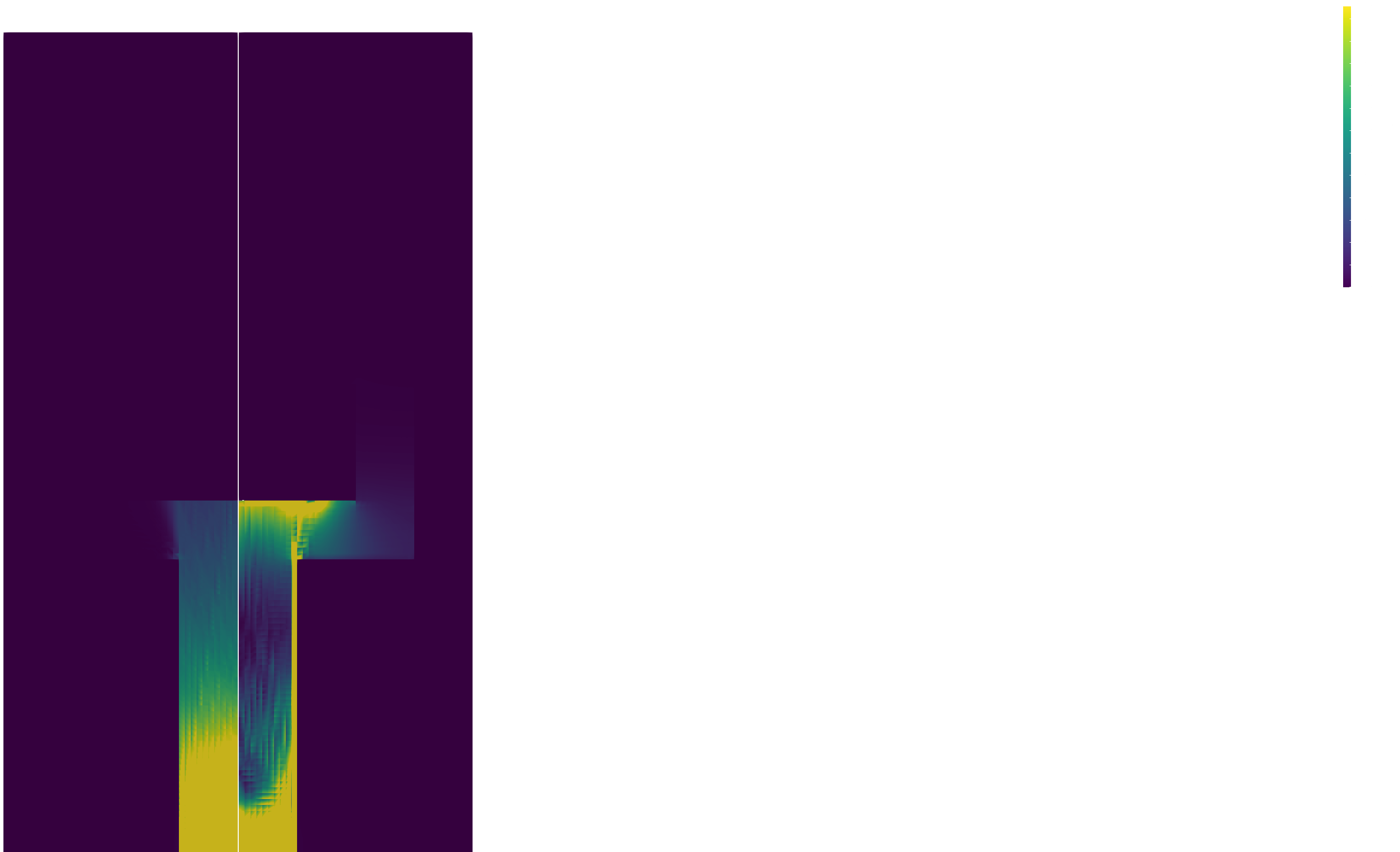}
		\put(6,92){\textcolor{white}{\scalebox{0.6}{Consistent}}}
		\put(32,92){\textcolor{white}{\scalebox{0.6}{Independent}}}
	\end{overpic}
	\caption{$t = \SI{0.1}{\sh}$}
\end{subfigure}
\begin{subfigure}{0.24\textwidth}
	\centering
	\begin{overpic}[width=\textwidth, trim={0 0 1825 0}, clip]{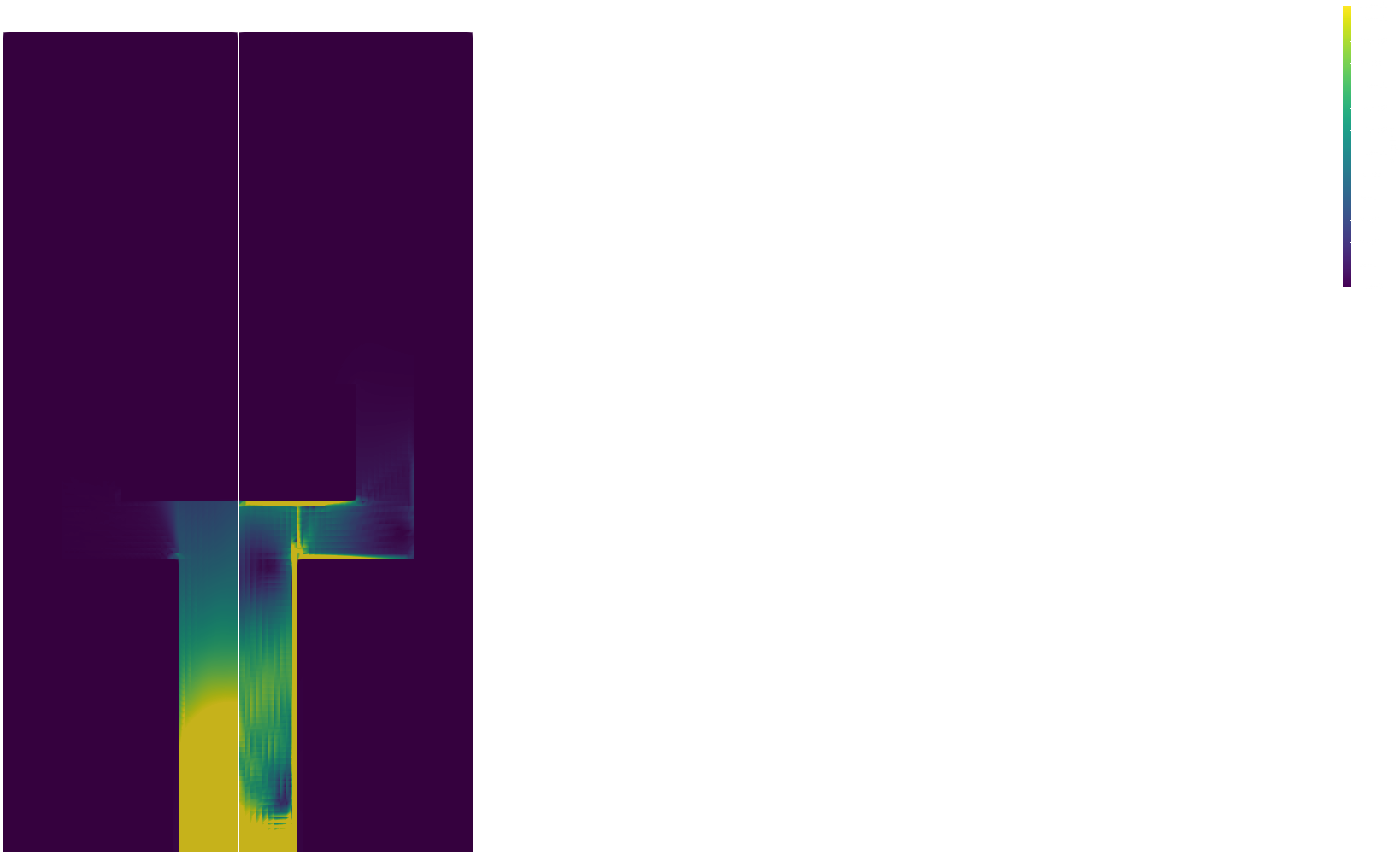}
		\put(6,92){\textcolor{white}{\scalebox{0.6}{Consistent}}}
		\put(32,92){\textcolor{white}{\scalebox{0.6}{Independent}}}
	\end{overpic}
	\caption{$t = \SI{0.5}{\sh}$}
\end{subfigure}
\begin{subfigure}{0.24\textwidth}
	\centering
	\begin{overpic}[width=\textwidth, trim={0 0 1825 0}, clip]{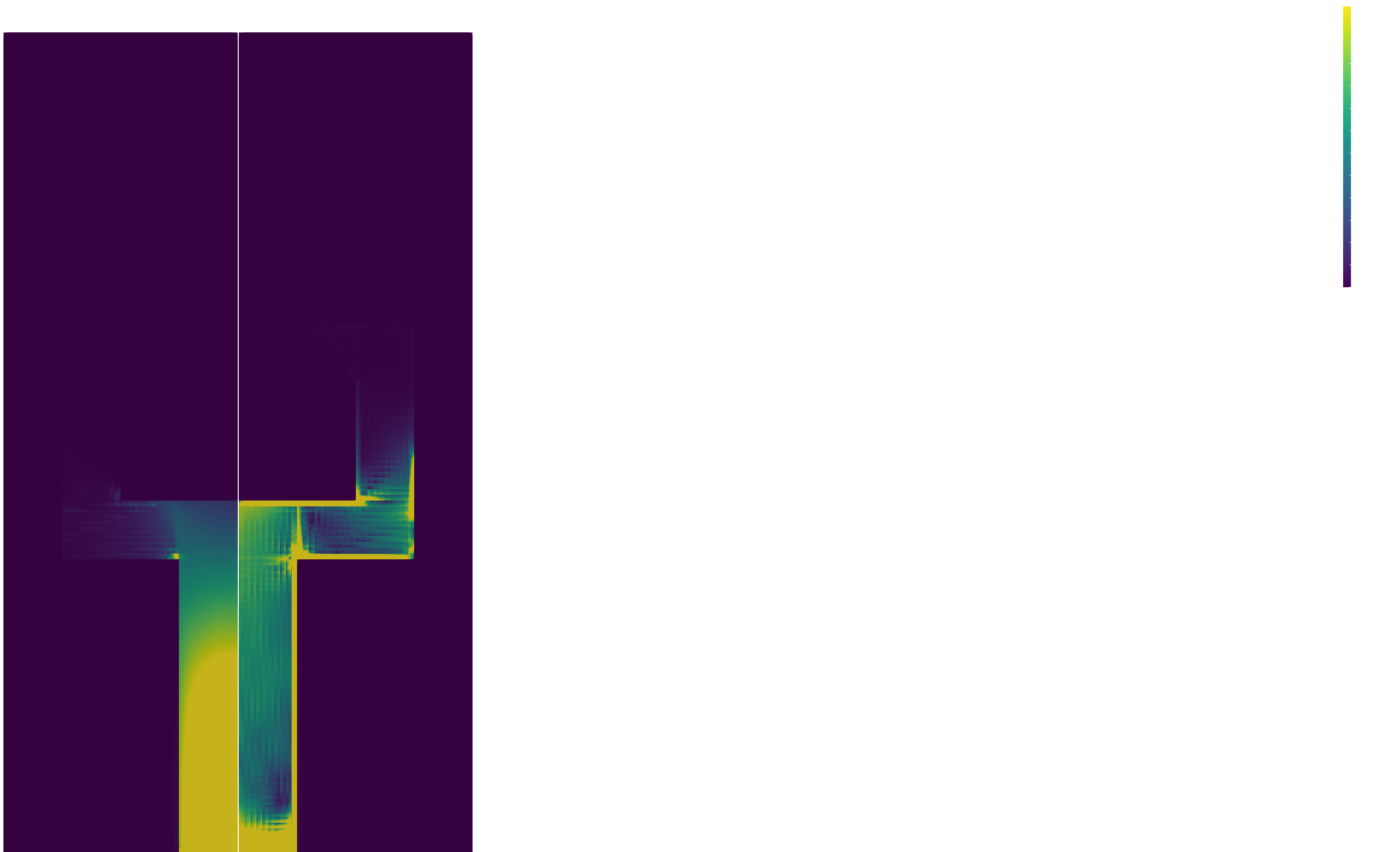}
		\put(6,92){\textcolor{white}{\scalebox{0.6}{Consistent}}}
		\put(32,92){\textcolor{white}{\scalebox{0.6}{Independent}}}
	\end{overpic}
	\caption{$t = \SI{1}{\sh}$}
\end{subfigure}
\begin{subfigure}{0.24\textwidth}
	\centering
	\begin{overpic}[width=\textwidth, trim={0 0 1825 0}, clip]{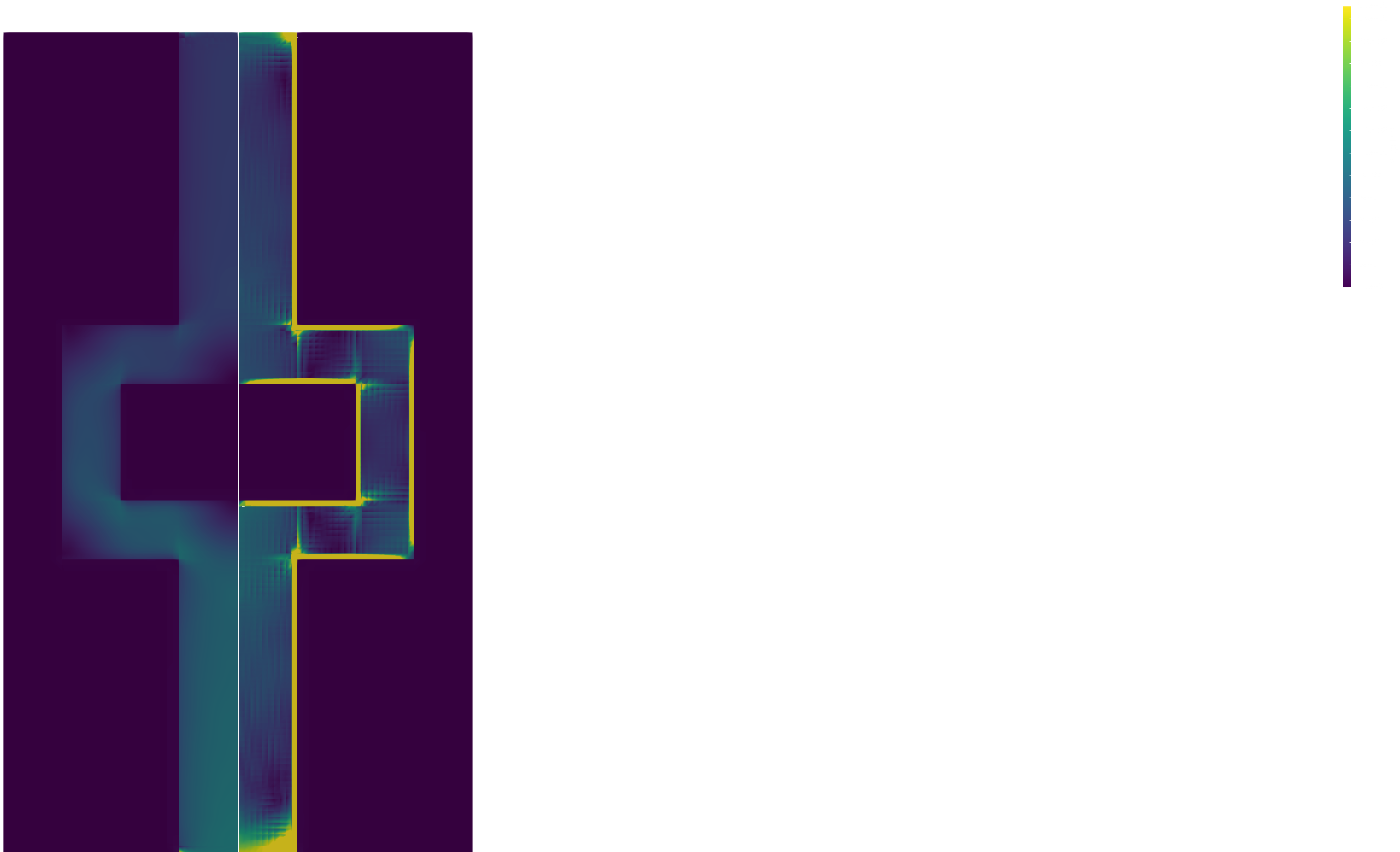}
		\put(6,92){\textcolor{white}{\scalebox{0.6}{Consistent}}}
		\put(32,92){\textcolor{white}{\scalebox{0.6}{Independent}}}
	\end{overpic}
	\caption{$t = \SI{50}{\sh}$}
\end{subfigure}
\caption{The magnitude of the flux over time from the consistent and independent methods on the coarse-mesh crooked pipe problem where the angular quadrature has been changed to $S_{12}$.}
\label{fig:gcp12_flux}
\end{figure}

We now investigate how the above comparison is affected by increasing the angular resolution. 
We consider the coarse-mesh crooked pipe problem setup where the $S_6$ square Chebyshev-Legendre quadrature rule with 36 angles has been replaced by the $S_{12}$ square Chebyshev-Legendre quadrature rule with 144 angles. 
All other parameters remain the same. 
In particular, we run the $S_{12}$ problem with the same number of MPI ranks and OpenMP threads as were used on the coarse mesh crooked pipe problem with $S_6$ angular quadrature rule. 
Figures \ref{fig:gcp12_temp} and \ref{fig:gcp12_flux} show the evolution of the temperature and radiative flux over time comparing the consistent and independent schemes. 
Compared to the coarse-mesh, $S_6$ solutions in Figs.~\ref{fig:gcp_temperature} and \ref{fig:gcp_flux}, the $S_{12}$ simulations show smoother filling of the pipe. 
However, despite the increased angular resolution, the independent solution shows similar issues with negativity during the filling of the pipe that are absent in the solutions generated by the consistent scheme. 
Likewise, the radiative flux shows the same hotspot behavior where the flux accumulates along the pipe-wall interface. 
The consistent method's temperature and flux are of similar spatial solution quality to the $S_6$ case on the same mesh. 

\begin{figure}
	\centering
	\foreach \x in {0,4,5,1,2,6,7,3}{
		\begin{subfigure}{0.24\textwidth}
			\centering
			\includegraphics[width=\textwidth]{gcp12_tracer\x.pdf}
			\caption{}
		\end{subfigure}	
	}
	\caption{Evolution of the temperature over time for the consistent and independent methods on the coarse-mesh crooked pipe problem where the angular quadrature rule varies between $S_6$ and $S_{12}$. Subfigures (a) -- (h) correspond to the spatial locations labeled in Fig.~\ref{fig:cp_diag}. }
	\label{fig:gcp12_tracer}
\end{figure}
Figure \ref{fig:gcp12_tracer} shows the solutions at the tracer locations over time. 
The $S_6$ solutions are reproduced to facilitate direct comparison. 
As with $S_6$, the independent method shows a non-physical cold spot at tracer location (b) where temperature flooring has been applied. 
However, with $S_{12}$, the duration the temperature spends at the floor temperature is reduced.
Notably, the $S_{12}$ independent solution shows a stronger dip below the inflow boundary condition of \SI{50}{\eV} at locations (d) and (e) than with $S_6$ suggesting that solution quality is not uniformly improved by angular refinement on a fixed mesh. 
The consistent solution shows less variance: where independent achieves a maximum tracer difference of \SI{200}{\eV} at tracer location (b), the consistent solution has maximum difference of only \SI{37}{\eV} at the same spatial location. 
The long-time behaviors for both methods are nearly identical with locations (f)--(h) all showing similar equilibrium solutions as the $S_6$ variants. 
This indicates that, in addition to the spatial sensitivity observed in \S\ref{sec:gcp_solution_quality}, the independent scheme's solution quality is also sensitive to angular resolution. 
The consistent method attains the high-order discretization's solution quality regardless of spatial mesh or angular quadrature rule and thus produces more robust solution quality than the independent method. 

\begin{figure}
	\centering
	\begin{subfigure}{0.4\textwidth}
		\centering 
		\includegraphics[width=\textwidth]{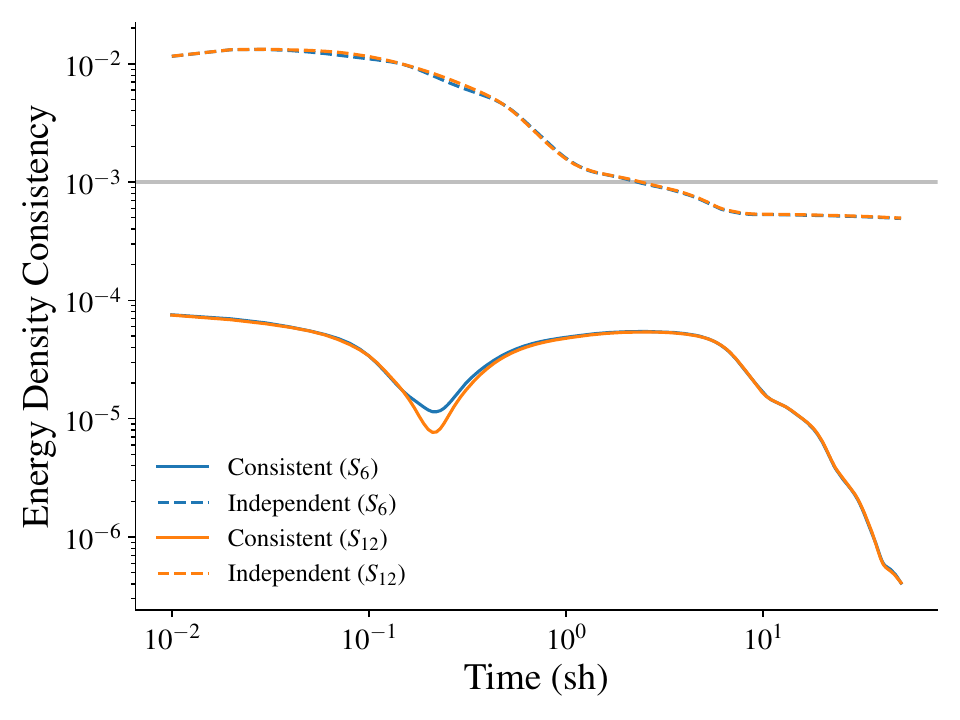}
		\caption{}
	\end{subfigure}
	\begin{subfigure}{0.4\textwidth}
		\centering 
		\includegraphics[width=\textwidth]{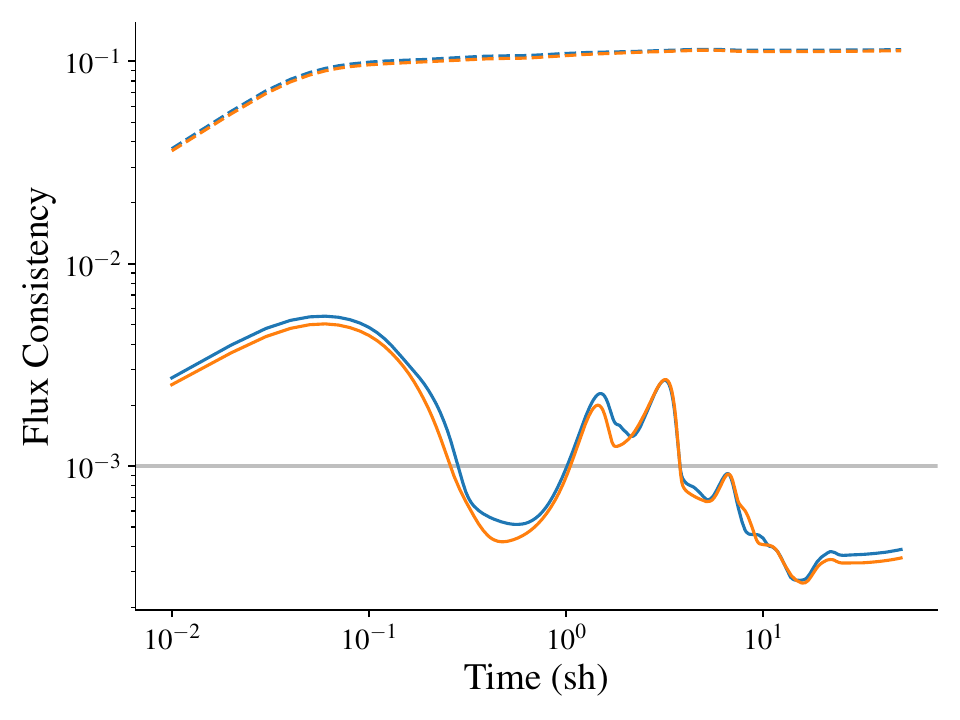}
		\caption{}
	\end{subfigure}
	\caption{The difference between the high and low-order energy density and flux as a function of simulation time for the consistent and independent SM methods on the coarse-mesh crooked pipe problem using the $S_6$ and $S_{12}$ angular quadrature rules.}
	\label{fig:gcp12_consistency}
\end{figure}
High-order low-order consistency is plotted over simulation time in Fig.~\ref{fig:gcp12_consistency}. 
The coarse-mesh, $S_6$ data is reproduced to facilitate direct comparison. 
Both consistent and independent have consistency largely independent of the angular resolution in both energy density and radiative flux. 
This indicates that the spatial inconsistency still dominates and is largely unaffected by the higher-resolution angular quadrature rule. 

\begin{table}
	\centering
	\caption{Summary of iterative performance and cost of the sweep, low-order nonlinear solve, low-order linear solve, and nonlinear temperature elimination for the consistent and independent SM methods on the coarse-mesh crooked pipe problem where the angular quadrature rule has been changed to $S_{12}$. The coarse-mesh, $S_6$ values from Table \ref{tab:summary_gcp} are reproduced for easier comparison.}
	\label{tab:summary_gcp12}
	\begin{tabular}{cccccccc}
\toprule
 &  & \multicolumn{2}{c}{$S_6$}  &  & \multicolumn{2}{c}{$S_{12}$} \\
\cmidrule{3-4}\cmidrule{6-7}
 & Metric & Consistent & Independent & & Consistent & Independent \\
\midrule
\multirow{6}{*}{\rotatebox{90}{Sweeps}} & Min & 3 & 3 & & 3 & 3 \\
 & Max & 11 & 5 & & 11 & 4 \\
 & Avg. & 6.47 & 3.47 & & 6.22 & 3.73 \\
 & Total & \num{32342} & \num{17363} & & \num{31125} & \num{18663} \\
 & Cost (s) & \num{572.5} & \num{311.8} & & \num{1978.0} & \num{1184.1} \\
 & Avg. Cost (ms) & \num{17.7} & \num{18.0} & & \num{63.5} & \num{63.4} \\
\addlinespace
\multirow{6}{*}{\rotatebox{90}{LO Nonlinear}} & Min & 2 & 2 & & 2 & 2 \\
 & Max & 3 & 4 & & 3 & 4 \\
 & Avg. & 2.05 & 2.07 & & 2.05 & 2.07 \\
 & Total & \num{66428} & \num{36240} & & \num{63990} & \num{38813} \\
 & Cost (s) & \num{2452.9} & \num{1287.5} & & \num{2709.2} & \num{1412.7} \\
 & Avg. Cost (ms) & \num{75.8} & \num{74.1} & & \num{87.0} & \num{75.7} \\
\addlinespace
\multirow{6}{*}{\rotatebox{90}{LO Linear}} & Min & 11 & 10 & & 11 & 10 \\
 & Max & 16 & 16 & & 16 & 16 \\
 & Avg. & 15.09 & 14.90 & & 15.00 & 14.96 \\
 & Total & \num{1002091} & \num{539978} & & \num{960163} & \num{580607} \\
 & Cost (s) & \num{685.1} & \num{370.7} & & \num{659.9} & \num{400.8} \\
 & Avg. Cost (ms) & \num{10.3} & \num{10.2} & & \num{10.3} & \num{10.3} \\
\addlinespace
\multirow{6}{*}{\rotatebox{90}{Energy Balance}} & Min & 1 & 1 & & 1 & 1 \\
 & Max & 5 & $40^*$ & & 5 & $40^*$ \\
 & Avg. & 3.31 & 4.05 & & 3.32 & 3.85 \\
 & Total & \num{220071} & \num{146673} & & \num{212414} & \num{149558} \\
 & Cost (s) & \num{419.2} & \num{226.9} & & \num{387.7} & \num{245.2} \\
 & Avg. Cost (ms) & \num{6.3} & \num{6.3} & & \num{6.1} & \num{6.3} \\
\addlinespace
 & Wall Time (s) & \num{3146.4} & \num{1665.0} & & \num{4860.7} & \num{2711.9} \\
\bottomrule
\end{tabular}
	$^*$ indicates maximum iterations reached
\end{table}
Table \ref{tab:summary_gcp12} compares cost and runtime of the consistent and independent schemes on the coarse-mesh crooked pipe problem for the $S_6$ and $S_{12}$ angular quadrature rules. 
As expected, the cost per sweep increases. 
The number of sweeps, nonlinear low-order iterations, linear low-order iterations, and temperature elimination iterations are insensitive to the angular quadrature rule; all iterative components of the algorithms remain rapidly converging. 
Note that the cost of computing the SM closures scales with the size of the angular quadrature rule. 
Thus, we see the average cost per inner iteration increase by 15\% and 2\% for the consistent and independent methods, respectively. 
The consistent SM source terms are more expensive to compute as they depend on half-range integration where the independent closures use full-range integration. 
Here, the low-order cost is now dwarfed by the high-order sweep cost as the low-order problem cost is less dependent on the number of angles than the sweep. 
Independent is much less expensive than the consistent method with speedups over consistent of $1.9\times$ and $1.8\times$ for $S_6$ and $S_{12}$, respectively. 
However, as with $S_6$, independent struggles with solution quality at this spatial resolution. 

\subsection{Multigroup ICF-inspired Lattice Problem}
Finally, we consider the two-dimensional, frequency-dependent lattice problem defined in \cite{osti_2280904,brunner_github}. 
This problem uses analytic opacities that are designed to mimic the ``edge'' and ``line'' features of real materials that arise in ICF simulations. 
The geometry is simplified to a lattice arrangement to make detailed mesh generation easier while still allowing a challenging array of materials and associated, frequency-dependent optical depths. 
The problem contains an initially hot material to hasten the onset of numerically stiff absorption-emission, allowing a shorter final time than the crooked pipe problem. 
As with the crooked pipe problem, resolving this problem in space-time is difficult so we focus on evaluating the algorithms in terms of their solution quality and performance. 

\subsubsection{Definition}
\begin{figure}
	\centering
	\includegraphics[width=0.65\textwidth]{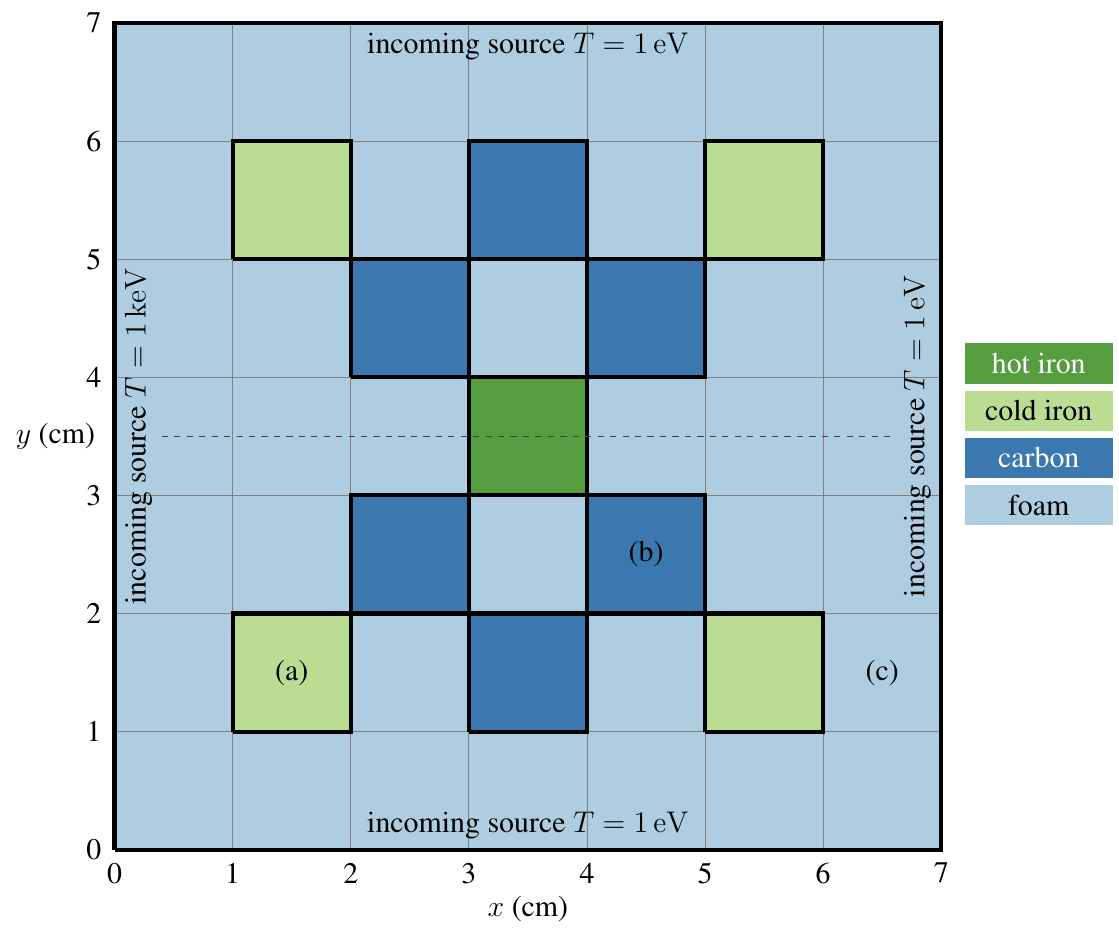}
	\caption{Problem description for a lattice problem with four materials. A reflection plane along $y=\SI{3.5}{\cm}$ is used. The solutions at the spatial locations (a), (b), and (c) are recorded at each time step. }
	\label{fig:brunner_diag}
\end{figure}
\begin{table}
	\centering
	\caption{Properties of the materials in the multigroup lattice problem.}
	\label{tab:brunner_mats}
	\begin{tabular}{c ccc}
		\toprule
		& Density (\si{\g\per\cm\cubed}) & Heat Capacity (\si{\erg\per\cm\cubed\per\eV}) & Initial Temperature (\si{\eV}) \\ \midrule
		Hot Iron & 8 & \num{5.4273e11} & 500 \\ 
		Cold Iron & 6 & \num{5.4273e11} & 1 \\ 
		Carbon & 2 & \num{2.41213e11} & 1 \\ 
		Foam & 0.2 & \num{2.41213e11} & 1 \\
		\bottomrule
	\end{tabular}
\end{table}
This problem consists of four materials in the $\SI{7}{\cm}\times \SI{7}{\cm}$ lattice geometry shown in Fig.~\ref{fig:brunner_diag}. 
Table \ref{tab:brunner_mats} defines the properties of each material. 
Table \ref{tab:brunner_mats} also lists the initial temperatures of the blocks: all materials begin at \SI{1}{\eV} except for the centermost block, labeled ``hot iron'', which has a temperature of \SI{500}{\eV}. 
The presence of an initially hot material reduces the lead time for stiff absorption-emission physics to develop in the problem, resulting in stiff physics throughout the simulation's runtime. 
We consider a two-dimensional Cartesian variant of the problem where the left face has a \SI{1}{\keV} boundary source and the bottom, top, and right faces have vacuum conditions mocked with a fictitious, small-valued source of \SI{1}{\eV}.
With these boundary conditions, the exact solution satisfies $\SI{1}{\eV} \leq T \leq \SI{1}{\keV}$. 
Since the lattice geometry and boundary conditions are symmetric about $y=\SI{3.5}{\cm}$, we apply a reflection boundary condition along $y=\SI{3.5}{\cm}$ to halve the computational domain. 
The simulations are run to a final time of \SI{50}{\ns}. 
As with the crooked pipe, the zero and scale negative flux fixup is used during the sweep to ensure positivity of the high-order solution. 

In \cite{osti_2280904}, analytic opacities have been designed to mimic the ``edge'' and ``line'' features of real opacities. 
Edges are modeled with a Heaviside function that raises the opacity by a fixed amount past a prescribed energy, $\nu_\text{edge}$. 
Lines are modeled with a series of thin Gaussians preceding an edge with width $\delta_w$ and spacing $\delta_s$. 
In addition, a frequency cutoff is included so that the opacity remains bounded as frequency approaches zero. 
With all these features incorporated, the opacities are of the form 
	\begin{subequations}
	\begin{equation} \label{eq:brunner_opac}
		\sigma(\nu,T,\rho) = C_0 \frac{\rho^2}{\sqrt{T} \hat{\nu}^3}\paren{1 - e^{-\hat{\nu}/T}}\bracket{1 + C_1 H(\hat{\nu} - \nu_\text{edge}) + \sum_{\ell=0}^{N_\ell-1} \frac{C_2}{N_\ell - \ell}\exp\!\paren{-\frac{1}{2}\bracket{\frac{\hat{\nu} - (\nu_\text{edge} - (\ell+1)\delta_s)}{\delta_w}}^2}} \,,
	\end{equation}
	\begin{equation}
		\hat{\nu} = \max(\nu_\text{min}, \nu) \,,
	\end{equation}
	\end{subequations}
where $C_0$, $C_1$, and $C_2$ control the magnitudes of the overall opacity, the edge, and the lines, respectively, $\nu_\text{min}$ enforces the frequency cutoff, and $N_\ell$ is the number of lines. 
Table \ref{tab:brunner_opac} defines the opacity parameters associated with Eq.~\ref{eq:brunner_opac} for each material. 
Detailed plots of the frequency and temperature dependence of the four opacities are provided in \cite[Figs.~3-6]{osti_2280904}. 
\begin{table}
	\centering
	\caption{Opacity parameters for each material.}
	\label{tab:brunner_opac}
	\begin{tabular}{c cccccccc}
		\toprule
		& $\nu_\text{min}$ (\si{\eV}) & $\nu_\text{edge}$ (\si{\eV}) & $C_0$ (\si{\cm^5\eV^{3.5}\per\g\squared}) & $C_1$ & $C_2$ & $N_\ell$ & $\delta_w$ (\si{\eV}) & $\delta_s$ (\si{\eV}) \\
		\midrule 
		Hot/Cold Iron & 50 & 7000 & \num{6.356178e11} & 1200 & 1200 & 5 & 10 & 200  \\
		Carbon & 40 & 1500 & \num{2.434954e+10} & 1200 & 30 & 1 & 10 & 1200 \\
		Foam & 40 & 300 & \num{6.324555e+10} & 400 & 0 & 0 & 0 & 0 \\ 
		\bottomrule
	\end{tabular}
\end{table}
17 groups are used in frequency according to the bounds specified in Table \ref{tab:brunner_bounds}. 
\begin{table}
	\centering 
	\caption{Energy bounds in \si{\eV}.}
	\label{tab:brunner_bounds}
	\begin{tabular}{r}
		\toprule
			\num{0} \\
			\num{0.1} \\
			\num{3} \\
			\num{10.95445} \\
			\num{40} \\
			\num{50} \\
			\num{78.25423} \\
			\num{122.4745} \\
			\num{191.6829} \\
			\num{300} \\
			\num{670.8204} \\
			\num{1500} \\
			\num{3240.37} \\
			\num{7000} \\
			\num{11146.2} \\
			\num{17748.24} \\
			\num{28260.76} \\
			\num{45000} \\
		\bottomrule
	\end{tabular}
\end{table}
Multigroup opacities are computed using the code contained in \cite{brunner_github} which has high-precision routines to integrate Eq.~\ref{eq:brunner_opac} weighted by the Planck and Rosseland spectra over each bin. 
In the results that follow, we use the Rosseland-weighted opacities generated by \cite{brunner_github} as the \emph{multigroup} absorption opacity, $\sigma_g$. 
The group structure does not resolve the line features so neither Planck nor Rosseland weighting will produce the group-converged solution. 
We choose to use Rosseland weighting to match the reference solutions provided in \cite{osti_2280904}. 
Note that the definitions of the gray opacities $\sigma_E$, $\sigma_F$, and $\sigma_P$ remain the same regardless of the weight function used within each group to generate $\sigma_g$. 

\begin{figure}
	\centering
	\begin{subfigure}{0.5\textwidth}
		\centering
		\includegraphics[width=\textwidth]{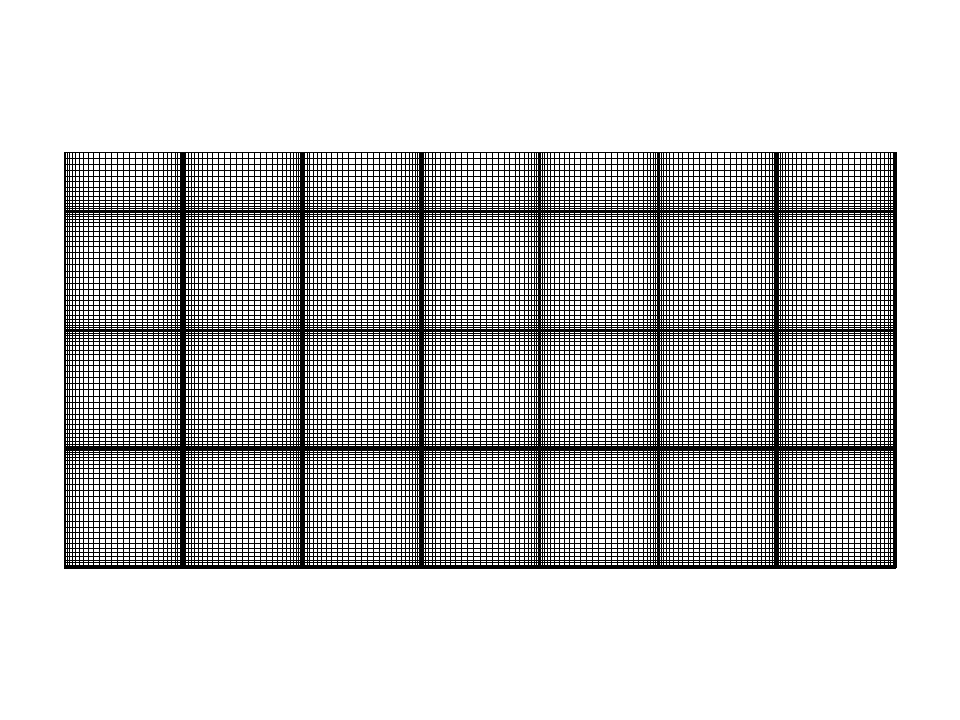}
		\caption{}
	\end{subfigure}
	\begin{subfigure}{0.4\textwidth}
		\centering
		\includegraphics[width=\textwidth]{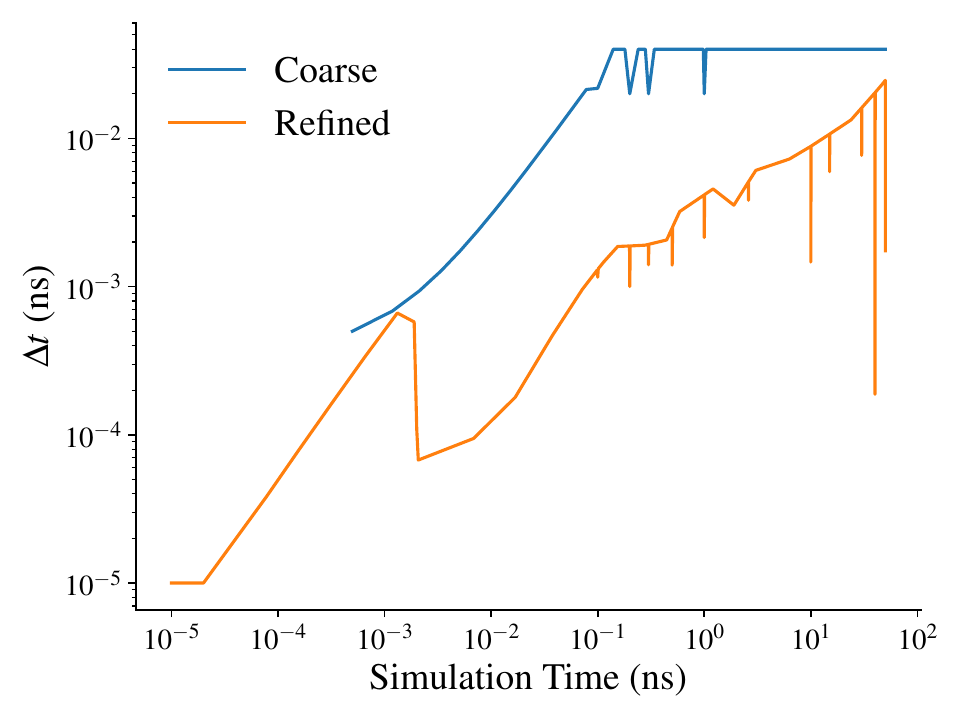}
		\caption{}
	\end{subfigure}
	\caption{The computational mesh used for the fine lattice problem and time step sizes used in the coarse and fine lattice problems.}
	\label{fig:brunner_mesh}
\end{figure}

An $S_8$ square Chebyshev-Legendre angular quadrature rule with 64 angles is used. 
Note that it is important to use a quadrature rule with an angle near parallel to the $x$-axis to properly transport energy along the top and bottom faces of the domain. 
As with the crooked pipe, we compare solution quality and performance on two meshes. 
The first is a uniform mesh with $10\times 10$ elements per $\SI{1}{\cm\squared}$ for a total of 2450 elements. 
The second is a refined mesh where each $\SI{1}{\cm}\times\SI{1}{\cm}$ block has a $30\times30$ submesh of elements distributed according to the Chebyshev points of the second kind so that elements cluster at the boundaries of each block. 
The fine mesh has \num{22050} elements and has characteristic mesh lengths ranging between \SI{2.739e-3}{\cm} and \SI{5.226e-2}{\cm}. 
On the coarse mesh, the time step begins at \SI{5e-4}{\ns} and uniformly ramps up to \SI{4e-2}{\ns} over 15 time steps. 
The fine problem uses a piecewise-linear approximation to the time step sizes taken by the deterministic method in \cite[Fig.~23]{osti_2280904}. 
The fine mesh and the coarse and fine time step sizes are shown in Fig.~\ref{fig:brunner_mesh}. 
Note that dips in the time step size occur due to the simulation adjusting the time step size before outputting snapshots in time to disk to ensure that comparisons are made at exact points in time. 
The coarse and fine problems are solved in parallel on one and four nodes with 44 processors per node. 

\subsubsection{Solution Quality}
\begin{figure}
	\centering
	\begin{subfigure}{0.45\textwidth}
		\centering
		\begin{overpic}[width=\textwidth, trim={0 0 25 250}, clip]{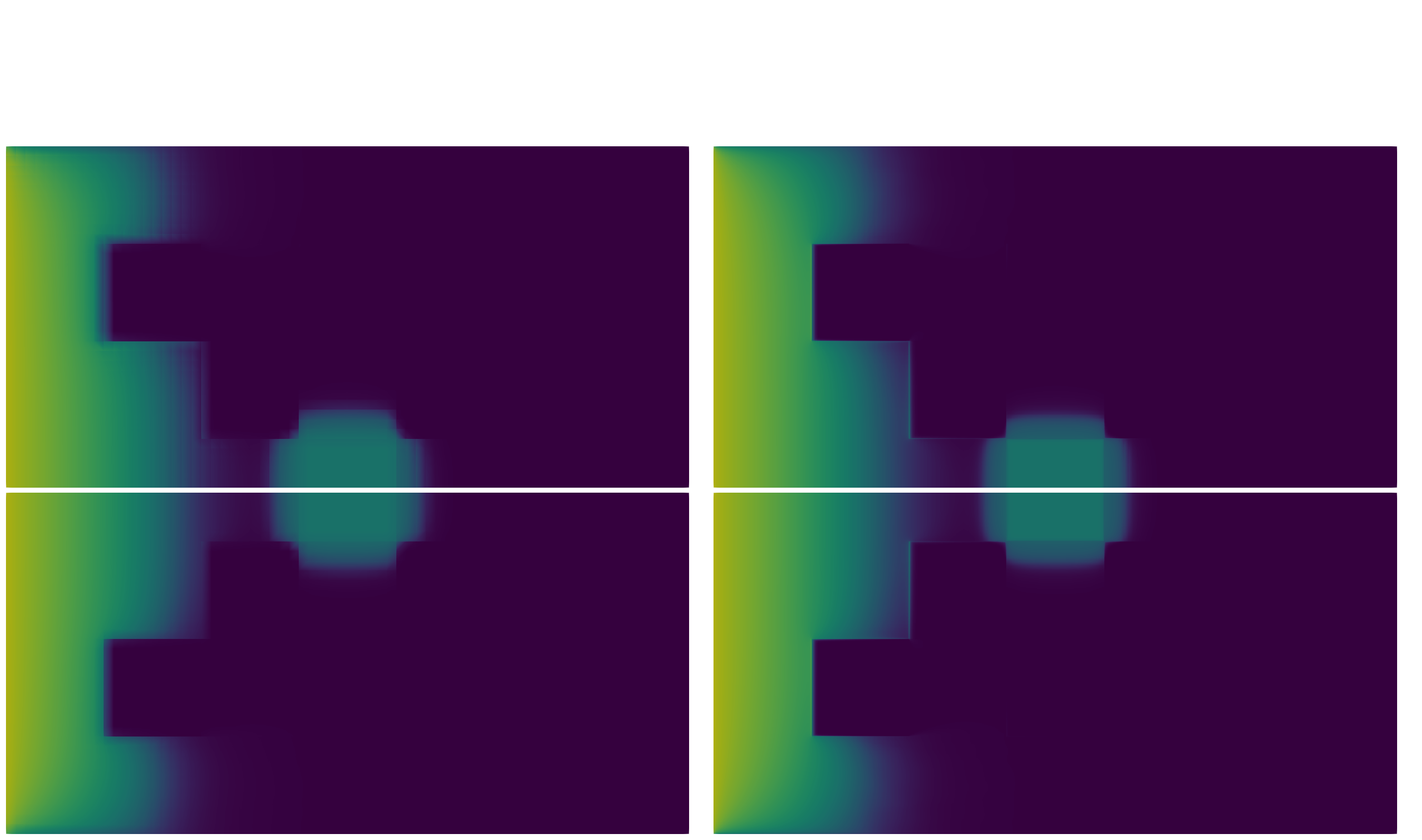}
			\put(1,1){\textcolor{white}{\scalebox{0.6}{Consistent}}}
			\put(1,26){\textcolor{white}{\scalebox{0.6}{Independent}}}
			\put(25,51){\makebox(0,0){\scalebox{0.6}{Coarse}}}

			\put(51,1){\textcolor{white}{\scalebox{0.6}{Consistent}}}
			\put(51,26){\textcolor{white}{\scalebox{0.6}{Independent}}}
			\put(75,51){\makebox(0,0){\scalebox{0.6}{Fine}}}
		\end{overpic}
		\caption{$t = \SI{0.2}{\ns}$}
	\end{subfigure}
	\quad
	\begin{subfigure}{0.45\textwidth}
		\centering
		\begin{overpic}[width=\textwidth, trim={0 0 25 250}, clip]{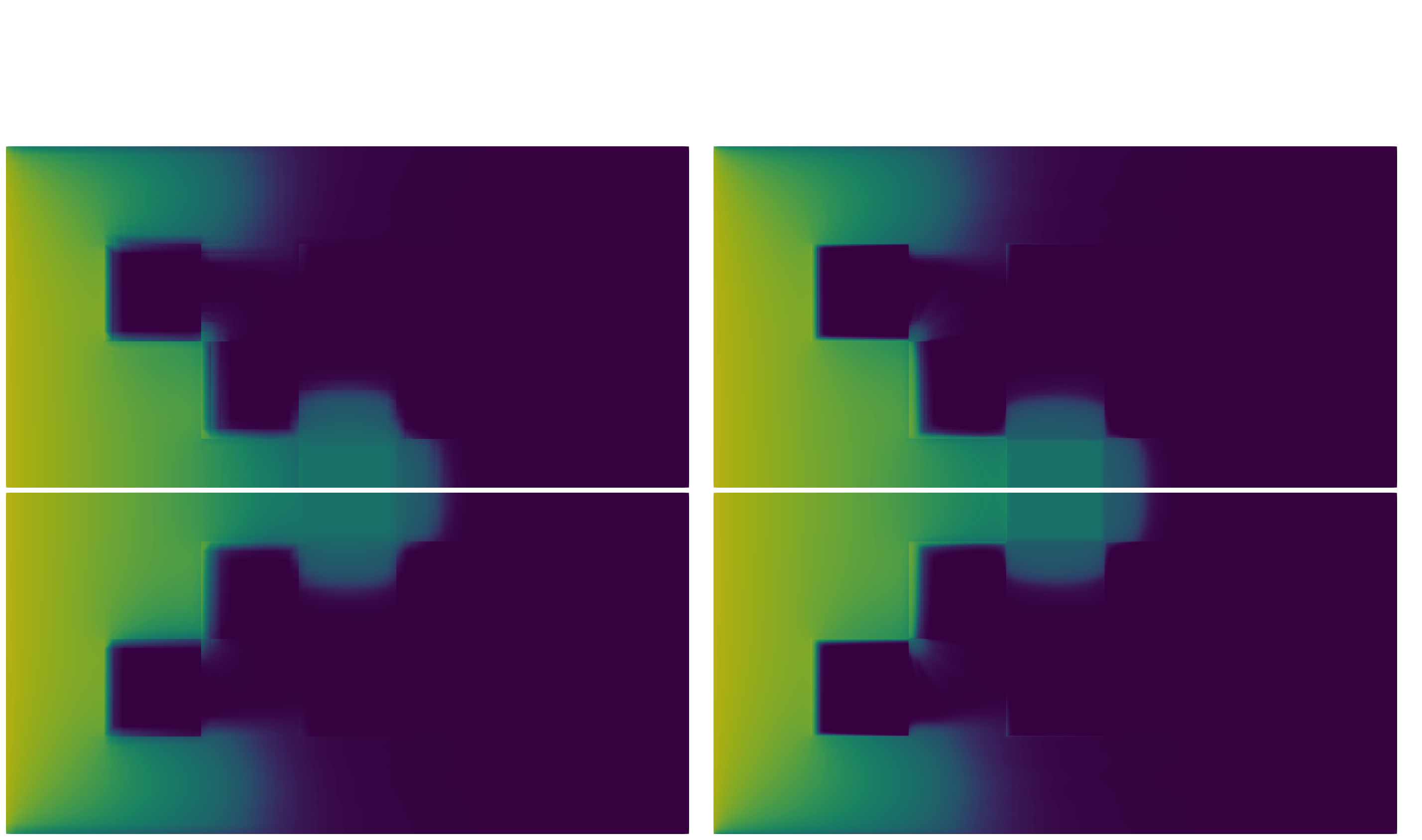}
			\put(1,1){\textcolor{white}{\scalebox{0.6}{Consistent}}}
			\put(1,26){\textcolor{white}{\scalebox{0.6}{Independent}}}
			\put(25,51){\makebox(0,0){\scalebox{0.6}{Coarse}}}

			\put(51,1){\textcolor{white}{\scalebox{0.6}{Consistent}}}
			\put(51,26){\textcolor{white}{\scalebox{0.6}{Independent}}}
			\put(75,51){\makebox(0,0){\scalebox{0.6}{Fine}}}
		\end{overpic}
		\caption{$t = \SI{0.5}{\ns}$}
	\end{subfigure}
	\begin{subfigure}{0.45\textwidth}
		\centering
		\begin{overpic}[width=\textwidth, trim={0 0 25 250}, clip]{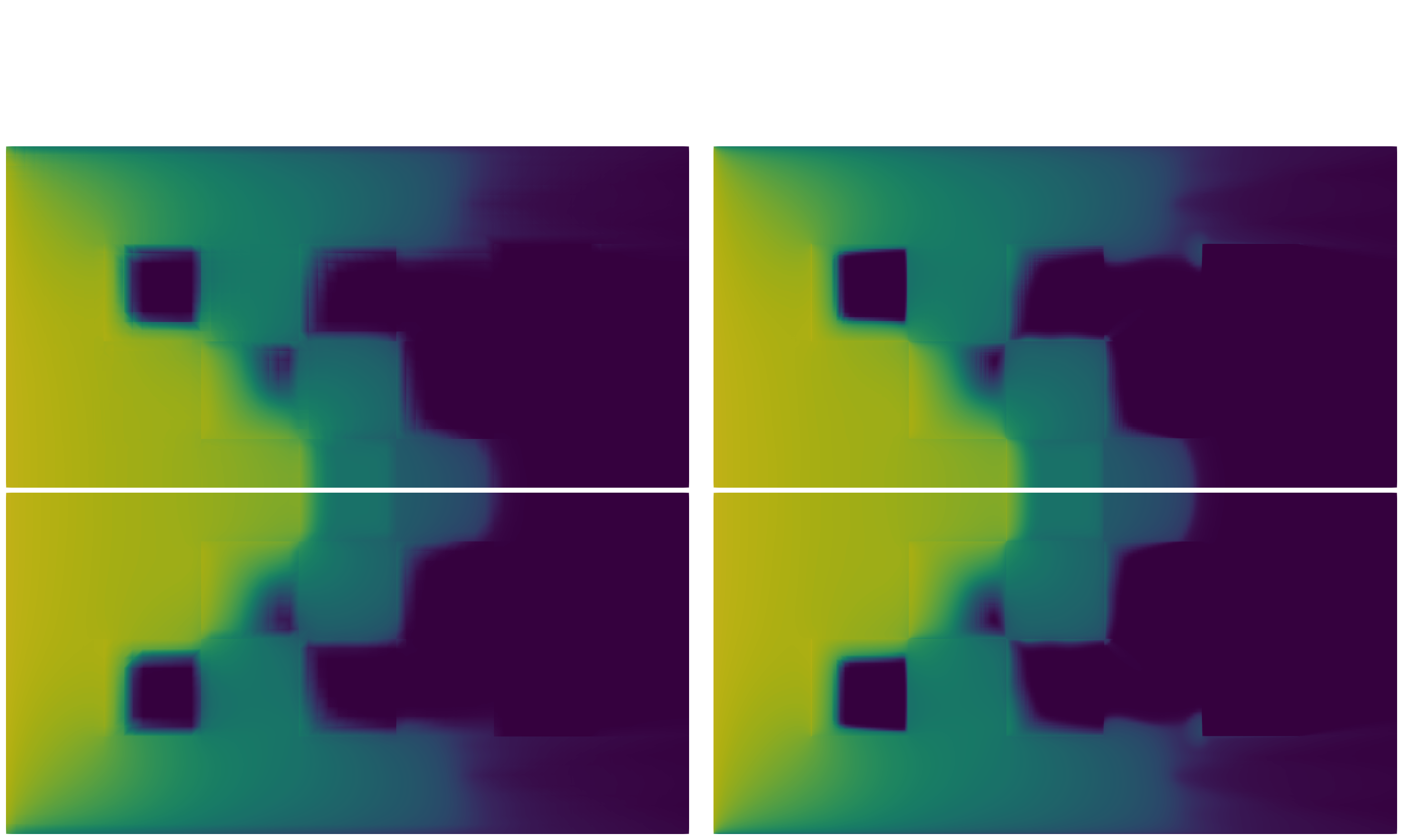}
			\put(1,1){\textcolor{white}{\scalebox{0.6}{Consistent}}}
			\put(1,26){\textcolor{white}{\scalebox{0.6}{Independent}}}
			\put(25,51){\makebox(0,0){\scalebox{0.6}{Coarse}}}

			\put(51,1){\textcolor{white}{\scalebox{0.6}{Consistent}}}
			\put(51,26){\textcolor{white}{\scalebox{0.6}{Independent}}}
			\put(75,51){\makebox(0,0){\scalebox{0.6}{Fine}}}
		\end{overpic}
		\caption{$t = \SI{2.6}{\ns}$}
	\end{subfigure}
	\quad
	\begin{subfigure}{0.45\textwidth}
		\centering
		\begin{overpic}[width=\textwidth, trim={0 0 25 250}, clip]{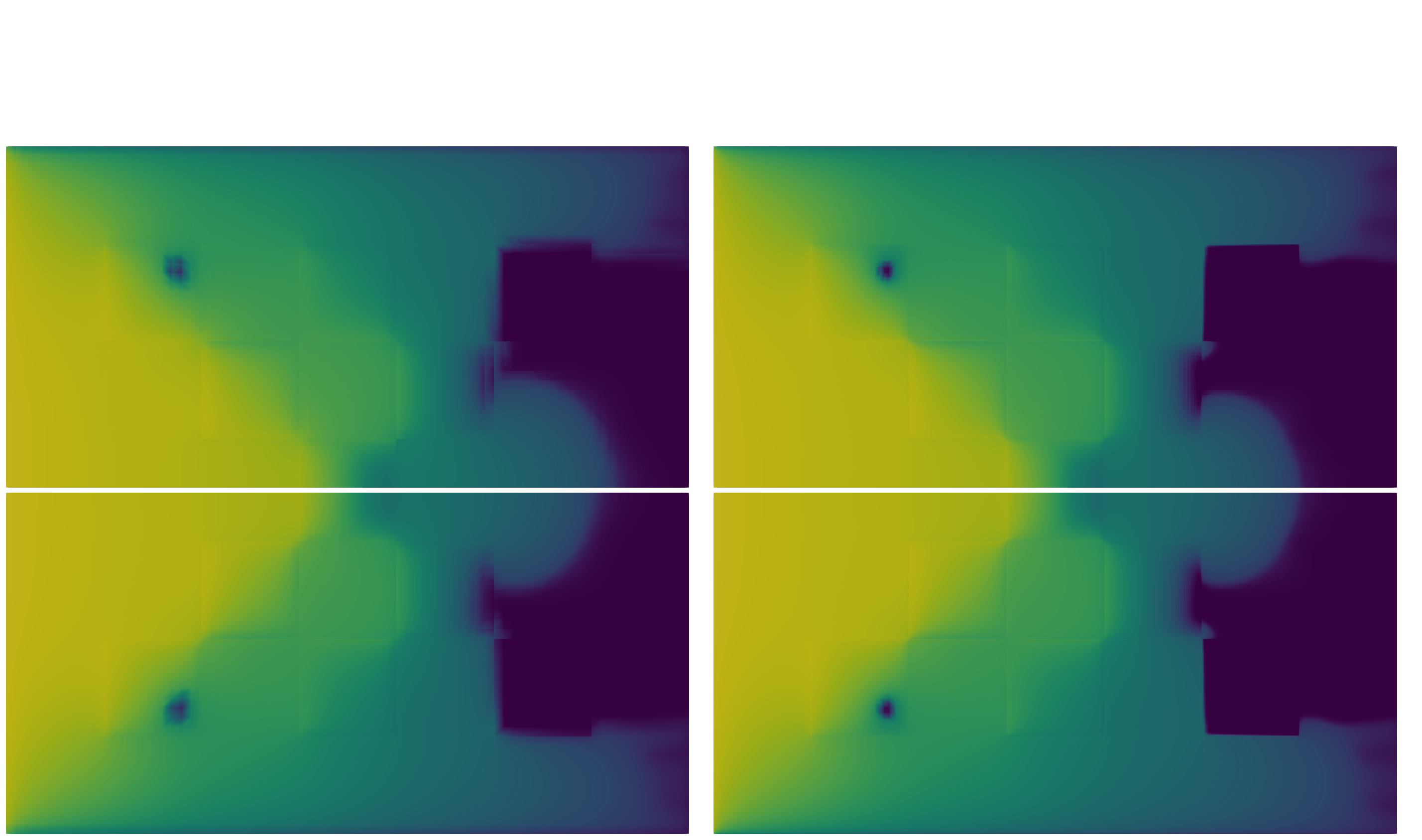}
			\put(1,1){\textcolor{white}{\scalebox{0.6}{Consistent}}}
			\put(1,26){\textcolor{white}{\scalebox{0.6}{Independent}}}
			\put(25,51){\makebox(0,0){\scalebox{0.6}{Coarse}}}

			\put(51,1){\textcolor{white}{\scalebox{0.6}{Consistent}}}
			\put(51,26){\textcolor{white}{\scalebox{0.6}{Independent}}}
			\put(75,51){\makebox(0,0){\scalebox{0.6}{Fine}}}
		\end{overpic}
		\caption{$t = \SI{10}{\ns}$}
	\end{subfigure}
	\begin{subfigure}{0.45\textwidth}
		\centering
		\begin{overpic}[width=\textwidth, trim={0 0 25 250}, clip]{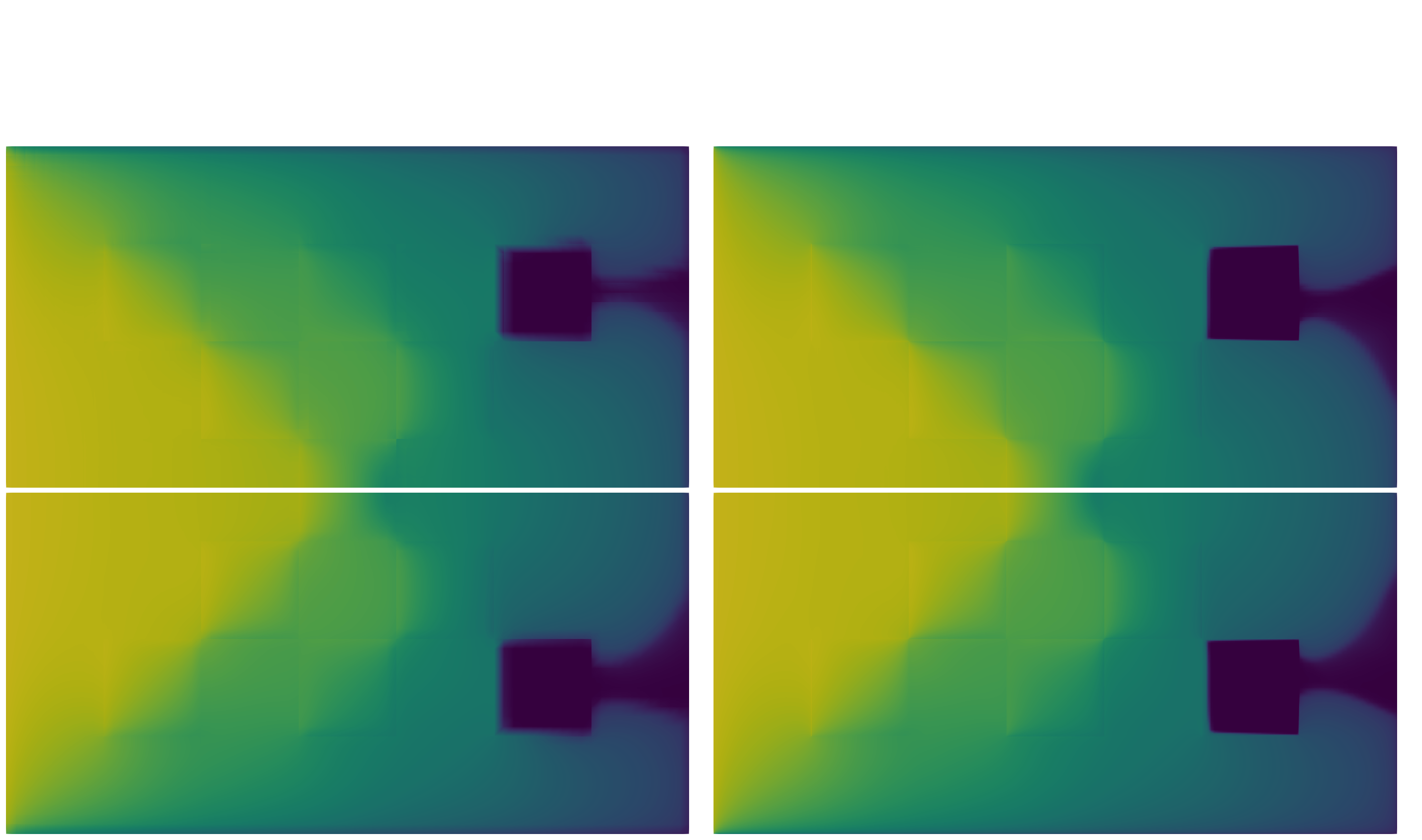}
			\put(1,1){\textcolor{white}{\scalebox{0.6}{Consistent}}}
			\put(1,26){\textcolor{white}{\scalebox{0.6}{Independent}}}
			\put(25,51){\makebox(0,0){\scalebox{0.6}{Coarse}}}

			\put(51,1){\textcolor{white}{\scalebox{0.6}{Consistent}}}
			\put(51,26){\textcolor{white}{\scalebox{0.6}{Independent}}}
			\put(75,51){\makebox(0,0){\scalebox{0.6}{Fine}}}
		\end{overpic}
		\caption{$t = \SI{15}{\ns}$}
	\end{subfigure}
	\quad
	\begin{subfigure}{0.45\textwidth}
		\centering
		\begin{overpic}[width=\textwidth, trim={0 0 25 250}, clip]{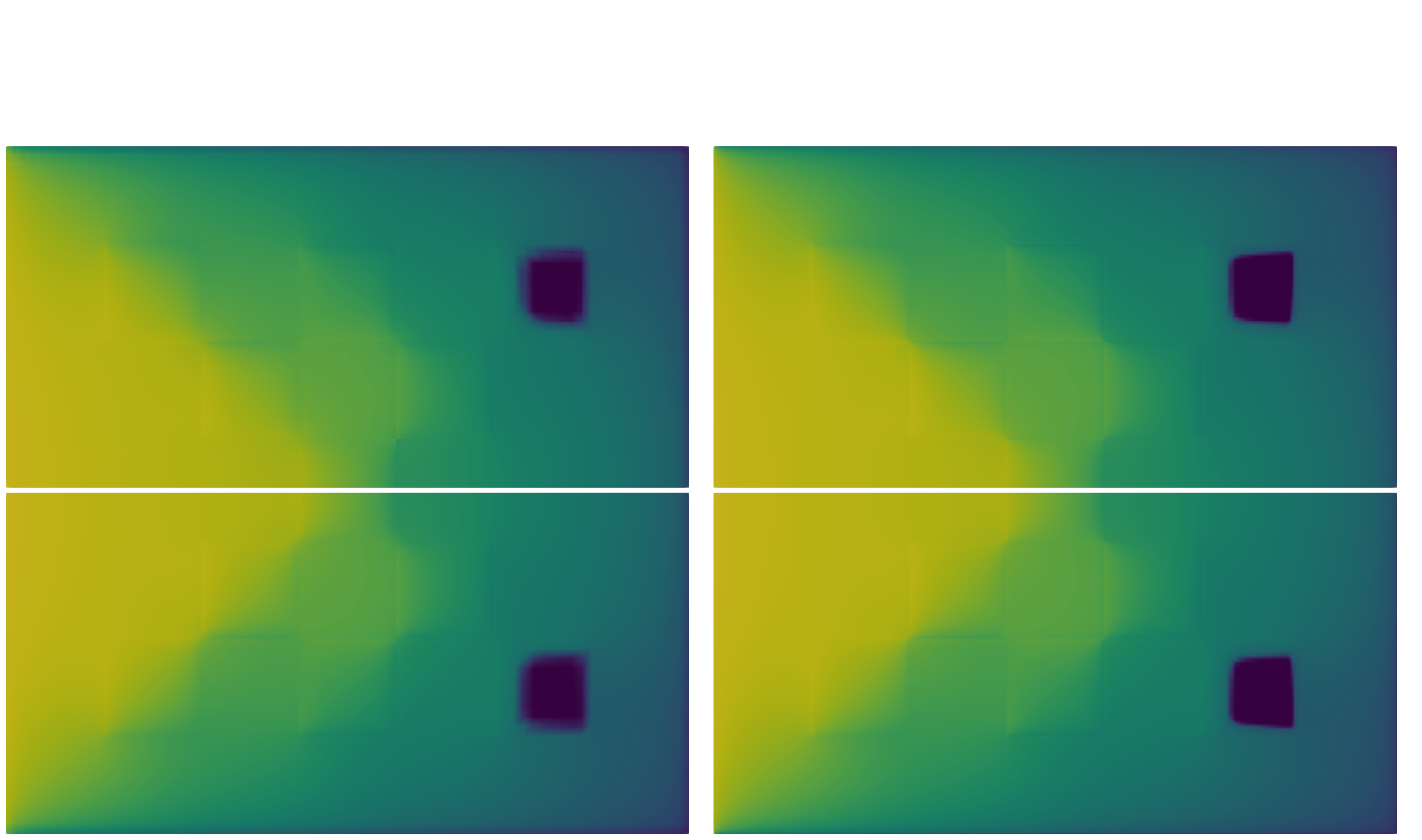}
			\put(1,1){\textcolor{white}{\scalebox{0.6}{Consistent}}}
			\put(1,26){\textcolor{white}{\scalebox{0.6}{Independent}}}
			\put(25,51){\makebox(0,0){\scalebox{0.6}{Coarse}}}

			\put(51,1){\textcolor{white}{\scalebox{0.6}{Consistent}}}
			\put(51,26){\textcolor{white}{\scalebox{0.6}{Independent}}}
			\put(75,51){\makebox(0,0){\scalebox{0.6}{Fine}}}
		\end{overpic}
		\caption{$t = \SI{50}{\ns}$}
	\end{subfigure}
	\caption{The temperature over time from the consistent and independent SM methods on the coarse and fine lattice problems. }
	\label{fig:brunner_temperature}
\end{figure}
\begin{figure}
	\centering
	\begin{subfigure}{0.45\textwidth}
		\centering
		\begin{overpic}[width=\textwidth, trim={0 0 25 250}, clip]{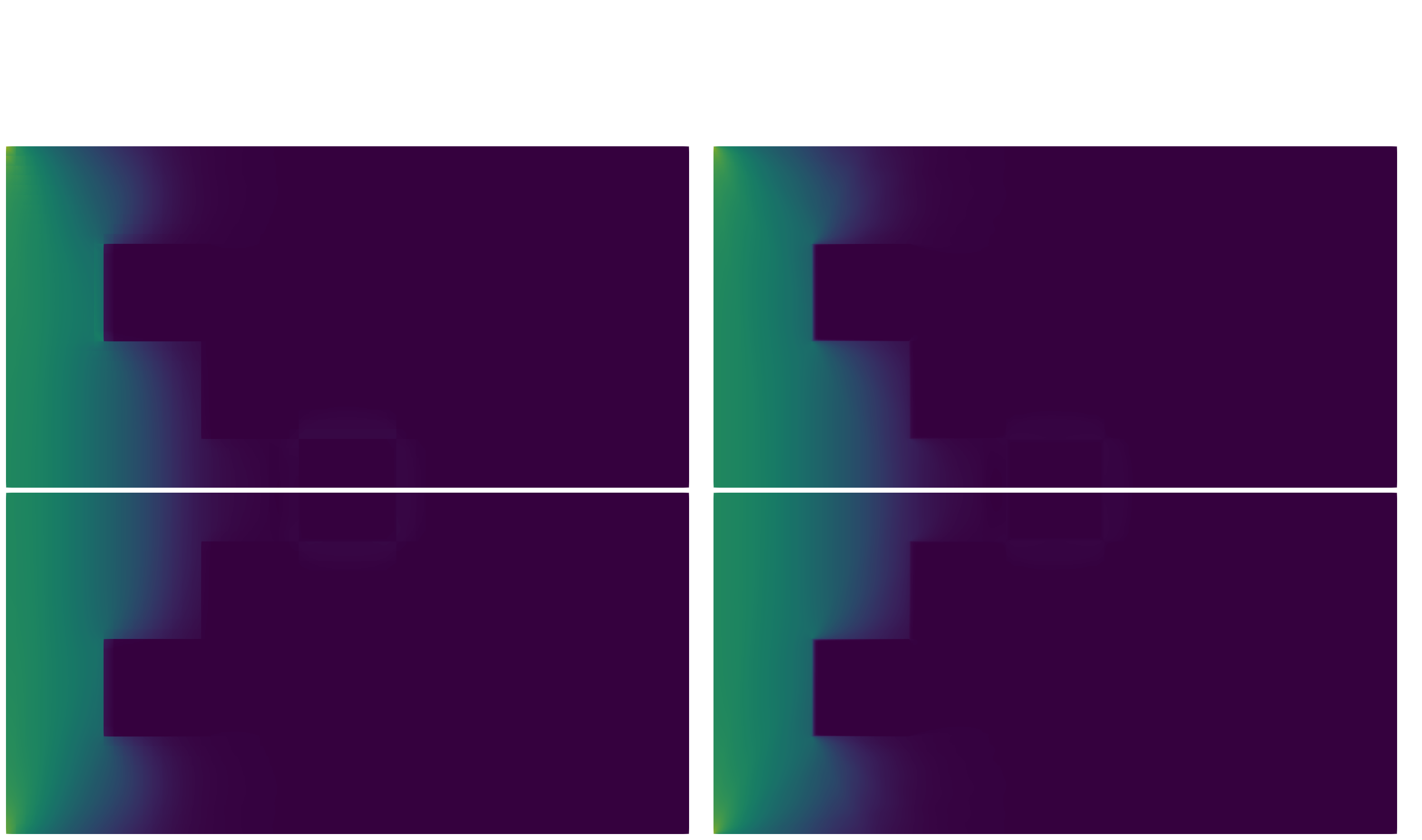}
			\put(1,1){\textcolor{white}{\scalebox{0.6}{Consistent}}}
			\put(1,26){\textcolor{white}{\scalebox{0.6}{Independent}}}
			\put(25,51){\makebox(0,0){\scalebox{0.6}{Coarse}}}

			\put(51,1){\textcolor{white}{\scalebox{0.6}{Consistent}}}
			\put(51,26){\textcolor{white}{\scalebox{0.6}{Independent}}}
			\put(75,51){\makebox(0,0){\scalebox{0.6}{Fine}}}
		\end{overpic}
		\caption{$t = \SI{0.2}{\ns}$}
	\end{subfigure}
	\quad
	\begin{subfigure}{0.45\textwidth}
		\centering
		\begin{overpic}[width=\textwidth, trim={0 0 25 250}, clip]{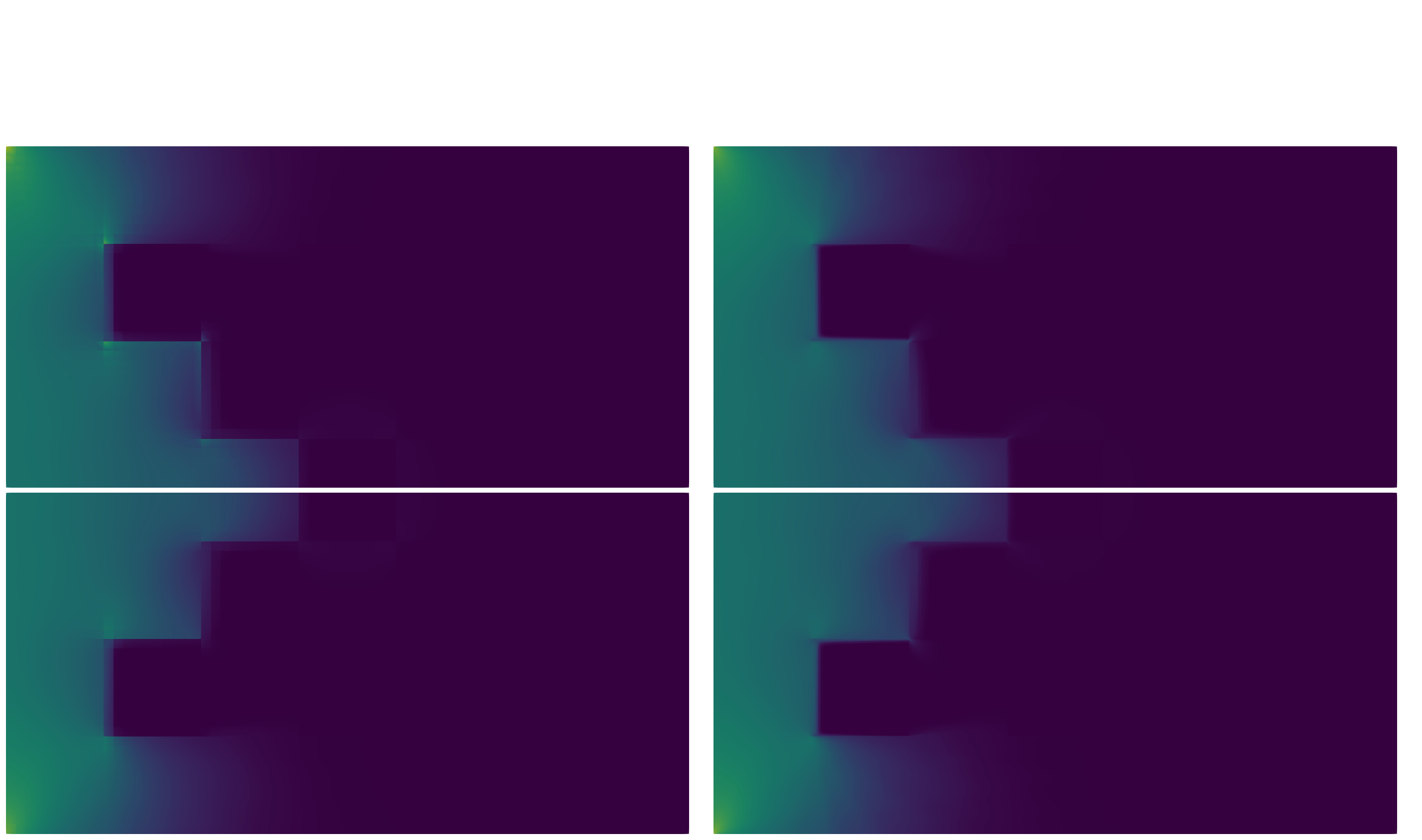}
			\put(1,1){\textcolor{white}{\scalebox{0.6}{Consistent}}}
			\put(1,26){\textcolor{white}{\scalebox{0.6}{Independent}}}
			\put(25,51){\makebox(0,0){\scalebox{0.6}{Coarse}}}

			\put(51,1){\textcolor{white}{\scalebox{0.6}{Consistent}}}
			\put(51,26){\textcolor{white}{\scalebox{0.6}{Independent}}}
			\put(75,51){\makebox(0,0){\scalebox{0.6}{Fine}}}
		\end{overpic}
		\caption{$t = \SI{0.5}{\ns}$}
	\end{subfigure}
	\begin{subfigure}{0.45\textwidth}
		\centering
		\begin{overpic}[width=\textwidth, trim={0 0 25 250}, clip]{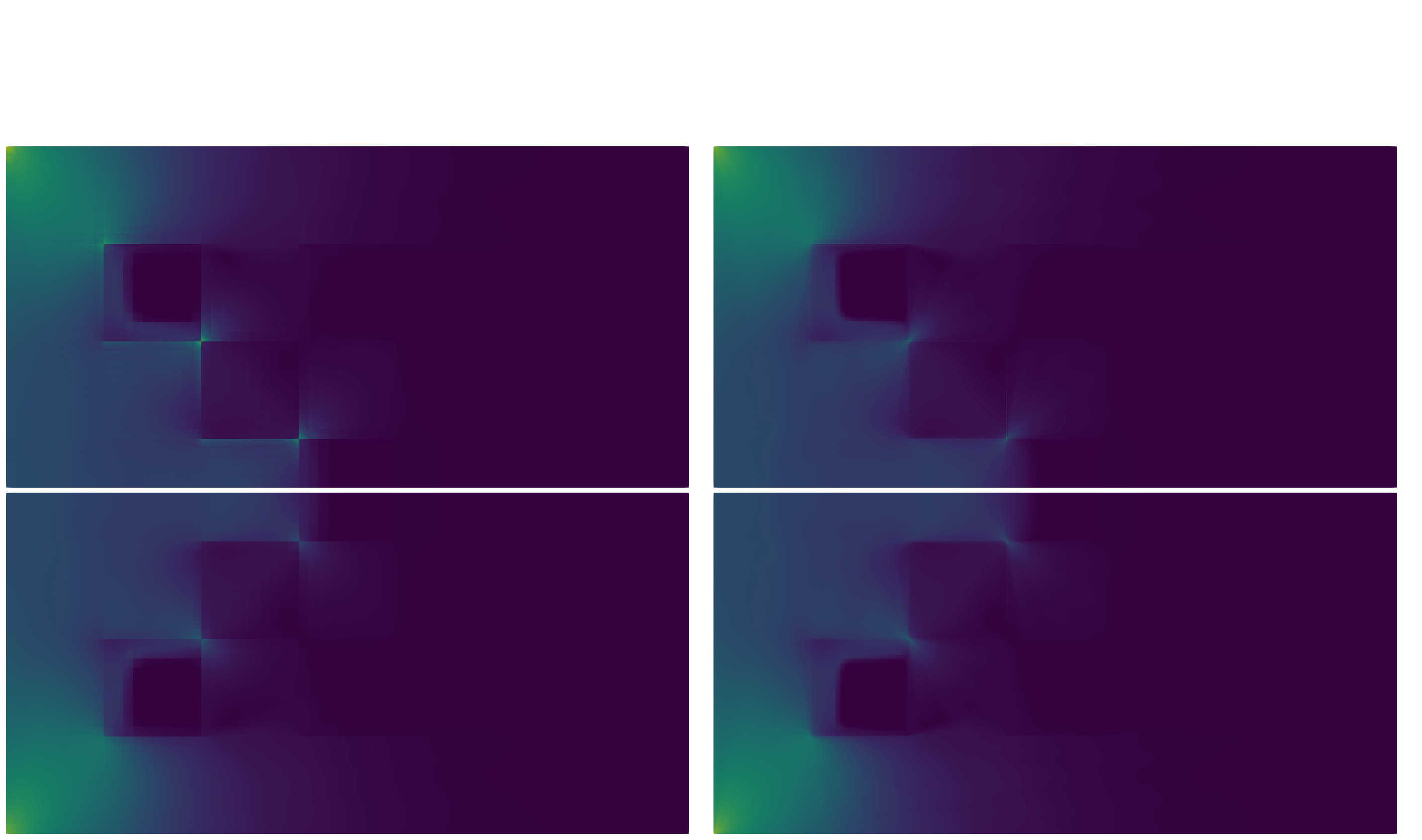}
			\put(1,1){\textcolor{white}{\scalebox{0.6}{Consistent}}}
			\put(1,26){\textcolor{white}{\scalebox{0.6}{Independent}}}
			\put(25,51){\makebox(0,0){\scalebox{0.6}{Coarse}}}

			\put(51,1){\textcolor{white}{\scalebox{0.6}{Consistent}}}
			\put(51,26){\textcolor{white}{\scalebox{0.6}{Independent}}}
			\put(75,51){\makebox(0,0){\scalebox{0.6}{Fine}}}
		\end{overpic}
		\caption{$t = \SI{2.6}{\ns}$}
	\end{subfigure}
	\quad
	\begin{subfigure}{0.45\textwidth}
		\centering
		\begin{overpic}[width=\textwidth, trim={0 0 25 250}, clip]{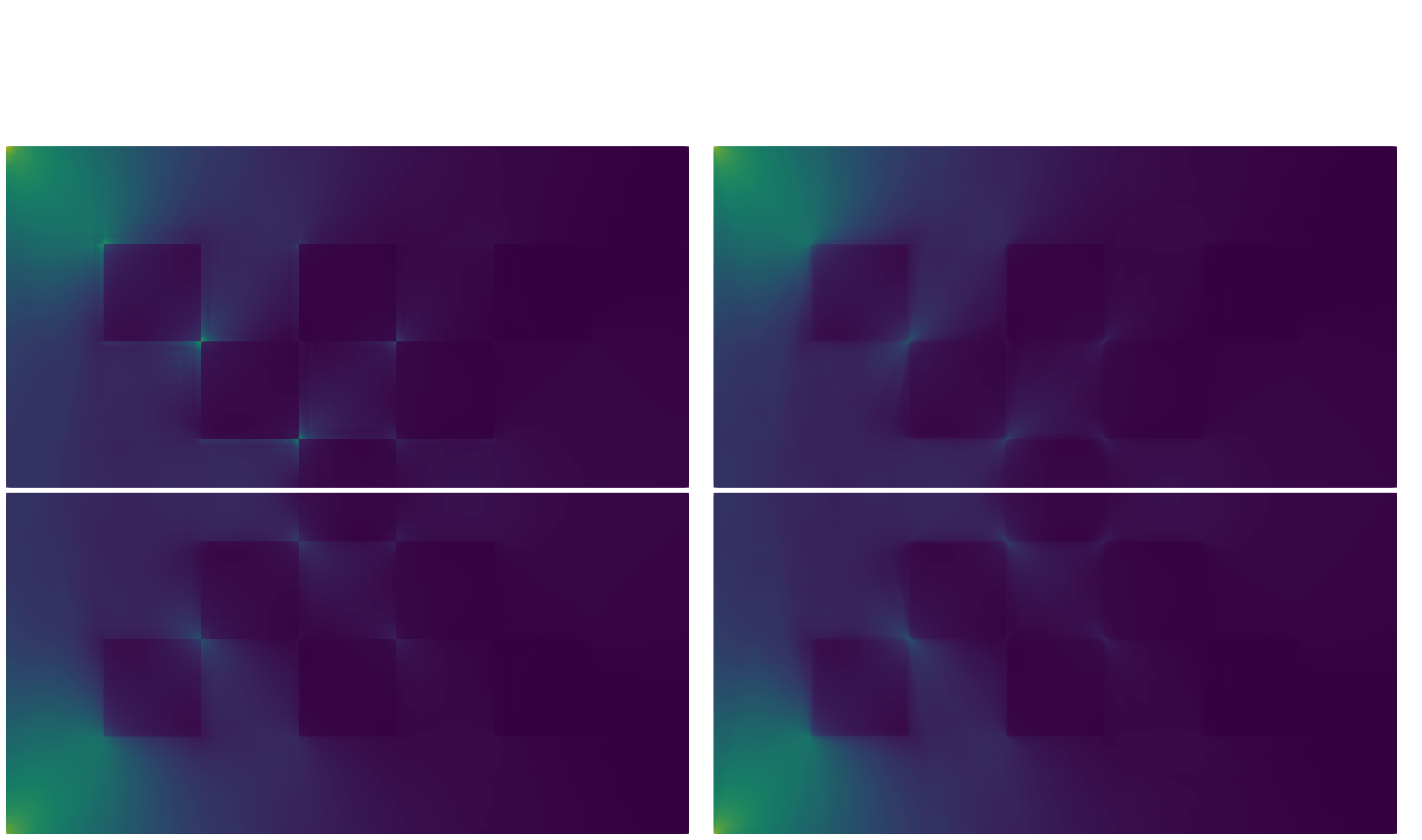}
			\put(1,1){\textcolor{white}{\scalebox{0.6}{Consistent}}}
			\put(1,26){\textcolor{white}{\scalebox{0.6}{Independent}}}
			\put(25,51){\makebox(0,0){\scalebox{0.6}{Coarse}}}

			\put(51,1){\textcolor{white}{\scalebox{0.6}{Consistent}}}
			\put(51,26){\textcolor{white}{\scalebox{0.6}{Independent}}}
			\put(75,51){\makebox(0,0){\scalebox{0.6}{Fine}}}
		\end{overpic}
		\caption{$t = \SI{50}{\ns}$}
	\end{subfigure}
	\caption{The magnitude of the flux over time from the consistent and independent SM methods on the coarse and fine lattice problems. }
	\label{fig:brunner_flux}
\end{figure}
The evolution of the temperature and magnitude of the flux are shown in Figs.~\ref{fig:brunner_temperature} and \ref{fig:brunner_flux} for the consistent and independent methods on the coarse and fine lattice problems. 
To save space, we omit the unaccelerated method from these plots as the consistent SM and unaccelerated methods produced solutions that were visually equivalent. 
Unlike the crooked pipe problem, the independent method is accurate even on the coarse lattice problem; the temperatures and fluxes are qualitatively in agreement between both methods and both resolutions. 
In particular, the flux does not have the hotspot and wave speed issues observed on the crooked pipe. 
This qualitative comparison is made quantitative with the three tracer plots in Fig.~\ref{fig:brunner_tracer} corresponding to the spatial locations depicted in Fig.~\ref{fig:brunner_diag}. 
At the final time, the consistent and independent temperatures deviate on the coarse problem by \SI{0.23}{\eV}, \SI{2.85}{\eV}, and \SI{2.52}{\eV} at locations (a), (b), and (c), respectively. 
For the fine problem, these deviations are reduced to \SI{5.97e-3}{\eV}, \SI{0.39}{\eV}, and \SI{0.36}{\eV}. 
This suggests that the coarse problem is well-resolved in space resulting in similar solution quality for both methods despite the fact that the coarse mesh used here is twice as coarse as the coarse crooked pipe problem. 
As with the crooked pipe, the unaccelerated method behaves similarly to the consistent SM solution. 
However, on this problem, the unaccelerated algorithm struggles in the thick ``cold iron'' blocks as shown in tracer location (a) where the coarse and fine solutions heat the center of the block approximately \SI{1.5}{\ns} and \SI{0.5}{\ns} earlier than the fine consistent solution, respectively. 
This is notable as both the coarse consistent and coarse independent solutions are closer to the fine mesh solution than even the fine unaccelerated solution. 
This suggests that the unaccelerated method may be struggling to locally converge in this thick iron block and that the SM algorithm provides a more robust solution in this challenging portion of the domain. 

\begin{figure}
	\centering
	\foreach \x in {0,4,8}{
		\begin{subfigure}{0.32\textwidth}
			\centering
			\includegraphics[width=\textwidth]{brun_tracer\x.pdf}
			\caption{}
		\end{subfigure}	
	}
	\caption{Evolution of the temperature over time for the consistent and independent SM methods on the coarse and fine lattice problems. Subfigures (a)--(c) correspond to the spatial locations in Fig.~\ref{fig:brunner_diag}. }
	\label{fig:brunner_tracer}
\end{figure}

\begin{figure}
	\centering
	\begin{subfigure}{0.4\textwidth}
		\centering 
		\includegraphics[width=\textwidth]{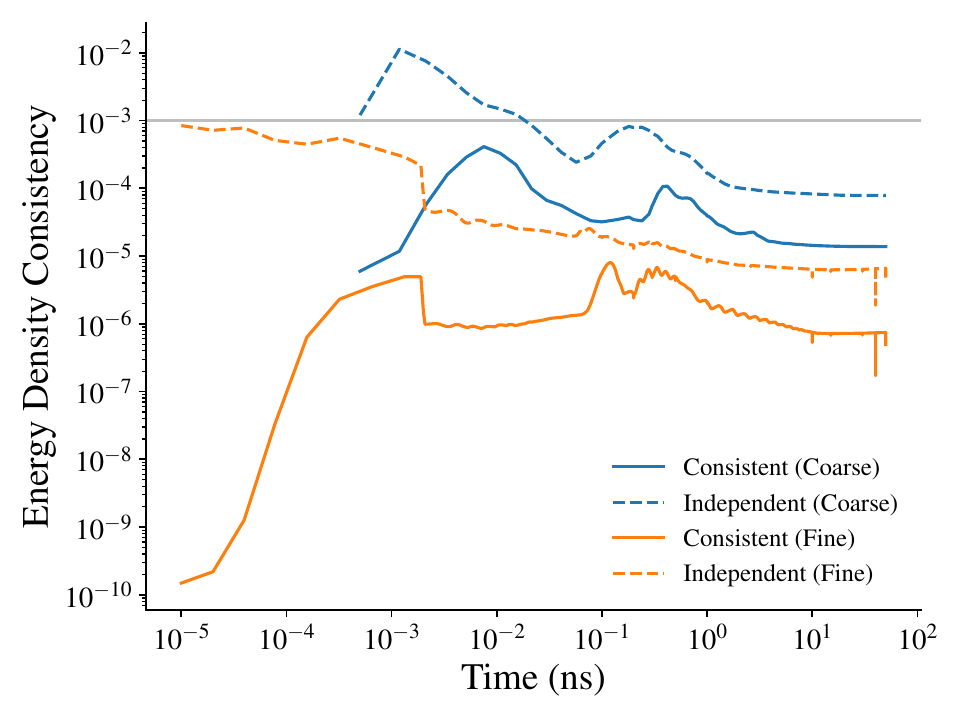}
		\caption{}
	\end{subfigure}
	\begin{subfigure}{0.4\textwidth}
		\centering 
		\includegraphics[width=\textwidth]{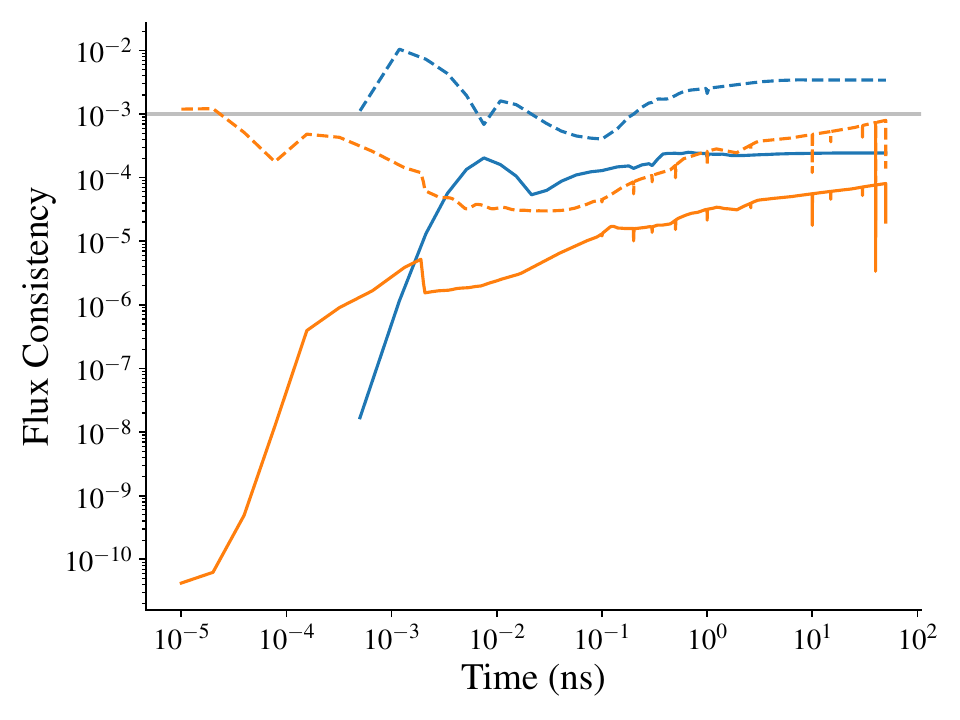}
		\caption{}
	\end{subfigure}
	\caption{The difference between the high and low-order energy density and flux as a function of simulation time for the consistent and independent SM methods on the coarse and fine lattice problems. }
	\label{fig:brunner_consistency}
\end{figure}
The consistency between the high and low-order solutions is plotted against simulation time in Fig.~\ref{fig:brunner_consistency}. 
Aside from the first few time steps on both the coarse and fine problems, the independent method is consistent with the high-order solution below the iterative tolerance of $10^{-3}$ in both the energy density and flux. 
The consistent method achieves consistency below the iterative tolerance on both problems for both moments. 
This is further evidence that the coarse problem is adequately resolved. 

\subsubsection{Performance}
\begin{figure} 
	\centering
	\begin{subfigure}{0.4\textwidth}
		\centering
		\includegraphics[width=\textwidth]{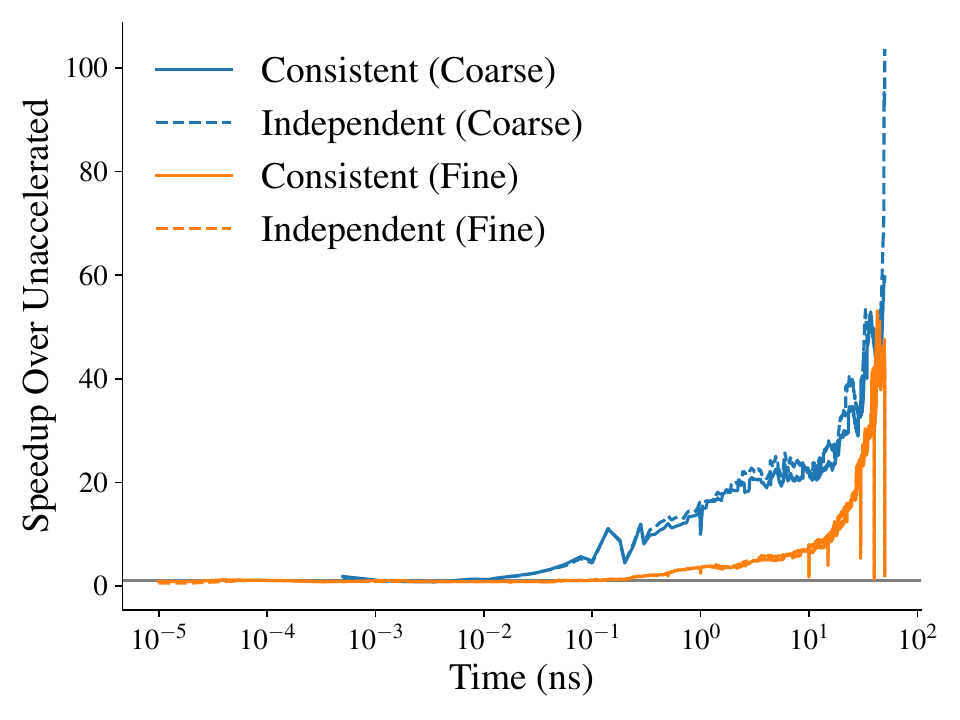}
		\caption{}
		\label{fig:brunner_performance_speedup}
	\end{subfigure}
	\begin{subfigure}{0.4\textwidth}
		\centering
		\includegraphics[width=\textwidth]{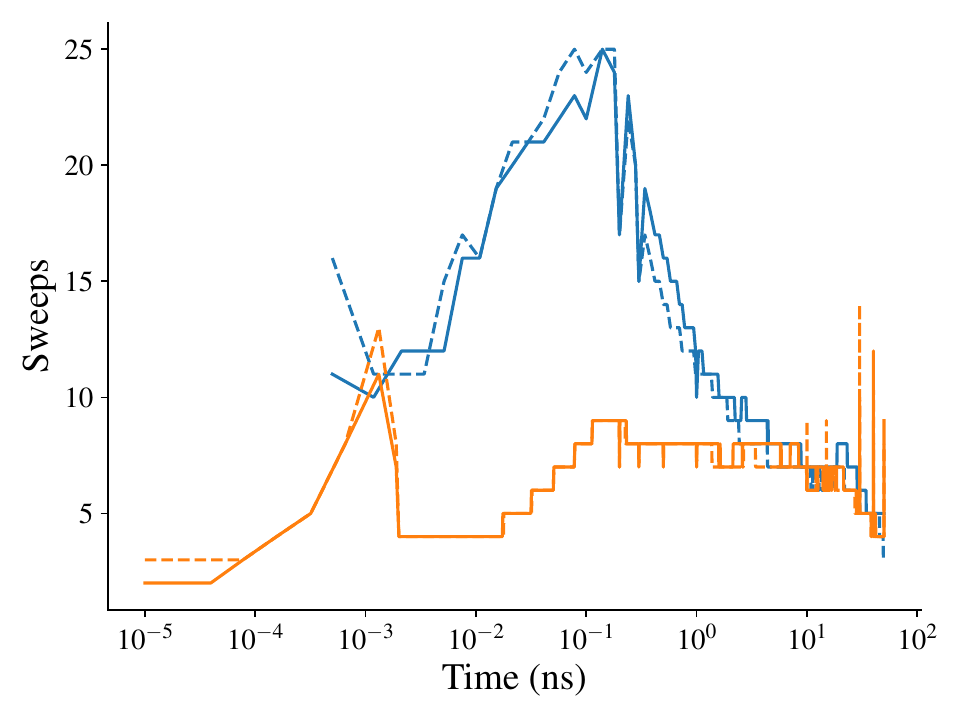}
		\caption{}
		\label{fig:brunner_performance_sweeps}
	\end{subfigure}
	\begin{subfigure}{0.4\textwidth}
		\centering
		\includegraphics[width=\textwidth]{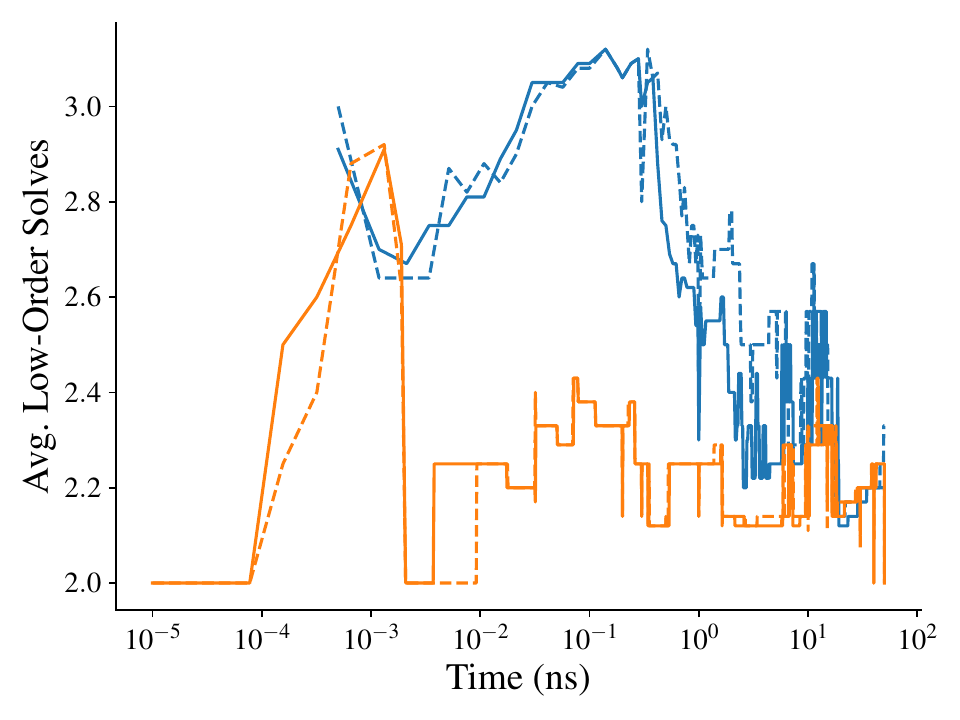}
		\caption{}
		\label{fig:brunner_performance_nonlinear}
	\end{subfigure}
	\begin{subfigure}{0.4\textwidth}
		\centering
		\includegraphics[width=\textwidth]{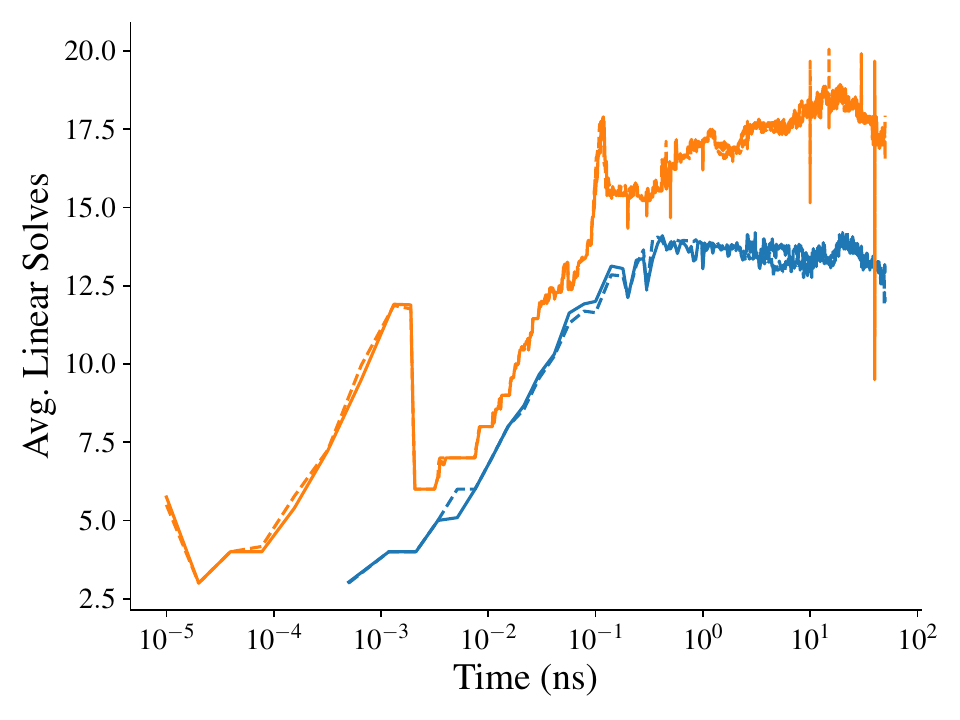}
		\caption{}
		\label{fig:brunner_performance_linear}
	\end{subfigure}
	\caption{Performance metrics as a function of simulation time for the consistent and independent SM methods on the coarse and fine lattice problems.}
	\label{fig:brunner_performance} 
\end{figure}
We conclude by presenting performance metrics. 
Figure \ref{fig:brunner_performance_speedup} shows the speedup against the unaccelerated method at each time step. 
Summaries of the iterative and computational performance of the sweep, low-order Newton solve, low-order linear solve, and nonlinear temperature elimination and the sweep savings and speedup with respect to the unaccelerated method are provided in Tables \ref{tab:it_summary_brunner} and \ref{tab:speedup_brunner}, respectively. 
At early times, the SM methods break even against the unaccelerated method achieving similar cycle costs. 
Recall that on the crooked pipe problem, the SM methods were slower than unaccelerated at early times. 
This change in behavior is due to the multigroup sweep dominating the cost of the gray low-order problem making the cost of the ``extra'' acceleration steps contribute relatively less to overall cost. 
As absorption-emission becomes stiff, the SM methods exhibit extreme speedups over unaccelerated with independent $103\times$ and consistent $60\times$ faster. 
Aside from the large discrepancy near the end of the simulation, consistent and independent show nearly identical speedups. 
In total time-to-solution, independent was $32\times$ faster and consistent $30\times$ with respect to the unaccelerated scheme's total runtime, respectively.

On the fine problem, the SM methods have the same cost as the unaccelerated method for much of the simulation runtime. 
The fine problem takes very small time steps allowing both the unaccelerated and accelerated schemes to converge rapidly. 
However, as the time step increases and the problem approaches the equilibrium diffusion limit, independent and consistent achieve speedups of $52\times$ and $53\times$, respectively. 
Due to the small time steps, the SM methods had lower overall speedups of $11.4\times$ and $12\times$ for consistent and independent, respectively. 
The sweeps to convergence, average nonlinear low-order iterations per outer iteration, and average low-order linear iterations per outer iteration versus time are shown in Figs.~\ref{fig:brunner_performance_sweeps}-\ref{fig:brunner_performance_linear}. 
In all cases, both SM methods converge at similar rates for outer, inner, and linear iterations. 
\begin{table}
	\centering
	\caption{Summary of iterative performance and cost of the sweep, low-order nonlinear solve, low-order linear solve, and nonlinear temperature elimination for the consistent SM, independent SM, and unaccelerated algorithm on the coarse and fine lattice problems.}
	\label{tab:it_summary_brunner}
	\begin{tabular}{cccccccccc}
\toprule
 &  & \multicolumn{3}{c}{Coarse}  &  & \multicolumn{3}{c}{Fine} \\
\cmidrule{3-5}\cmidrule{7-9}
 & Metric & Consistent & Independent & Unaccelerated & & Consistent & Independent & Unaccelerated \\
\midrule
\multirow{6}{*}{\rotatebox{90}{Sweeps}} & Min & 5 & 3 & 13 & & 2 & 3 & 2 \\
 & Max & 25 & 25 & 417 & & 11 & 14 & 225 \\
 & Avg. & 6.90 & 6.24 & 278.31 & & 6.48 & 6.13 & 83.14 \\
 & Total & \num{8713} & \num{7878} & \num{351511} & & \num{30601} & \num{28937} & \num{392403} \\
 & Cost (s) & \num{876.5} & \num{816.9} & \num{35817.8} & & \num{12245.5} & \num{11661.5} & \num{169378.2} \\
 & Avg. Cost (ms) & \num{100.6} & \num{103.7} & \num{101.9} & & \num{400.2} & \num{403.0} & \num{431.6} \\
\addlinespace
\multirow{6}{*}{\rotatebox{90}{LO Nonlinear}} & Min & 2 & 2 & -- & & 2 & 2 & -- \\
 & Max & 5 & 5 & -- & & 3 & 3 & -- \\
 & Avg. & 2.27 & 2.30 & -- & & 2.20 & 2.21 & -- \\
 & Total & \num{20054} & \num{18528} & -- & & \num{67276} & \num{64045} & -- \\
 & Cost (s) & \num{473.1} & \num{416.9} & -- & & \num{3458.4} & \num{3246.1} & -- \\
 & Avg. Cost (ms) & \num{54.3} & \num{52.9} & -- & & \num{113.0} & \num{112.2} & -- \\
\addlinespace
\multirow{6}{*}{\rotatebox{90}{LO Linear}} & Min & 3 & 3 & -- & & 3 & 3 & -- \\
 & Max & 15 & 16 & -- & & 22 & 21 & -- \\
 & Avg. & 13.29 & 13.16 & -- & & 17.41 & 17.31 & -- \\
 & Total & \num{266576} & \num{243880} & -- & & \num{1170986} & \num{1108735} & -- \\
 & Cost (s) & \num{166.3} & \num{150.6} & -- & & \num{1207.5} & \num{1162.4} & -- \\
 & Avg. Cost (ms) & \num{8.3} & \num{8.1} & -- & & \num{17.9} & \num{18.1} & -- \\
\addlinespace
\multirow{6}{*}{\rotatebox{90}{Energy Balance}} & Min & 1 & 1 & 4 & & 1 & 1 & 4 \\
 & Max & $40^*$ & $40^*$ & 8 & & $40^*$ & $40^*$ & 5 \\
 & Avg. & 8.16 & 8.64 & 4.51 & & 11.24 & 10.84 & 4.53 \\
 & Total & \num{163680} & \num{160033} & \num{1583855} & & \num{755998} & \num{694191} & \num{1778154} \\
 & Cost (s) & \num{77.0} & \num{70.7} & \num{6207.3} & & \num{574.1} & \num{560.5} & \num{14620.1} \\
 & Avg. Cost (ms) & \num{3.8} & \num{3.8} & \num{17.7} & & \num{8.5} & \num{8.8} & \num{37.3} \\
\addlinespace
 & Wall Time (s) & \num{1479.4} & \num{1354.9} & \num{44133.8} & & \num{16718.3} & \num{15818.0} & \num{190405.5} \\
\bottomrule
\end{tabular}
	$^*$ indicates maximum iterations reached
\end{table}
\begin{table}
	\centering
	\caption{The min, max, average, and total speedups and cost in cycle time over the unaccelerated scheme on the lattice problem. }
	\label{tab:speedup_brunner}
	\begin{tabular}{cccccccc}
\toprule
 &  & \multicolumn{2}{c}{Coarse}  &  & \multicolumn{2}{c}{Fine} \\
\cmidrule{3-4}\cmidrule{6-7}
 & Metric & Consistent & Independent & & Consistent & Independent \\
\midrule
\multirow{4}{*}{\rotatebox{90}{Sweeps}} & Min & 1.08 & 1.13 & & 0.80 & 0.67 \\
 & Max & 83.40 & 139.00 & & 56.25 & 56.25 \\
 & Avg. & 44.47 & 50.03 & & 15.01 & 15.60 \\
 & Total & 40.34 & 44.62 & & 12.82 & 13.56 \\
\addlinespace
\multirow{4}{*}{\rotatebox{90}{Cost}} & Min & 0.78 & 0.84 & & 0.68 & 0.57 \\
 & Max & 59.83 & 103.62 & & 53.13 & 51.98 \\
 & Avg. & 32.96 & 36.58 & & 13.42 & 13.95 \\
 & Total & 29.83 & 32.57 & & 11.39 & 12.04 \\
\bottomrule
\end{tabular}
\end{table}

\section{Conclusions} \label{sec:conclusions}
We have extended the consistent and independent local Discontinuous Galerkin (LDG)-based Second Moment (SM) methods from Olivier et~al.\,\cite{csmm} and Olivier and Haut \cite{smm_olivier} developed for steady-state, linear transport problems to the time-dependent, multi-frequency thermal radiative transfer (TRT) equations. 
The high-order system is discretized with an upwind discontinuous Galerkin discretization in space, multigroup in frequency, discrete ordinates in direction, and semi-implicit backward Euler in time. 
We consider two spatial discretizations for the SM low-order problem: a consistent method where the spatial discretization is designed to be algebraically consistent with the moments of the high-order transport equation and an independent method where the low-order discretization is chosen independently of the high-order equation. 
Both discretizations are based on the LDG discretization of radiation diffusion and were designed to leverage existing linear solver technology to efficiently invert the low-order system. 
Use of LDG for both the consistent and independent methods ensures that the consistent and independent approaches are compared fairly. 

The algorithms use a \emph{gray} SM low-order problem to accelerate convergence of the stiff, nonlinear radiation-material coupling. 
The high-order transport equation is solved to compute closures for the low-order SM system. 
A Newton-type iteration is combined with local, nonlinear elimination of the temperature to nonlinearly resolve the coupled low-order and energy balance equations. 
Since the low-order system is independent of both direction and frequency, it is practical to compute and invert the Jacobian of the low-order and energy balance equations. 
The collapse to a gray low-order problem is achieved with a novel, \emph{additive} correction that makes the frequency-integrated high-order and gray low-order problems consistent. 
An additive correction allows the correction term to contribute only to the right-hand side, preserving the symmetric and positive definite structure of the left-hand-side diffusion operator.  
Furthermore, we show that when the spatial discretization is lumped, the linearization of the energy balance equation results in a small diagonal perturbation to the diffusion system. 
Thus, absorption-emission is resolved nonlinearly while only requiring the inversion of gray, symmetric, and positive definite diffusion systems. 
Due to this, the Jacobian can be inverted efficiently with algebraic multigrid-preconditioned conjugate gradient. 
The temperature from the simultaneous low-order and energy balance solve is used to compute a new emission source for the high-order equation and this process is repeated until convergence. 

Space-time accuracy studies were conducted on gray and multigroup, one-dimensional problems. 
Using a spatially and temporally resolved simulation as a discrete reference, the consistent and independent methods were shown to converge with first-order accuracy in both space and time with some degraded convergence seen on the more challenging multigroup problem. 
We postulate that the first-order spatial convergence is reduced from the expected second-order convergence of the spatial discretization due to the use of piecewise-constant opacity data or the lack of regularity in the solution limiting the maximum convergence rate to first-order (e.g.~see the steep temperature gradient observed in Fig.~\ref{fig:larsen_opacity_a}). 
As seen in \cite{csmm}, the consistent scheme was more accurate but also less iteratively efficient than the independent method. 
However, both methods were similar in accuracy per cost suggesting that the increased cost of the consistent method is sufficiently balanced by improved accuracy. 
Gray acceleration remained effective across a range of time step sizes on the multigroup problem, resulting in the SM methods requiring up to $10\times$ fewer sweeps than an unaccelerated scheme that does not use a low-order diffusion system to accelerate absorption-emission. 
We also showed that optimal convergence in time is achieved when using either the high or low-order solution as the previous time step's moment solution in the temporal integration of the low-order system. 
In addition, piecewise-constant \emph{gray} opacities were not sufficient to preserve solution quality on multigroup problems even though the \emph{multigroup} opacity was represented as piecewise-constant. 
This is due to the spatial variance of the weight functions used to collapse the multigroup opacity. 
We found that representing the gray opacities as piecewise-linear was sufficient to preserve the solution quality delivered by the high-order discretization in isolation. 

Solution quality and performance of the schemes were investigated on gray crooked pipe and multigroup lattice problems in two spatial dimensions. 
Performance was grouped into two behaviors depending on the mesh resolution. 
On spatially under-resolved meshes, the independent method was up to $2\times$ faster than the consistent method. 
However, this increased performance was paired with poor solution quality. 
On the under-resolved crooked pipe problem, the independent method produced non-physical cold spots and heavily relied on energy-conservation-violating floors of the temperature to prevent crashes which were not required by the consistent method. 
We note that this behavior was only observed on the severely under-resolved crooked pipe problem where the consistent method was also inaccurate compared to resolved simulations. 
This behavior persisted even when the angular resolution was increased from $S_6$ to $S_{12}$ suggesting that spatial error is responsible for the reduced solution quality of the independent method. 
On spatially resolved meshes, the cost and solution quality of the two schemes converged with independent marginally more efficient.
The SM methods were between $12\times$ and $30\times$ faster than the unaccelerated algorithm indicating the additional low-order acceleration computations are effective. 
The iterative efficiencies of the inner, Newton-type iteration and the linear solvers were not dependent on the choice of low-order diffusion discretization, both achieving robust performance and scalability with problem size. 

Through this study, we have found that the behavior of the consistent and independent methods established for linear transport in \cite{csmm} largely extends to TRT: the consistent method was more accurate and more expensive than the independent method but the increased cost was generally commensurate with the associated improvement in accuracy and robustness.
Thus, of the schemes we have studied here, we recommend the consistent method. 
The consistent scheme did not exhibit the solution quality issues we observed on the spatially under-resolved crooked pipe problem while being capable of near-equivalent performance to the independent method on spatially resolved problems. 
It is possible that improvements to the independent method's low-order diffusion discretization -- such as through the design of numerical fluxes tailored for TRT problems, use of a low-order discretization with higher order-of-accuracy than the high-order system, or the use of slope limiting -- could improve accuracy or mitigate oscillations on spatially under-resolved problems, resulting in a more robust independent method that still retains its performance advantage over the consistent method. 
In particular, analysis of the discretization-dependent correction terms in the consistent scheme may indicate the deficiencies of the LDG-based independent method. 
Note that, while we found that consistency improves accuracy and solution quality, this is likely only the case for high-fidelity discretizations of the transport equation such as the upwind DG method used here. 
In the future, we plan to compare the performance of the consistent scheme to existing efficient algorithms for TRT and extend these methods to radiation-hydrodynamics algorithms.

\section*{Data Availability}
No external data sets were used in this study. 
Numerical results were generated with the authors' open-source code \emph{allium} \cite{allium}. 
The data supporting the findings of this study can be reproduced using this code and are therefore not separately archived. 

\section*{Acknowledgments}
S.O. was supported by the U.S. Department of Energy as a Nicholas C. Metropolis Fellow under the Laboratory Directed Research and Development program of the Los Alamos National Laboratory.
This research used resources provided by the Darwin testbed at Los Alamos National Laboratory (LANL) which is funded by the Computational Systems and Software Environments subprogram of LANL's Advanced Simulation and Computing program. 
LANL is operated by Triad National Security, LLC, for the National Nuclear Security Administration of the U.S. Department of Energy (Contract No.~89233218CNA000001). 

\bibliographystyle{IEEEtranN}
\bibliography{references}

\begin{thebibliography}{39}
\providecommand{\natexlab}[1]{#1}
\providecommand{\url}[1]{#1}
\csname url@samestyle\endcsname
\providecommand{\newblock}{\relax}
\providecommand{\bibinfo}[2]{#2}
\providecommand{\BIBentrySTDinterwordspacing}{\spaceskip=0pt\relax}
\providecommand{\BIBentryALTinterwordstretchfactor}{4}
\providecommand{\BIBentryALTinterwordspacing}{\spaceskip=\fontdimen2\font plus
\BIBentryALTinterwordstretchfactor\fontdimen3\font minus
  \fontdimen4\font\relax}
\providecommand{\BIBforeignlanguage}[2]{{%
\expandafter\ifx\csname l@#1\endcsname\relax
\typeout{** WARNING: IEEEtranN.bst: No hyphenation pattern has been}%
\typeout{** loaded for the language `#1'. Using the pattern for}%
\typeout{** the default language instead.}%
\else
\language=\csname l@#1\endcsname
\fi
#2}}
\providecommand{\BIBdecl}{\relax}
\BIBdecl

\bibitem[Mihalas and Mihalas(1984)]{mihalas1999foundations}
D.~Mihalas and B.~W. Mihalas, \emph{Foundations of Radiation
  Hydrodynamics}.\hskip 1em plus 0.5em minus 0.4em\relax Oxford University
  Press, New York, 1984.

\bibitem[{V. Ya. Gol'din}(1964)]{goldin}
{V. Ya. Gol'din}, ``A quasi-diffusion method of solving the kinetic equation,''
  \emph{USSR Computational Mathematics and Mathematical Physics}, vol.~4, pp.
  136--149, 1964.

\bibitem[Chac\'{o}n et~al.(2017)Chac\'{o}n, Chen, Knoll, Newman, Park, Taitano,
  Willert, and Womeldorff]{CHACON201721}
\BIBentryALTinterwordspacing
L.~Chac\'{o}n, G.~Chen, D.~A. Knoll, C.~Newman, H.~Park, W.~Taitano, J.~A.
  Willert, and G.~Womeldorff, ``Multiscale high-order/low-order ({HOLO})
  algorithms and applications,'' \emph{Journal of Computational Physics}, vol.
  330, pp. 21--45, 2017. [Online]. Available:
  \url{https://doi-org.libproxy.berkeley.edu/10.1016/j.jcp.2016.10.069}
\BIBentrySTDinterwordspacing

\bibitem[Lewis and {W. Miller Jr.}(1976)]{lewis_miller}
E.~Lewis and {W. Miller Jr.}, ``A comparison of p$_1$ synthetic acceleration
  techniques,'' in \emph{Transactions of the American Nuclear Society 23},
  1976.

\bibitem[Anistratov et~al.(1996)Anistratov, Aristova, and
  Gol'din]{anistratov1996nonlinear}
D.~Y. Anistratov, E.~N. Aristova, and V.~Y. Gol'din, ``A nonlinear method for
  solving the problems of radiation transfer in medium,'' \emph{Matematicheskoe
  modelirovanie}, vol.~8, no.~12, pp. 3--28, 1996.

\bibitem[Hardy et~al.(2025)Hardy, Morel, and Vermaak]{Hardy03042025}
Z.~K. Hardy, J.~E. Morel, and J.~I.~C. Vermaak, ``A second moment method for
  k-eigenvalue acceleration with continuous diffusion and discontinuous
  transport discretizations,'' \emph{Nuclear Science and Engineering}, vol.
  199, no.~4, pp. 599--612, 2025.

\bibitem[Anistratov(2021)]{anistratov2020nonlinear}
\BIBentryALTinterwordspacing
D.~Y. Anistratov, ``Nonlinear iterative projection methods with multigrid in
  photon frequency for thermal radiative transfer,'' \emph{Journal of
  Computational Physics}, vol. 444, p. 110568, 2021. [Online]. Available:
  \url{https://www.sciencedirect.com/science/article/pii/S0021999121004630}
\BIBentrySTDinterwordspacing

\bibitem[Anistratov(2019)]{ANISTRATOV2019186}
------, ``Stability analysis of a multilevel quasidiffusion method for thermal
  radiative transfer problems,'' \emph{Journal of Computational Physics}, vol.
  376, pp. 186--209, 2019.

\bibitem[Ghassemi and Anistratov(2020)]{GHASSEMI2020109315}
P.~Ghassemi and D.~Y. Anistratov, ``Multilevel quasidiffusion method with
  mixed-order time discretization for multigroup thermal radiative transfer
  problems,'' \emph{Journal of Computational Physics}, vol. 409, p. 109315,
  2020.

\bibitem[Coale and Anistratov(2023)]{COALE2023108458}
\BIBentryALTinterwordspacing
J.~M. Coale and D.~Y. Anistratov, ``Reduced order models for thermal radiative
  transfer problems based on moment equations and data-driven approximations
  of the eddington tensor,'' \emph{Journal of Quantitative Spectroscopy and
  Radiative Transfer}, vol. 296, p. 108458, 2023. [Online]. Available:
  \url{https://www.sciencedirect.com/science/article/pii/S0022407322003934}
\BIBentrySTDinterwordspacing

\bibitem[Yee et~al.(2021)Yee, Olivier, Southworth, Holec, and Haut]{yee_mc21}
B.~Yee, S.~Olivier, B.~Southworth, M.~Holec, and T.~Haut, ``A new scheme for
  solving high-order {DG} discretizations of thermal radiative transfer using
  the variable {E}ddington factor method,'' in \emph{Proceedings of the
  International Conference on Mathematics and Computational Methods Applied to
  Nuclear Science and Engineering (M\&C 2021)}, 2021.

\bibitem[Southworth et~al.(2024)Southworth, Olivier, Park, and
  Buvoli]{SOUTHWORTH2024113349}
\BIBentryALTinterwordspacing
B.~S. Southworth, S.~Olivier, H.~Park, and T.~Buvoli, ``One-sweep moment-based
  semi-implicit-explicit integration for gray thermal radiation transport,''
  \emph{Journal of Computational Physics}, vol. 517, p. 113349, 2024. [Online].
  Available:
  \url{https://www.sciencedirect.com/science/article/pii/S0021999124005977}
\BIBentrySTDinterwordspacing

\bibitem[Lou and Morel(2021)]{LOU2021110393}
\BIBentryALTinterwordspacing
J.~Lou and J.~E. Morel, ``A variable eddington factor method with different
  spatial discretizations for the radiative transfer equation and the
  hydrodynamics/radiation-moment equations,'' \emph{Journal of Computational
  Physics}, vol. 439, p. 110393, 2021. [Online]. Available:
  \url{https://www.sciencedirect.com/science/article/pii/S0021999121002886}
\BIBentrySTDinterwordspacing

\bibitem[Wollaber et~al.(2017)Wollaber, Park, Lowrie, Rauenzahn, and
  Cleveland]{Wollaber02012017}
A.~B. Wollaber, H.~Park, R.~B. Lowrie, R.~M. Rauenzahn, and M.~A. Cleveland,
  ``Multigroup radiation hydrodynamics with a high-order–low-order method,''
  \emph{Nuclear Science and Engineering}, vol. 185, no.~1, pp. 117--129, 2017.

\bibitem[Olivier et~al.(2025)Olivier, Southworth, Warsa, and Park]{csmm}
S.~Olivier, B.~S. Southworth, J.~S. Warsa, and H.~Park, ``Consistent second
  moment methods with scalable linear solvers for radiation transport,''
  \emph{SIAM Journal on Scientific Computing}, vol.~47, no.~3, pp.
  A1383--A1407, 2025.

\bibitem[Olivier and Haut(2023{\natexlab{a}})]{smm_olivier}
S.~Olivier and T.~S. Haut, ``High-order finite element second moment methods
  for linear transport,'' \emph{Nuclear Science and Engineering}, 2023.

\bibitem[Larsen(1988)]{LARSEN1988459}
\BIBentryALTinterwordspacing
E.~Larsen, ``A grey transport acceleration method for time-dependent radiative
  transfer problems,'' \emph{Journal of Computational Physics}, vol.~78, no.~2,
  pp. 459--480, 1988. [Online]. Available:
  \url{https://www.sciencedirect.com/science/article/pii/0021999188900605}
\BIBentrySTDinterwordspacing

\bibitem[Morel et~al.(2007)Morel, {Brian Yang}, and Warsa]{MOREL2007244}
\BIBentryALTinterwordspacing
J.~E. Morel, T.-Y. {Brian Yang}, and J.~S. Warsa, ``Linear multifrequency-grey
  acceleration recast for preconditioned krylov iterations,'' \emph{Journal of
  Computational Physics}, vol. 227, no.~1, pp. 244--263, 2007. [Online].
  Available:
  \url{https://www.sciencedirect.com/science/article/pii/S0021999107003294}
\BIBentrySTDinterwordspacing

\bibitem[Park et~al.(2012)Park, Knoll, Rauenzahn, Wollaber, and
  Densmore]{doi:10.1080/00411450.2012.671224}
\BIBentryALTinterwordspacing
H.~Park, D.~A. Knoll, R.~M. Rauenzahn, A.~B. Wollaber, and J.~D. Densmore, ``A
  consistent, moment-based, multiscale solution approach for thermal radiative
  transfer problems,'' \emph{Transport Theory and Statistical Physics},
  vol.~41, no. 3-4, pp. 284--303, 2012. [Online]. Available:
  \url{https://doi.org/10.1080/00411450.2012.671224}
\BIBentrySTDinterwordspacing

\bibitem[Olivier and Morel(2017)]{me}
S.~Olivier and J.~E. Morel, ``Variable {E}ddington factor method for the
  {$S_N$} equations with lumped discontinuous {G}alerkin spatial discretization
  coupled to a drift-diffusion acceleration equation with mixed finite-element
  discretization,'' \emph{Journal of Computational and Theoretical Transport},
  vol.~46, no. 6-7, pp. 480--496, 2017.

\bibitem[Olivier et~al.(2023)Olivier, Pazner, Haut, and Yee]{dgvef_olivier}
\BIBentryALTinterwordspacing
S.~Olivier, W.~Pazner, T.~S. Haut, and B.~C. Yee, ``A family of independent
  variable eddington factor methods with efficient preconditioned iterative
  solvers,'' \emph{Journal of Computational Physics}, vol. 473, p. 111747,
  2023. [Online]. Available:
  \url{https://www.sciencedirect.com/science/article/pii/S0021999122008105}
\BIBentrySTDinterwordspacing

\bibitem[Olivier and Haut(2023{\natexlab{b}})]{rtvef_olivier}
S.~Olivier and T.~S. Haut, ``High-order mixed finite element variable eddington
  factor methods,'' \emph{Journal of Computational and Theoretical Transport},
  vol.~52, no.~2, pp. 79--142, 2023.

\bibitem[Anistratov and Warsa(2018)]{dima_dfem}
D.~Y. Anistratov and J.~S. Warsa, ``Discontinuous finite element
  quasi-diffusion methods,'' \emph{Nuclear Science and Engineering}, vol. 191,
  no.~2, pp. 105--120, 2018.

\bibitem[Warsa and Anistratov(2018)]{two-level-independent-warsa}
J.~Warsa and D.~Anistratov, ``Two-level transport methods with independent
  discretization,'' \emph{Journal of Computational and Theoretical Transport},
  vol.~47, no. 4-6, pp. 424--450, 2018.

\bibitem[Ganapol and Pomraning(1983)]{GANAPOL1983311}
\BIBentryALTinterwordspacing
B.~Ganapol and G.~Pomraning, ``The non-equilibrium marshak wave problem: A
  transport theory solution,'' \emph{Journal of Quantitative Spectroscopy and
  Radiative Transfer}, vol.~29, no.~4, pp. 311--320, 1983. [Online]. Available:
  \url{https://www.sciencedirect.com/science/article/pii/0022407383900493}
\BIBentrySTDinterwordspacing

\bibitem[{Gentile}(2001)]{2001JCoPh.172..543G}
N.~A. {Gentile}, ``{Implicit Monte Carlo Diffusion{\textemdash}An Acceleration
  Method for Monte Carlo Time-Dependent Radiative Transfer Simulations},''
  \emph{Journal of Computational Physics}, vol. 172, no.~2, pp. 543--571, Sep.
  2001.

\bibitem[Brunner(2023{\natexlab{a}})]{osti_2280904}
\BIBentryALTinterwordspacing
T.~A. Brunner, ``A family of multi-dimensional thermal radiative transfer test
  problems,'' Lawrence Livermore National Laboratory (LLNL), Livermore, CA
  (United States), Tech. Rep., 12 2023. [Online]. Available:
  \url{https://www.osti.gov/biblio/2280904}
\BIBentrySTDinterwordspacing

\bibitem[Clark(1987)]{CLARK1987311}
\BIBentryALTinterwordspacing
B.~A. Clark, ``Computing multigroup radiation integrals using
  polylogarithm-based methods,'' \emph{Journal of Computational Physics},
  vol.~70, no.~2, pp. 311--329, 1987. [Online]. Available:
  \url{https://www.sciencedirect.com/science/article/pii/0021999187901859}
\BIBentrySTDinterwordspacing

\bibitem[Miften and Larsen(1993)]{QDBC}
M.~M. Miften and E.~W. Larsen, ``The quasi-diffusion method for solving
  transport problems in planar and spherical geometries,'' \emph{Transport
  Theory and Statistical Physics}, vol.~22, no. 2-3, pp. 165--186, 1993.

\bibitem[Adams and Larsen(2002)]{AL}
\BIBentryALTinterwordspacing
M.~L. Adams and E.~W. Larsen, ``Fast iterative methods for discrete-ordinates
  particle transport calculations,'' \emph{Progress in Nuclear Energy},
  vol.~40, no.~1, pp. 3--159, 2002. [Online]. Available:
  \url{https://www.sciencedirect.com/science/article/pii/S0149197001000233}
\BIBentrySTDinterwordspacing

\bibitem[Castillo(2006)]{CASTILLO20061307}
\BIBentryALTinterwordspacing
P.~Castillo, ``A review of the local discontinuous galerkin (ldg) method
  applied to elliptic problems,'' \emph{Applied Numerical Mathematics},
  vol.~56, no.~10, pp. 1307--1313, 2006, selected Papers from the First Chilean
  Workshop on Numerical Analysis of Partial Differential Equations (WONAPDE
  2004). [Online]. Available:
  \url{https://www.sciencedirect.com/science/article/pii/S0168927406000651}
\BIBentrySTDinterwordspacing

\bibitem[Arnold et~al.(2002)Arnold, Brezzi, Cockburn, and Marini]{Arnold2002}
D.~N. Arnold, F.~Brezzi, B.~Cockburn, and L.~D. Marini, ``Unified analysis of
  discontinuous {G}alerkin methods for elliptic problems,'' \emph{SIAM Journal
  on Numerical Analysis}, vol.~39, no.~5, pp. 1749--1779, 2002.

\bibitem[Cockburn and Dong(2007)]{10.1007/s10915-007-9130-3}
B.~Cockburn and B.~Dong, ``An analysis of the minimal dissipation local
  discontinuous {G}alerkin method for convection-diffusion problems,''
  \emph{Journal of Scientific Computing}, vol.~32, no.~2, p. 233–262, Aug.
  2007.

\bibitem[MFEM()]{mfem-web}
``{MFEM}: Modular finite element methods {[Software]},''
  \url{https://mfem.org}, 2010.

\bibitem[Andrej et~al.(2024)Andrej, Atallah, Bäcker, Camier, Copeland, Dobrev,
  Dudouit, Duswald, Keith, Kim, Kolev, Lazarov, Mittal, Pazner, Petrides,
  Shiraiwa, Stowell, and Tomov]{mfem-2024}
J.~Andrej, N.~Atallah, J.-P. Bäcker, J.-S. Camier, D.~Copeland, V.~Dobrev,
  Y.~Dudouit, T.~Duswald, B.~Keith, D.~Kim, T.~Kolev, B.~Lazarov, K.~Mittal,
  W.~Pazner, S.~Petrides, S.~Shiraiwa, M.~Stowell, and V.~Tomov,
  ``High-performance finite elements with {MFEM},'' \emph{The International
  Journal of High Performance Computing Applications}, vol.~38, no.~5, pp.
  447--467, 2024.

\bibitem[Olivier(2026)]{allium}
\BIBentryALTinterwordspacing
S.~Olivier, ``lanl/allium: v1.0,'' Mar. 2026. [Online]. Available:
  \url{https://doi.org/10.5281/zenodo.19118775}
\BIBentrySTDinterwordspacing

\bibitem[Falgout and Yang(2002)]{hypre}
R.~D. Falgout and U.~M. Yang, ``Hypre: A library of high performance
  preconditioners,'' in \emph{Proceedings of the International Conference on
  Computational Science-Part III}, ser. ICCS '02.\hskip 1em plus 0.5em minus
  0.4em\relax Berlin, Heidelberg: Springer-Verlag, 2002, p. 632–641.

\bibitem[Hamilton et~al.(2009)Hamilton, Benzi, and Warsa]{hamilton2009negative}
S.~Hamilton, M.~Benzi, and J.~Warsa, ``Negative flux fixups in discontinuous
  finite element {$S_N$} transport,'' in \emph{International Conference on
  Mathematics, Computational Methods and Reactor Physics (M\&C 2009), American
  Nuclear Society, LaGrange Park, Illinois, USA}.\hskip 1em plus 0.5em minus
  0.4em\relax Citeseer, 2009.

\bibitem[Brunner(2023{\natexlab{b}})]{brunner_github}
\BIBentryALTinterwordspacing
T.~A. Brunner, ``{AnalyticMGOpac}: An analytic opacity integrator for
  multigroup calculations,'' 2023, github repository. [Online]. Available:
  \url{https://github.com/LLNL/AnalyticMGOpac}
\BIBentrySTDinterwordspacing

\end{thebibliography}

\end{document}